%% file: main.tex
\documentclass[11pt]{article}
\usepackage{fullpage}

\usepackage[utf8]{inputenc} 
\usepackage[T1]{fontenc}    
\usepackage{hyperref}       
\usepackage{url}            
\usepackage{booktabs}       
\usepackage{amsfonts}       
\usepackage{nicefrac}       
\usepackage{microtype}      
\usepackage{tcolorbox}
\usepackage{diagbox}
\usepackage{enumerate}
\usepackage[shortlabels]{enumitem}
\usepackage{algorithm}
\usepackage[algo2e]{algorithm2e}
\usepackage{subfigure}
\usepackage{graphicx} 
\usepackage{caption}
\usepackage{amsmath}
\usepackage{amsthm}
\usepackage{amssymb}
\usepackage{tikz}
\usepackage{tablefootnote}
\usepackage{multirow}
\usepackage{enumerate}
\usepackage{color}
\usepackage{xcolor}
\usepackage{natbib}
\usepackage{overpic} 

\allowdisplaybreaks[4]

\usepackage{mathrsfs}
\usepackage{hyperref}
\usepackage{bm,bbm,todonotes}
\allowdisplaybreaks

\newtheorem{Rema}{Remark}[section]

\def\hyz#1 {\textcolor{red}{Hyz: #1 }}
\def\lx#1 {\textcolor{magenta}{Lx: #1 }}

\hypersetup{
	colorlinks=true,
	filecolor=blue,
	citecolor = blue,
	urlcolor=cyan,
}

\input{tex/math_commands}

\title{A Regularized Online Newton Method for Stochastic Convex Bandits with 
 Linear Vanishing Noise
}

\author{
Jingxin Zhan 
\thanks{School of Mathematical Sciences, Peking University; email: \texttt{bjdxzjx@pku.edu.cn}.}
\and
Yuchen Xin 
\thanks{School of Mathematical Sciences, Peking University; email: \texttt{2301110087@pku.edu.cn}. } 
\and
Kaicheng Jin
\thanks{School of Mathematical Sciences, Peking University; email: \texttt{kcjin@pku.edu.cn}. } 
\and
Zhihua Zhang 
\thanks{School of Mathematical Sciences, Peking University; email: \texttt{zhzhang@math.pku.edu.cn}. }
}

\begin{document}

\maketitle

\begin{abstract}%
We study a stochastic convex bandit problem where the subgaussian noise parameter is assumed to decrease linearly as the learner selects actions closer and closer to the minimizer of the convex loss function. Accordingly, we propose a Regularized Online Newton Method (\textbf{RONM}) for solving the problem, based on the Online Newton Method (\textbf{ONM}) of 
    \cite{fokkema2024onlinenewtonmethodbandit}. Our \textbf{RONM} reaches a polylogarithmic regret in the time horizon $n$ when the loss function grows quadratically in the constraint set, which recovers the results of \cite{lumbreras2024linear} in linear bandits. Our analyses rely on the growth rate of the precision matrix $\Sigma_t^{-1}$ in \textbf{ONM} and we find that linear growth solves the question exactly. These analyses also help us obtain 
better convergence rates 
when the loss function grows faster. We also study and analyze two new bandit models: stochastic convex bandits with 
noise scaled to a subgaussian parameter function and convex bandits with stochastic multiplicative noise.
\end{abstract}

\section{Introduction}\label{sec:intro}
\input{tex/main/intro}

\input{tex/main/notion_assump}

\section{Preliminaries}\label{Preliminaries}
\input{tex/main/assumption}
\subsection{Online Newton method for bandit convex optimization}\label{pre onm}
\input{tex/main/preliminaries_for_ONM}

\section{Main Results}\label{main result}
In this section, first, we present our algorithm that we call the Regularized Online Newton Method (\textbf{RONM}). The basic idea of our method is to introduce a regularized term built on the online Newton method for bandit convex optimization  
 (\textbf{ONM}) \citep{fokkema2024onlinenewtonmethodbandit}. And then we present our theoretical results.
\subsection{The regularized online Newton method}\label{Regularized online newton method}
\input{tex/main/algorithm}

\subsection{Theoretical results}\label{main theory}
\input{tex/main/maintrsult}

\section{Sketch of Proof for Theorem \ref{fqg}}\label{Sketch of proof for Theorem for fqg}
\input{tex/main/proof_for_theorem_fqg}

\section{Faster Convergence Rates}\label{section 2/q}
\input{tex/main/convergence_rate}

\section{Concluding Remarks}\label{concluding remarks}
\input{tex/main/open}
\newpage
\appendix

\bibliography{ref}
\bibliographystyle{plainnat}

\section{Proof for Theorem \ref{fqg} and Theorem \ref{rate >1}}\label{proof main}
\input{tex/append/Algorithm_and_proof}

\section{\texorpdfstring{The case when $f(x)$ is $(\beta,1)$ convex}{The case when f(x) is (beta,1) convex}}\label{q=1 proof}
\input{tex/append/q=1}

\section{Proof for Corollary \ref{equi f sigma} and Corollary \ref{multiplicative noise}}\label{proof sigma scaled}
\input{tex/main/noise_scaled}

\section{Some Useful Facts}\label{useful}
\input{tex/append/useful_facts}

\section{Proofs for Lemmas}\label{proof for lemmas}
\input{tex/append/proofforlemmas}

\section{
Generalized Stein's Lemma
}\label{exchange}
\input{tex/append/exchange}

\section{Computation of Hessian Matrices}\label{convex app}
\input{tex/append/hessian}

\section{Surrogate Function}\label{surro app}
\input{tex/append/Surrogatefunction}

\section{Concentration Bounds}
\input{tex/append/concentration}

\section{Auxiliary Lemma}\label{aux lem}
\input{tex/append/auxi}

\section{Constraints for Constants}\label{cons app}
\input{tex/append/constraint}

\section{Discussion of the Exponential Dependence on Dimension}\label{discussion of exponential}
\input{tex/append/exponential_dependence}

\end{document}

%% file: tex/math_commands.tex
\newcommand{\indep}{\perp\!\!\!\perp}		
\newcommand{\N}{\mathcal{N}}				
\newcommand{\1}{\mathbbm{1}}				
\renewcommand{\P}{\mathbb{P}}				
\newcommand{\E}{\mathbb{E}}					
\newcommand{\lip}{\operatorname{lip}}		
\newcommand{\Reg}{\operatorname{Reg}}		
\newcommand{\agp}[1]{\langle #1 \rangle}	
\newcommand{\sgn}[1]{\operatorname{sgn}(#1)}	

\newcommand{\normmm}[1]{\lvert\kern-0.25ex\lvert\kern-0.25ex\lvert #1 \rvert\kern-0.25ex\rvert\kern-0.25ex\rvert}

\newcommand{\tr}[1]{\operatorname{tr}(#1)}		
\newcommand{\diag}[1]{\operatorname{diag}(#1)}	
\newcommand{\R}{\mathbb{R}}		

\renewcommand*{\d}{\mathop{}\!\mathrm{d}}	
\newcommand{\pd}{\partial}					

\newtheorem{Def}{Definition}[section]
\newtheorem{Lem}[Def]{Lemma}
\newtheorem{Thm}[Def]{Theorem}
\newtheorem{Cor}[Def]{Corollary}
\newcommand{\bc}{\ell}
\newcommand{\idm}{\mathbf{I}}

%% file: tex/main/intro.tex
 Bandit convex optimization  \citep[see, e.g., ][]{lattimore40bandit} can be regarded as an online version of the zeroth-order optimization problem, a fundamental issue in optimization with many applications in operations research and other fields. 
In the stochastic convex  bandit problem, the learner picks an arm from a convex action set $\mathcal{K}\subset \R^d$. At the beginning, the environment secretly chooses a convex loss function $f(x)\colon \mathcal{K}\to [0,1]$, and in every round $t$ the learner picks $X_t$ from $\mathcal{K}$ and then suffers from a loss 
$Y_t=f(X_t)+\varepsilon_t$, where $\varepsilon_t$ is a conditionally  zero-mean subgaussian
random variable. The goal in the stochastic convex bandit  problem is to control the regret over $n$ rounds:
\[
\Reg_n=\sup _{x \in \mathcal{K}} \; \sum_{t=1}^n\left(f\left(X_t\right)-f\left(x\right)\right).
\]

Recently, \cite{lumbreras2024linear} proposed  a vanishing noise model where they 
assumed that $\varepsilon_t$ is conditionally $\sigma_t$-subgaussian with
$
\sigma_t\leq \|X_t-x_\star\|_2.
$
Here $x_\star$ is the minimizer of $f(x)$ in $\mathcal{K}$. Such a  vanishing noise model is intuitively interesting in certain contexts. For instance, in recommendation systems, it is reasonable to assume that a user's decision becomes more confident as the recommendation (action) aligns more closely with the user's preference (maximizer/minimizer). The reward/loss might even become deterministic when the recommendation perfectly matches the preference. More concretely, in quantum mechanics, measurement outcomes are random, as determined by Born's rule. However, the variance of these probabilistic outcomes decreases quadratically for projections that are aligned with the unknown pure state (\cite{Lumbreras_2022}).

 \cite{lumbreras2024linear} also proposed an algorithm based on LinUCB which achieves a
polylogarithmic regret.
However, they only studied the setting of both the linear loss and the unit sphere as the action set. 
In this paper, we  extend this setting to  general convex bandits  with a convex loss $f$ and a convex action set $\mathcal{K}$, and
would answer the following natural question:

\emph{Does and when does there exist an algorithm achieving polylogarithmic regret for stochastic convex bandits with linear vanishing noise?}
\paragraph{Contribution}
Our contributions can be summarized as follows:
\begin{enumerate}
    \item We first 
    consider the linear vanishing noise model in convex bandits
    and propose a new algorithm Regularized Online Newton Method (\textbf{RONM}) based on \textbf{ONM} (\cite{fokkema2024onlinenewtonmethodbandit}). We prove that it reaches polylogarithmic regret when $f(x)$ grows quadratically on $\mathcal{K}$.
    \item We provide a new analysis of the growth rate of the precision matrix $\Sigma_t^{-1}$ in time $t$ in the \textbf{ONM} algorithm. We find that \textbf{ONM} achieves polylogarithmic regret in this setting when $\Sigma_t^{-1}$ grows linearly. We also find that when $f(x)$ grows faster, $\Sigma_t^{-1}$ also grows faster, which yields a faster convergence rate of $\|X_t-x_\star\|_2.$
    \item We first propose and analyze stochastic convex bandits with noise scaled to $\sigma(x)$ and convex bandits with stochastic multiplicative noise. Here $\sigma(x)$ is a subgaussian parameter function. 
\end{enumerate}

%% file: tex/main/notion_assump.tex
\paragraph{Notation}
Given a function $f\colon \mathbb{R}^d \rightarrow \mathbb{R}$, with a slight abuse of notation, we write $\nabla f(x)$ for its gradient at $x$ and $\nabla^2f(x)$ for the Hessian. Given a set $U$, we let $\operatorname{lip}_U(f)=\sup \left\{\frac{f(x)-f(y)}{\|x-y\|_2}: x, y \in U, x \neq y\right\}$ and especially when $U=\R^d$ we will omit 
the subscript. The ball centered at $x$ and of radius $r$  in $\R^d$ is $\mathbb{B}_r^d(x)=\left\{y \in \mathbb{R}^d:\|y-x\|_2 \leq r\right\}$ and we omit $x$ when $x=0$. The identity matrix of dimension $d$ is denoted by $\idm_d$, and the $\ell_p$ norm denoted by $\|\cdot\|_p$. Given a square matrix $A$ we use the notation $\|x\|_A^2=x^{\top} A x$. The operator norm of  matrix $A$ is $\|A\|=\max _{x \neq 0
}\|A x\|_2 /\|x\|_2$. Given a symmetric matrix $A$, we use $\lambda_{\operatorname{min}}(A)$ and $\lambda_{\operatorname{max}}(A)$ to denote its minimum and maximum eigenvalues, respectively. Given two  matrices $A$ and $B$, we denote  $A\succeq B$ if $A-B$ is positive semidefinite. We use $\mathbb{P}$ to refer to the probability measure on some measurable space carrying all the random variables associated with the learner/environment interaction, including actions, losses, noise, and any exogenous randomness incurred by the learner. The associated expectation operator is $\mathbb{E}$ and filtration is $(\mathscr{F}_t)_{t=1}^n$. We denote the conditional expectation with respect to $\mathscr{F}_t$ as $\E_t[\cdot]=\E[\cdot|\mathscr{F}_t]$. We denote by $\mathcal{N}(\mu, \Sigma)$ a Gaussian vector with mean $\mu$ and covariance $\Sigma$. For a zero-mean random variable $W$, we say it to be $\sigma$-subgaussian if for all $\lambda\in\R$, $\E[\exp(\lambda W)]\leq\exp(\lambda^2\sigma^2/2)$. And let $\|W\|_{\psi_2}:=\inf \left\{t>0: \mathbb{E}\left[\exp \left(W^2 / t^2\right)\right] \leq 2\right\}$ be the subgaussian norm. 
\paragraph{Assumptions}
In this paper, we consider a stochastic convex bandit on a convex and compact action set $\mathcal{K}\subset\R^d$, which should have a non-empty interior. We also assume that $\mathcal{K}$ contains $\mathbb{B}_r^d$ and for simplicity, in our main results, we assume that $R=\sup_{x\in\mathcal{K}}\|x\|_2=1$. The environment chooses a convex loss function $f(x)\colon \mathcal{K}\to[0,1]$, which is also assumed to be $G$-Lipschitz continuous on $\mathcal{K}$, i.e., 
$\lip_{\mathcal{K}}(f)\leq G$. 
The learner interacts with the environment over $n$ rounds. At time  $t$ the learner picks an arm $X_t$ from $\mathcal{K}$ and then bears and only observes a loss $Y_t=f(X_t)+\varepsilon_t,$ where $\varepsilon_t$ is a zero-mean $\sigma_t$-subgaussian
random variable, conditioning on $\mathscr{F}_{t-1}$ and $X_t$. Here $\mathscr{F}_{t-1}
$, generated by $X_1,Y_1, \ldots, X_{t-1}, Y_{t-1}$, is the filtration of all information up to time $t$ before $X_t$ is chosen.

\paragraph{Related works}
Initialed by \cite{FK05} and \cite{Kle04}, the literature on convex bandits is increasingly extensive. In the stochastic setting, \cite{AFHK13} first obtained $\operatorname{poly}(d) \sqrt{n}$ regret using classical zeroth-order optimization techniques. Their results were improved by \cite{LG21a}, which yielded a better dependence on dimension. In the adversarial setting, $\sqrt{n}$ regret was also achieved by \cite{HL14} with an online Newton method  when the loss function is both smooth and strongly convex. And if only boundedness of the loss function holds, \cite{BDKP15}  showed that $\sqrt{n}$ regret is also possible when the dimension $d=1$. The case for general $d$ was later proved by \cite{BE18}. Their non-constructive information-theoretic approaches 
were followed by \cite{Lat20-cvx}. The current state of art is \cite{fokkema2024onlinenewtonmethodbandit}, for which the regret is $d^{3.5}\sqrt{n}$ in the adversarial setting while $d^{1.5}\sqrt{n}$ (informally) in the stochastic setting using the online Newton method. There are also many other variants of the online Newton method in bandit convex optimization \citep[see, e.g., ][]{suggala2024second,mhammedi2024online}.

\paragraph{Organization} The remainder of this paper is organized as follows. Section \ref{Preliminaries} gives preliminaries. Section \ref{main result} presents our main results, including the algorithm and theoretical results. Section \ref{Sketch of proof for Theorem for fqg} provides a proof sketch of our main theorem. Section \ref{section 2/q} analyzes the convergence rate. We conclude the remaining open problems in Section \ref{concluding remarks}. All the proof details are deferred to the appendices. 

%% file: tex/main/assumption.tex
In this section, we present some necessary definitions and
then give a brief review of the online Newton method in bandit convex optimization \citep{fokkema2024onlinenewtonmethodbandit}, which will be abbreviated as \textbf{ONM}.
For more details, the readers are also recommended to read the  literature \citep{lattimore40bandit}. 

\subsection{Vanishing noise}

In this paper, we consider a noise model such that the subgaussian parameter $\sigma_t$ satisfies 
$
\sigma_t\leq \|X_t-x_\star\|_2.
$
Here $x_\star$ is the minimizer of $f(x)$ in $\mathcal{K}$, which is unique due to our assumption about $f(x)$ later. We call this bandit model \emph{a stochastic convex bandit with linear vanishing noise}.
\paragraph{Noise scaled to $\sigma(x)$}\label{scaled noise}
We also investigate a stronger assumption, that is, there exists a function $\sigma(x)\colon \mathcal{K}\to\R^+$ such that $\varepsilon_t=\sigma(X_t)\cdot\bar{\varepsilon}_t.$ We then assume that $\{\bar{\varepsilon}_t\}_{t=1}^n$ are independent and identically distributed $1$-subgaussian 
non-degenerate random variables. Hence it is clear that $\sigma_t\leq\sigma(X_t)$ and we call this bandit model \emph{a stochastic convex bandit with noise scaled to $\sigma(x)$}. Note that when $\sigma(x)\leq\|x-x_\star\|_2$ , it is also a stochastic convex bandit with linear vanishing noise.
Especially, when $\sigma(x)=f(x)$ and $f(x_\star)=0$, the feedback at time step $t$ is
\[
Y_t=f(X_t) +\sigma(X_t)\cdot\bar{\varepsilon}_t=f(X_t)(1+\bar{\varepsilon}_t).
\]
Hence, such a model becomes a special case of \emph{a convex bandit with stochastic multiplicative noise}, 
which is also an important theme 
because multiplicative noise models are underappreciated and appear in many physical problems \citep[see, e.g., ][]{hodgkinson2020multiplicativenoiseheavytails}.
\subsection{The quadratic growth condition}
In this paper, we assume $f(x)$ grows fast enough for the player to reach polylogarithmic regret and hence we define the quadratic growth condition. One should note that it is a \emph{local} property.
\begin{Def}[$\rho$-Quadratic Growth (QG)]\label{qg def}
    We say that $f(x)$ has a $\rho$-quadratic growth property on $\mathcal{K}$ if $f(x)$ is convex, has a unique minimizer $x_\star$ in $\mathcal{K}$, and there exists a constant $\rho>0$ such that for all $x\in\mathcal{K}$,
\[
f(x)-f(x_\star)\ge \frac{\rho}{2}\|x-x_\star\|_2^2.
\]
\end{Def}
Quadratic growth is a common regular condition in optimization theory and has close connections with other conditions \citep[see, e.g.,][]{karimi2020linearconvergencegradientproximalgradient, drusvyatskiy2018error}.
\paragraph{Example}\label{example}
Let $\mathcal{K}=\mathbb{B}_1^d$. Then for any $\theta \in\mathbb{B}_1^d$, consider 
$
f(x)=\agp{x,\theta},\forall x\in\mathbb{B}_1^d,
$
where we ignore the assumption that $f(x)\in[0,1]$ for simplicity. Clearly, now $x_\star=-\theta$ and then for all $x\in\mathbb{B}_1^d$,
\[
f(x)-f(x_\star)=\agp{x+\theta,\theta}=1+\agp{x,\theta}\ge\frac{1}{2}\|x-x_\star\|_2^2,
\]
where we used that $\|x\|_2\leq 1.$ This implies that $f(x)$ is $1$-QG on $\mathcal{K}$. Hence, the linear loss function in \cite{lumbreras2024linear} also has the QG property. 
In addition, there are many other examples, and in the following, we are mostly concerned with a definition of $(\beta,\bc)$-convexity.
\begin{Def}[$(\beta,\bc)$-Convexity]\label{beta q convex}
     We say that $f(x)$ is $(\beta,\bc)$-convex on $\mathcal{K}$ if there exist constants $\beta>0$ and $\bc\ge1$ such that $f(x)-\beta\|x-x_{\star}\|_2^\bc$ is convex on $\mathcal{K}$.
\end{Def}
By Lemma \ref{q qg}, when $1<\bc\leq 2$, it is easy to see that $f(x)$ grows faster near $x_\star$ and also grows quadratically. Especially, when $\bc=2$ , clearly $f(x)$ is just $2\beta$-strongly convex.

%% file: tex/main/preliminaries_for_ONM.tex
Here we review the sketch of \textbf{ONM} and its proof. Recall that the online Newton method is a second-order method and needs loss functions' gradients and Hessian matrices for updating. However, bandit convex optimization is zeroth-order. To solve this, we need a surrogate loss function by Gaussian convolution, with which, if the player samples a new arm from a properly chosen normal distribution, the player can immediately obtain unbiased estimators for the first- and second-order information of the surrogate loss using the feedback. Finally, since the player picks points outside $\mathcal{K}$ with high probability, we apply the convex extension in \cite{fokkema2024onlinenewtonmethodbandit}.
\paragraph{Convex extension}\label{extension section}
Suppose that $\mathbb{B}_r^d \subset \mathcal{K} \subset \mathbb{B}_R^d$, and let $\pi(x)=\inf \{t>0: x \in t \mathcal{K}\}$ be the Minkowski functional of $\mathcal{K}$. Then the convex extension (denoted $e$) of $f$ is defined as
\begin{equation}\label{extension defi}
    e(x) :=\pi^+(x) f\left(\frac{x}{\pi^+(x)}\right)+G R(\pi^+(x)-1),
\end{equation}
where $\pi^+(x)=\max (1, \pi(x)).$
Such an extension makes it possible to obtain the value of $e(x)$ using only a single evaluation for $f(x)$ when $x\notin\mathcal{K}$. Its properties are summarized in Lemma \ref{extension}, which states that $e(x)$ is convex on $\R^d$ and equals $f(x)$ on $\mathcal{K}$.

When $X$ is chosen by the learner, the learner actually picks $\frac{X}{\pi^+(X)}$, and the bandit outputs $f(\frac{X}{\pi^+(X)})+\varepsilon$ as feedback. By simply substituting $f(\frac{X}{ \pi^+(X)})$ with $f(\frac{X}{ \pi^+(X)})+\varepsilon$ in 
Eq.~\eqref{extension defi}, we can feed the player with the loss
\begin{equation}\label{real loss}
        Y=e(X)+\pi^+(X)\varepsilon.
\end{equation}
Then from the learner's perspective, the loss function is $e(x)$ and noise is $\xi:=\pi^+(X)\varepsilon.$ 

\paragraph{Surrogate loss}\label{gtht defi}
Given a Gaussian distribution $\N(\mu, \Sigma)$ and a parameter $\lambda\in(0,1)$, the Gaussian optimistic smoothing surrogate function of 
$e(x)$ is defined as
\[
s(x)=\mathbb{E}\left[\left(1-\frac{1}{\lambda}\right) e(X)+\frac{1}{\lambda} e\left((1-\lambda) X+\lambda x \right)\right],
\]
where $X\sim\N(\mu,\Sigma)$. The surrogate has been widely used in bandit convex optimization \citep{bubeck2021kernel,LG21a,lattimore2023secondordermethodstochasticbandit,fokkema2024onlinenewtonmethodbandit}.

The algorithm keeps track of an iterate mean vector $\mu_t$ and covariance matrix $\Sigma_t$. Hence, we denote by $s_t$ the surrogate loss at round $t$ with $\N(\mu_t,\Sigma_t)$. At time step $t$, the learner samples a new arm $X_t$ from $\N(\mu_t,\Sigma_t)$
 and then the learner can construct estimators for $\nabla s_t(\mu_t)$ and $\nabla^2 s_t(\mu_t)$, denoted $g_t$ and $H_t$. We will give their explicit expressions  in Eq.~\eqref{Ht}. Their properties can be found in \cite{lattimore40bandit} and we include them in Appendix \ref{surro app} for our use.  

\paragraph{Online Newton method}
Generally, the online Newton method \citep{hazan2007logarithmic} is an online algorithm for a sequence of quadratic loss functions. In the bandit convex optimization problem, one typically considers a proper quadratic approximation to $s_t$. Let $q_t$ be such an approximation and $\hat{q}_t$ be the estimator of $q_t$. That is, they are defined as 
\[
q_t(x)=\left\langle \nabla s_t(\mu_t), x-\mu_t\right\rangle+\frac{1}{4}\|x-\mu_t\|_{\nabla^2 s_t(\mu_t)}^2, \quad\hat{q}_t(x)=\left\langle g_t, x-\mu_t\right\rangle+\frac{1}{4}\|x-\mu_t\|_{H_t}^2,
\]
where $g_t$ and $H_t$ are the estimators of $\nabla s_t(\mu_t)$ and $\nabla^2 s_t(\mu_t)$, respectively. Accordingly,  one implements the online Newton method with the $\left(\hat{q}_t\right)_{t=1}^n$  on $\mathcal{K}$. At every round $t$, the online Newton method plays $\mu_t$ and updates with that
\begin{equation}\label{update rule}
    \begin{aligned}
& \Sigma_{t+1}^{-1}=\Sigma_t^{-1}+\eta \nabla^2\hat{q}_t(\mu_t),
& \mu_{t+1}=\arg \min _{x \in \mathcal{K}}\left\|x-\left[\mu_t-\eta \Sigma_{t+1} \nabla\hat{q}_t(\mu_t)\right]\right\|_{\Sigma_{t+1}^{-1}}^2,
\end{aligned}
\end{equation}
where 
$\mu_1\in\mathcal{K}$ and $\Sigma_{1}=\sigma^2\idm_d.$ 
\begin{Rema}
   In summary, in \textbf{ONM}, at time step $t$, the learner samples $X_t$ from $\N(\mu_t,\Sigma_t)$ and gets a loss by Eq.~\eqref{real loss}. Using this feedback, the learner computes $g_t$ and $H_t$, with which the learner updates $\mu_t$ and $\Sigma_t$ by Eq.~\eqref{update rule}. 
\end{Rema}

\paragraph{Sketch of the regret analysis}\label{overview} 
For ease of exposition, here we give a sketch of the regret analysis, which relies on a stopping time (denoted $\tau$) that 
promises that for all $t\leq\tau$, $\frac{1}{2}\|\mu_t-x_\star\|_{\Sigma_t^{-1}}^2\leq\frac{1}{2\lambda^2L^2}$. It would be shown that $\tau=n$ holds with high probability. By Lemma \ref{regret for onm}, we have the following decomposition
\begin{equation}\label{decomposition equation}
    \begin{aligned}
          \frac{1}{2}\left\|\mu_{\tau+1}-x_\star\right\|_{\Sigma_{\tau+1}^{-1}}^2 &\leq \frac{1}{2}\left\|\mu_1-x_\star\right\|_{\Sigma_1^{-1}}^2+\frac{\eta^2}{2} \sum_{t=1}^\tau\left\|g_t\right\|_{\Sigma_{t+1}}^2-\eta \operatorname{Reg}_\tau(x_\star)\\
          & \quad +\eta\underbrace{\left(\operatorname{qReg}_\tau(x_\star)-\widehat{\operatorname{qReg}}_\tau(x_\star)\right)}_{\text{Estimation Error}}+\eta\underbrace{\left(\operatorname{Reg}_\tau(x_\star)-\operatorname{qReg}_\tau(x_\star)\right)}_{\text{Approximation Error}},
    \end{aligned}
\end{equation}
where $
 \operatorname{q Reg}_\tau(x_\star):=\sum_{t=1}^\tau\left(q_t\left(\mu_t\right)-q_t(x_\star)\right)
 $ and $\widehat{\operatorname{q Reg}}_\tau(x_\star):=\sum_{t=1}^\tau\left(\hat{q}_t\left(\mu_t\right)-\hat{q}_t(x_\star)\right)$. Note that $\operatorname{Reg}_\tau(x_\star)\ge 0$ and one can show that the sum of other terms will be less than $\frac{1}{2\lambda^2L^2}$ with high probability for carefully chosen constants. 
Then this simultaneously shows that $\tau=n$ and $\operatorname{Reg}_\tau(x_\star)\leq\frac{1}{2\lambda^2L^2\eta}$ with high probability.

%% file: tex/main/algorithm.tex
We present \textbf{RONM} in Algorithm~\ref{algorithm 1}. For the $\rho$-QG $f(x)$, we set the constants as follows:
\begin{equation}\label{alg1 1}
    \begin{aligned}
    \sigma=\frac{r}{5\sqrt2 d}
    ,\quad\lambda=\frac{1}{HdL^3},\quad\eta=\frac{\gamma}{100H^2d^4L^{5}},\quad\gamma=\rho.
\end{aligned}
\end{equation}
Especially, if $f(x)$ is $(\beta,\bc)$-convex on $\mathcal{K}$, $1<\ell\leq 2$, we set $\gamma=2^{\bc-1}\beta$.

\begin{algorithm}[H]
    \caption{\textbf{RONM} for stochastic convex bandits with linear vanishing noise} \label{algorithm 1}
    Require: $\eta,\lambda,\sigma,\gamma>0$

    Set $\mu_1=0$, $\Sigma_1=\sigma^2\idm_d$ and $Y_0=0$

    \For{$t=1,2,\cdots,n$}{
        \vspace{1mm}
        ~sample $X_t$ from $\N(\mu_t,\Sigma_t)$ with density $p_t$
        \vspace{1mm}
        
        ~observe $Y_t=\pi^+(X_t) \left[f\left(\frac{X_t}{\pi^+(X_t)}\right)+\varepsilon_t\right]+G R(\pi^+(X_t)-1)$
        \vspace{1mm}

        ~let~ $R_t=\frac{p_t\left(\frac{X_t-\lambda\mu_t}{1-\lambda}\right)}{(1-\lambda)^dp_t(X_t)}$ and $Z_t=Y_t-Y_{t-1}$
        \vspace{1mm}

        ~compute $g_t=\frac{R_tZ_t\Sigma_t^{-1}(X_t-\mu_t)}{(1-\lambda)^2}$
        \vspace{1mm}

        ~compute $H_t=\frac{\lambda R_tZ_t}{(1-\lambda)^2}\left[\frac{\Sigma_t^{-1}(X_t-\mu_t)(X_t-\mu_t)^\top\Sigma_t^{-1}}{(1-\lambda)^2}-\Sigma_t^{-1}\right]$
        \vspace{1mm}

        ~$\Sigma_{t+1}^{-1}\gets\Sigma_t^{-1}+ \eta
       \left( \frac{1}{2}H_t+\gamma\idm_d\right)$
       \vspace{1mm}

        ~$\mu_{t+1}\gets\arg\min_{\mu\in \mathcal{K}}\|\mu-[\mu_t-\eta\Sigma_{t+1}g_t]\|_{\Sigma_{t+1}^{-1}}$
       
    }
\end{algorithm}

We have made two important modifications to the original \textbf{ONM}. First, the estimators $g_t$ and $H_t$ for $\nabla s_t(\mu_t)$ and $\nabla^2 s_t(\mu_t)$ are given as 
    \begin{equation}\label{Ht}
    \begin{aligned}
g_t & =\frac{Z_t R_t(\mu_t)}{1-\lambda} \Sigma_t^{-1}\left[\frac{X_t-\lambda \mu_t}{1-\lambda}-\mu_t\right], \\
H_t & =\frac{\lambda Z_t R_t(\mu_t)}{(1-\lambda)^2}\left(\Sigma_t^{-1}\left[\frac{X_t-\lambda \mu_t}{1-\lambda}-\mu_t\right]\left[\frac{X_t-\lambda \mu_t}{1-\lambda}-\mu_t\right]^{\top} \Sigma_t^{-1}-\Sigma_t^{-1}\right),
\end{aligned}
\end{equation}
where $R_t(z)=\frac{p_t\left(\frac{X_t-\lambda z}{1-\lambda}\right)}{(1-\lambda)^d p_t\left(X_t\right)}$ and $p_t$ is the density of $\N(\mu_t,\Sigma_t)$. It is different from \cite{fokkema2024onlinenewtonmethodbandit} because we replace $Y_t$ with $Z_t:=Y_{t}-Y_{t-1}$ and $ Y_0=0$ for technical reasons, which is also used in \cite{lattimore2023secondordermethodstochasticbandit}. By Lemma \ref{unbias}, they are both unbiased. 

Second, we introduce a regularized term $\frac\gamma2\|x-\mu_t\|_2^2$ to the sequence of quadratic loss functions 
$\left(\hat{q}_t\right)_{t=1}^n$ in the online Newton method. That is,
\[
\hat{q}^\gamma_t(x)=\hat{q}_t(x)+\frac{\gamma}{2}\|x-\mu_t\|_2^2, \quad \hat{q}_t(x)=\left\langle g_t, x-\mu_t\right\rangle+\frac{1}{4}\|x-\mu_t\|_{H_t}^2.
\]
Then we implement the method with $\left(\hat{q}_t^\gamma\right)_{t=1}^n$, and by Eq.~\eqref{update rule}, 
$\Sigma_{t+1}^{-1} =\Sigma_t^{-1}+ \eta   \left( \frac{1}{2}H_t+\gamma\idm_d\right)$.

\paragraph{Intuition for regularization}
Our proof for polylogarithmic regret in stochastic bandits with linear vanishing noise relies on the linear growth of $\Sigma_t^{-1}$, which is very similar to \cite{lumbreras2024linear}, where they also show that the design matrix grows linearly. However, there is no guarantee for the growth rate of $\Sigma_t^{-1}$ in \textbf{ONM}. Hence, regularization is a very natural idea for this goal because it's clear that now in every time step $t$, $\Sigma^{-1}_t$ is added by $\eta\gamma\idm_d$ and thus $\Sigma^{-1}_t\succeq\eta\gamma(t-1)\idm_d=\Omega(t)\idm_d$, which grows linearly.

%% file: tex/main/maintrsult.tex
We present theoretical results of our \textbf{RONM} method. In the setting of linear vanishing noise, we show that when $f(x)$ has the $\rho$-QG property, \textbf{RONM} can achieve a polylogarithmic regret in the horizon $n$ and we also analyze the convergence rate of $\left\|\frac{X_t}{ \pi^+(X_t)}-x_\star\right\|_2$. Here and later, $C$ and $C'$ are sufficiently large universal constants
and due to the extension, $\Reg_n= \sum_{t=1}^n\left(f\left(\frac{X_t}{ \pi^+(X_t)}\right)-f\left(x_{\star}\right)\right)$.
\begin{Thm}\label{fqg}
    If $f(x)$ has the $\rho$-QG property on $\mathcal{K}$, then with probability at least $1-\delta$, the regret of Algorithm \ref{algorithm 1} is bounded by 
    \[
    \Reg_n=\mathcal{O}(H^4d^6L^{10}/\rho),
    \]
    where $L=C[1+\log \max (n, d, H, 1/\rho, 1/\delta)], \,\delta=\operatorname{Poly}(1/n,1/d,1/H)\in(0,1)
    $ and $
    H=C'\max(G/r,1/r)
    $.  Moreover, we have that for all $t\leq n$, 
    $
    \|\frac{X_t}{ \pi^+(X_t)}-x_\star\|_2=\widetilde{\mathcal{O}}\left(t^{-\frac12}\right).
    $
\end{Thm}
\noindent Theorem \ref{fqg} recovers the polylogarithmic regret in \cite{lumbreras2024linear} by recalling Section \ref{example}. Note that their results highly rely on the linear setting and can't be applied to general convex loss functions. We also have new results about faster convergence rates.
\begin{Thm}\label{rate >1}
    If $f(x)$ is $(\beta,\bc)$-convex on $\mathcal{K}$ for $1<\bc\leq 2$, then with probability at least $1-\delta$, the regret of Algorithm \ref{algorithm 1} is bounded by 
    \[
    \Reg_n=\mathcal{O}(H^4d^6L^{10}/\beta),
    \]
    where $L=C[1+\log \max (n, d, H, 1/\beta, 1/\delta)], \,\delta=\operatorname{Poly}(1/n,1/d,1/H)\in(0,1)
    $ and $H=C'\max(G/r,1/r)$. Moreover, we have that for all $t\leq n$,
    \[
    \left\|\frac{X_t}{ \pi^+(X_t)}-x_\star\right\|_2=\widetilde{\mathcal{O}}\left(\min\bigg(t^{-\frac12}, \Big(\frac{\sqrt{2}}{r}\Big)^{\frac{(d-1)}{\bc}}t^{-\frac{1}{\bc}}\bigg)\right).
    \]
\end{Thm}
\noindent
In the convergence rates in Theorems~\ref{fqg} and  \ref{rate >1} we hide polylogarithmic terms,  
polynomial terms of $d$ and dependence on other parameters for simplicity. One should note that the convergence rate of $t^{-\frac{1}{\bc}}$ in Theorem \ref{rate >1} is faster than $t^{-\frac12}$ in Theorem \ref{fqg}, which makes sense because a $(\beta,\bc)$-convex $f(x)$ grows faster than a quadratic function for $1<\bc<2$. However, the price is an extra exponential dependence on the dimension $d$, which can be removed when $\bc=2$ and we discuss this in Appendix \ref{discussion of exponential}. Their proofs can be found in Appendix \ref{proof main}.

In the setting of  noise scaled to $\sigma(x)$, we find that $f(x)$ and $\sigma(x)$ play very similar roles in the sense of the following corollary:
\begin{Cor}\label{equi f sigma}
    In the setting of  noise scaled to $\sigma(x)$, if for all $x\in\mathcal{K}$, $\sigma(x)\leq\|x-x_\star\|_2$, then when at least one of $\sigma(x)$ and $f(x)$ has the $\rho$-QG property (resp. is $(\beta,\bc)$-convex, $1<\bc<2$) on $\mathcal{K}$, the convergence rate in Theorem \ref{fqg} (resp. Theorem \ref{rate >1}) also holds. Especially, if there exists $C>0$ such that $\sigma(x)\ge C(f(x)-f(x_\star))$ for all $x\in\mathcal{K}$, then the regret bound also holds.
\end{Cor}
\begin{Rema}
    Note that when $\sigma(x)$ has the QG property, there are \emph{not} any restrictions for $f(x)$ (e.g., not necessarily convex).
\end{Rema}
\noindent
Then, since convex bandits with stochastic multiplicative noise are special cases, we have
\begin{Cor}\label{multiplicative noise}
    In the setting of convex bandits with stochastic multiplicative noise, if $f(x)$ has the $\rho$-QG property (resp. is $(\beta,\bc)$-convex, $1<\bc<2$) on $\mathcal{K}$ and $f(x_\star)=0$, then the results in Theorem \ref{fqg} (resp. Theorem \ref{rate >1}) also hold. 
\end{Cor}
\noindent
Here we say ``also hold'' meaning that orders of regrets and convergence rates in $n$ and $d$ remain the same while omitting the dependence on other parameters. Their proofs can be found in Appendix \ref{proof sigma scaled}.

\subsection{\texorpdfstring{Results with the case $\bc=1$}{The case when q=1}}\label{main result q=1}

Intuitively, when $\bc=1$, $f(x)$ grows fastest and the algorithm also seems to have the fastest convergence rate. Nevertheless, this case is much harder and we will explain the technical challenges in Section \ref{section 2/q}. Technically, to achieve that we need some other assumptions:
\begin{description}
    \item[Assumption 1:] $f(x)$ is $(\beta,1)$-convex on $\R^d$ with its minimizer $x_\star$ in $\mathbb{B}_1^d$, i.e., $\mathcal{K}=\mathbb{B}_1^d$, and $G$-Lipschitz continuous in $\R^d$ and $\sup_{x\in\mathcal{K}}|f(x)|\leq 1;$
    \item[Assumption 2:] Any queries outside $\mathcal{K}$ are also allowed. In other words, every time the player picks $X_t$ from $\R^d$, the player will get feedback $Y_t=f(X_t)+\varepsilon_t$, where $\varepsilon_t$ is conditionally $\|X_t-x_\star\|_{2}$-subgaussian;
    \item[Assumption 3:] $f(x)$ is symmetric with respect to $x_\star$, i.e., for all $x\in\R^d$, $f(x_\star+x)=f(x_\star-x).$
\end{description}
Assumption 2 is also used in \cite{lattimore2023secondordermethodstochasticbandit} and can be regarded as unconstrained stochastic convex bandits. The remaining two assumptions are for technical reasons. Note that any function of $\|x-x_\star\|_2$ meets Assumption 3. We show that under these assumptions, the convergence rate for \textbf{ONM} can be arbitrarily close to $t^{-1}$. 
\begin{Thm}\label{q=1 result}
   If Assumptions 1-3 are satisfied, 
   then for all $\kappa\in(0,1]$, with probability at least $1-\delta$, in Algorithm \ref{algorithm 2}, for all $t\leq n$,
    \[
    \|X_t-x_\star\|_2=\widetilde{\mathcal{O}}(6^{\frac{d(2-\kappa)}{2\kappa}}t^{-1+\frac{\kappa}{2}}),
    \]
     where $
     \,\delta=\operatorname{Poly}(1/n,1/d,1/H)\in(0,1)
     $ and $H=C'\max(G,1)$.
\end{Thm}
The algorithm and proof are given in Appendix \ref{q=1 proof}. Note that we also omit the secondary terms in the convergence rate, and this result does not have a polylogarithmic regret bound as the price of a fast convergence rate (If only the polylogarithmic regret is needed, one can simply apply \textbf{RONM} with a slower convergence rate). Note that the case $\bc=1$ can be similarly modified as Corollary \ref{equi f sigma}, \ref{multiplicative noise} and we omit it for simplicity.

\begin{Rema}
    By analyses in Section \ref{section 2/q} and Lemma \ref{nabla2}, Theorem \ref{rate >1} and Theorem \ref{q=1 result} are still true if the $\|\cdot\|_2$ norm in Definition \ref{beta q convex} is replaced by the $\|\cdot\|_p$ norm for $1<p<2$.
\end{Rema}

%% file: tex/main/proof_for_theorem_fqg.tex
In this section, we present the sketch of the proof for Theorem \ref{fqg}. We follow the proof in Section \ref{overview} with some significant adaptations. The details can be found in Appendix \ref{proof main}.

\paragraph{Extension}
The convex extension keeps the property of linear vanishing in noise. Recalling Section \ref{extension section}, in the player's perspective, a noise $\xi=\pi^{+}(X)\varepsilon$ is fed when $X$ is chosen. Then by Lemma \ref{xi lip}, $\xi$ is conditionally $\sigma'$-subgaussian, where
\begin{equation}\label{extension noise}
    \sigma'\leq\pi^+(X)\left\|\frac{X}{\pi^+(X)}-x_\star\right\|_2 \leq (1+R/r)\|X-x_\star\|_2,
\end{equation}
because the real choice is $\frac{X}{\pi^+(X)}.$

\paragraph{Bound for error terms}
The main difference takes place in Estimation Error. We improve it from $\widetilde{O}(\sqrt{n})$ to $\widetilde{O}(1)$ in stochastic convex bandits with linear vanishing noise, which just yields polylogarithmic regret. The key is that we will show that, informally, the Estimation Error in Eq.~\eqref{decomposition equation} is $\widetilde{\mathcal{O}}(\sqrt{\sum_{t=1}^n Z_t^2})$ 
in Lemma \ref{sk e3} and we note that if $\Sigma_t^{-1}$ grows linearly, as in \textbf{RONM}, then $\sqrt{\sum_{t=1}^n Z_t^2}=\widetilde{O}(1)$, which just follows from that
\[
 |Z_t|\leq |e(X_t)-e(x_\star)|+|e(X_{t-1})-e(x_\star)|+|\xi_t|+|\xi_{t-1}|=\tilde{\mathcal{O}}(t^{-1/2})\quad  \text{w.h.p.}
 \]
Here we used that $\xi_t$ is conditionally $2R/r\|X_t-x_\star\|_2$-subgaussian by Eq.~\eqref{extension noise} and
\begin{equation}\label{rate eg}
    \|X_t-x_\star\|_2^2\lesssim\|X_t-x_\star\|_{\Sigma_t^{-1}}^2/t\sim\|\mu_t-x_\star\|_{\Sigma_t^{-1}}^2/t\leq\frac{1}{2\lambda^2L^2 t}=\widetilde{\mathcal{O}}(t^{-1}),
\end{equation}
for properly chosen $\lambda$, because $X_t\sim\N(\mu_t,\Sigma_t)$. Note that this is not true for $\sum_{t=1}^n Y_t^2$ when $f(x_\star)>0$, which explains why we use $Z_t$ in the estimators $g_t$ and $H_t$ in Eq.~\eqref{Ht} instead of $Y_t$. There are also some small changes in Approximation Error and the details are deferred in Appendix \ref{proof main}.

\paragraph{New decomposition}
Due to the regularized term in \textbf{RONM}, the decomposition~\eqref{decomposition equation} 
now becomes
\begin{equation}\label{equation decomposition 1}
    \begin{aligned}
    \frac{1}{2}\left\|\mu_{\tau+1} {-} x_\star\right\|_{\Sigma_{\tau+1}^{-1}}^2 
    \leq& \frac{1}{2}\left\|\mu_1 {-} x_\star\right\|_{\Sigma_1^{-1}}^2+\frac{\eta^2}{2} \sum_{t=1}^\tau\left\|g_t\right\|_{\Sigma_{t+1}}^2-\eta \widehat{\operatorname{q^\gamma Reg}}_\tau(x_\star)\\
    =& \frac{1}{2}\left\|\mu_1 {-} x_\star\right\|_{\Sigma_1^{-1}}^2+\frac{\eta^2}{2} \sum_{t=1}^\tau\left\|g_t\right\|_{\Sigma_{t+1}}^2+\frac{\eta\gamma}{2}\sum_{t=1}^\tau\|\mu_t {-} x_\star\|_2^2-\eta\widehat{\operatorname{q Reg}}_\tau(x_\star)
    \\
    \leq&\frac{1}{2}\left\|\mu_1 {-} x_\star\right\|_{\Sigma_1^{-1}}^2+\frac{\eta^2}{2} \sum_{t=1}^\tau\left\|g_t\right\|_{\Sigma_{t+1}}^2+\frac{\eta\gamma}{2}\sum_{t=1}^\tau\|\mu_t {-} x_\star\|_2^2-\eta \operatorname{Reg}_\tau(x_\star)+\operatorname{Error},
\end{aligned}
\end{equation}
where in the first inequality $\widehat{\operatorname{q^\gamma Reg}}_\tau(x_\star):=\sum_{t=1}^\tau\left(\hat{q}_t^\gamma\left(\mu_t\right)-\hat{q}_t^\gamma(x_\star)\right)=\widehat{\operatorname{q Reg}}_\tau(x_\star)-\frac{\gamma}{2}\sum_{t=1}^\tau\|\mu_t {-} x_\star\|_2^2$ and 
in the second inequality the $\operatorname{Error}$ term is just the Approximation Error and Estimation Error defined in Eq.~\eqref{decomposition equation}.  
However, we are not able to control the right hand to be less than $\frac{1}{2\lambda^2L^2}$ because we can only bound $\eta\gamma\|\mu_t-x_\star\|_2^2/2$ by $\frac{1}{2\lambda^2L^2 t}$, whose sum in $t$ is in the order of $\frac{1}{\lambda^2 L}$ and is larger than $\frac{1}{2\lambda^2 L^2}$. To see this, recall that
in \textbf{RONM}, $\Sigma_t^{-1}\sim\eta\gamma t\idm_d$ and then $\eta\gamma\|\mu_t-x_\star\|_2^2/2\leq\|\mu_t-x_\star\|_{\Sigma_t^{-1}}^2/2t\leq \frac{1}{2\lambda^2L^2 t}$.
To fix this, we should use a different way of decomposition. Similar to Eq.~\eqref{equation decomposition 1}, we have
\begin{equation}\label{equation decomposition 2}
    \begin{aligned}
    \frac{1}{2}\left\|\mu_{\tau+1} {-} x_\star\right\|_{\Sigma_{\tau+1}^{-1}}^2 \leq&\frac{1}{2}\left\|\mu_1 {-} x_\star\right\|_{\Sigma_1^{-1}}^2+\frac{\eta^2}{2} \sum_{t=1}^\tau\left\|g_t\right\|_{\Sigma_{t+1}}^2+\frac{\eta\gamma}{2}\sum_{t=1}^\tau\|\mu_t {-} x_\star\|_2^2-\eta \widetilde{\operatorname{e Reg}}_\tau(x_\star)+\operatorname{Error},
\end{aligned}
\end{equation}
where $\widetilde{\operatorname{e Reg}}_\tau(x_\star):=\sum_{t=1}^\tau(e(\mu_t)-e(x_\star))$ and then the Approximation Error becomes $\widetilde{\operatorname{e Reg}}_\tau(x_\star)-\operatorname{qReg}_\tau(x_\star)$ and the Estimation Error remains the same. Note that, when $f(x)$ has the $\rho$-QG property,
\begin{equation}\label{ereg lower}
    \widetilde{\operatorname{e Reg}}_\tau(x_\star)=\sum_{t=1}^\tau(f(\mu_t)-f(x_\star))\ge\frac{\rho}{2}\sum_{t=1}^\tau\|\mu_t-x_\star\|_2^2,
\end{equation} 
hence, choosing $\gamma=\rho$ and combining Eq.~\eqref{equation decomposition 2} and Eq.~\eqref{ereg lower}, we have
\begin{equation}\label{equation decomposition 3}
    \begin{aligned}
    \frac{1}{2}\left\|\mu_{\tau+1}-x_\star\right\|_{\Sigma_{\tau+1}^{-1}}^2 \leq&\frac{1}{2}\left\|\mu_1-x_\star\right\|_{\Sigma_1^{-1}}^2+\frac{\eta^2}{2} \sum_{t=1}^\tau\left\|g_t\right\|_{\Sigma_{t+1}}^2+\operatorname{Error}.
\end{aligned}
\end{equation}
Now we can pick proper constants to make sure the right hand of Eq.~\eqref{equation decomposition 3} is under $\frac{1}{2\lambda^2L^2}$, which implies that $\tau=n$ with high probability. Recall that our goal is to bound $\operatorname{Reg}_\tau(x_\star)$ and it suffices to use Eq.~\eqref{equation decomposition 1} again. That is $\operatorname{Reg}_n(x_\star)=\mathcal{O}(\frac{1}{\lambda^2L\eta})$.


%% file: tex/main/convergence_rate.tex
Similar to Eq.~\eqref{rate eg}, it's clear that if there exists $k>0$ such that $\Sigma_t^{-1}=\Omega(t^{k})\idm_d$, then $\|X_t-x_\star\|_{2}$ is $\tilde{\mathcal{O}}(t^{-k/2}).$ In this section, we present the sketch of proof for Theorem \ref{rate >1} by showing that $\Sigma^{-1}_t$ grows in the order of $t^{\frac2\bc}$ when $f(x)$ is $(\beta,\bc)$-convex, $1<\bc\leq 2$ (i.e., Lemma \ref{growth rate >1}), which contributes to faster convergence rates, $t^{-\frac1\bc}$, of $\|X_t-x_\star\|_2$. Then we introduce the technical challenge for Theorem \ref{q=1 result}.
\begin{Rema}
    Recall that the real action is $\frac{X_t}{\pi^+(X_t)}$. Actually,  $\left\|\frac{X_t}{\pi^+(X_t)}-x_\star\right\|_2$ has the same convergence rate by Lemma \ref{xi lip} and that $\pi^+(X_t)\ge1.$
\end{Rema}
\begin{Lem}\label{growth rate >1}
    If $f(x)$ is $(\beta,\bc )$-convex on $\mathcal{K}$, $1<\bc\leq 2$, $\frac{r}{\sqrt{2}\sigma}\ge 5d$, $\lambda\leq\frac{1}{10dL}$ and $\sigma^{-2}\ge \Theta$, where $\Theta=\left(\frac{\bc-1}{30}\right)^{2/\bc}\beta^{\frac2\bc}d^{-\frac1\bc}\left(\frac{r}{\sqrt{2}R}\right)^{\frac{2(d-1)}{\bc}}\eta^{\frac2\bc}\lambda^{\frac6\bc-2}L^{\frac4\bc-2}$, then in Algorithm \ref{algorithm 1},
     $\Sigma_t^{-1}\succeq\frac{\Theta}{16}t^{\frac2\bc}\idm_d$ for all $t\leq\tau$. 
\end{Lem}
\paragraph{Sketch of proof of Lemma \ref{growth rate >1}}
By the updating rule of $\Sigma_t^{-1}$, informally,
 \[
 \Sigma^{-1}_{t+1}-\Sigma^{-1}_{t}\gtrsim \eta H_t\approx\eta\E_{t-1}[H_t]= \eta\nabla^2 s_t(\mu_t).
 \]
Then recalling the definition of $H_t$ in Eq.~\eqref{Ht} and by Lemma \ref{convex lip}, we have
\begin{equation}\label{EtHt}
    \begin{aligned}
        \E_{t-1}[H_t]=\frac{\lambda}{(1 {-} \lambda)^2}\E_{t-1}\left[e(\widetilde{X}_t)\left\{\Sigma_t^{-1}(\widetilde{X}_t {-} \mu_t)(\widetilde{X}_t {-} \mu_t)^{\top}\Sigma^{-1}_t/(1-\lambda)^2-\Sigma_t^{-1}\right\}\right]\succeq\lambda\E_{t-1}[\nabla^2 e(\widetilde{X}_t)],
    \end{aligned}
\end{equation}
where $\widetilde{X}_t\sim\N(\mu_t,(1-\lambda)^2\Sigma_t)$ and $e(x)$ is the convex extension of $f(x)$. Lemma \ref{convex lip} is a general version of Stein's Lemma \citep{stein1981estimation} for convex functions, which is evidently true if we can exchange expectation and differential freely. Its proof is given in Appendix \ref{exchange}.

Note that $e(x)$ equals $f(x)$ on $\mathcal{K}$ and when $f(x)$ is $(\beta, \bc)$-convex, $f(x)-\beta\|x-x_\star\|_2^\bc$ is convex on $\mathcal{K}$. Then by Eq.\eqref{EtHt}, we have
\begin{equation}\label{xtk}
    \E_{t-1}[H_t]\succeq\lambda\beta\E_{t-1}[\nabla^2 \|\widetilde{X}_t-x_\star\|_2^\bc\cdot\1_{\{\widetilde{X}_t\in\mathcal{K}\}}]\succeq\lambda\beta \bc(\bc-1)\E_{t-1}[\|\widetilde{X}_t-x_\star\|_2^{\bc-2}\idm_d\cdot\1_{\{\widetilde{X}_t\in\mathcal{K}\}}],
\end{equation}
where the second inequality follows from Lemma \ref{nabla2}.
Let $U=(1 {-}\lambda)^{-1}\Sigma_t^{-1/2}(\widetilde{X}_t {-} \mu_t)\sim\N(0,\idm_d)$. Then
\[
\|\widetilde{X}_t-x_\star\|_2\leq \|\mu_t-x_{\star}\|_2+(1-\lambda)\|\Sigma_t\|^{1/2}\|U\|_2\leq \|\Sigma_t\|^{1/2}\left(\widetilde{\mathcal{O}}(1)+\|U\|_2\right),
\]
where in the second inequality we used that $\|\mu_t-x_{\star}\|_{\Sigma_t^{-1}}^2=\widetilde{\mathcal{O}}(1)$, which is larger than $\|\mu_t-x_{\star}\|_2^2\|\Sigma_t\|^{-1}$ because $\|\Sigma_t\|^{-1}$ equals the minimal eigenvalues of $\Sigma_t^{-1}$. Then when $\bc\leq 2$, by Eq.\eqref{xtk},
\begin{equation}\label{xtkc}
    \E_{t-1}[H_t]\succeq\lambda\beta \bc(\bc-1)\|\Sigma_t\|^{\frac{\bc-2}{2}}\idm_d\E_{t-1}[(\|U\|_2+\widetilde{\mathcal{O}}(1))^{\bc-2}\cdot \1_{\{U\in\widetilde{\mathcal{K}}\}}],
\end{equation}
where $\widetilde{\mathcal{K}}=(1-\lambda)^{-1}\Sigma_t^{-1/2}(\mathcal{K}-\mu_t)$. For brevity, let's ignore $\1_{\{U\in\widetilde{\mathcal{K}}\}}$ and logarithmic terms for a while. Then it's easy to see that there exists a constant $C_{d,\bc}>0$
such that 
\begin{equation}\label{ignore}
    C_{d,\bc}\lesssim \bc(\bc-1)\E_{t-1}[(\|U\|_2+\widetilde{\mathcal{O}}(1))^{\bc-2}].
\end{equation}
Therefore, informally,
\begin{equation}\label{sigmat+1-sigmat}
    \Sigma^{-1}_{t+1}-\Sigma^{-1}_{t}\gtrsim\lambda\eta\beta C_{d,\bc}\|\Sigma_t\|^{\frac{\bc-2}{2}}\idm_d,
\end{equation}
hence, 
\begin{equation}\label{sigma min appro}
    \lambda_{\operatorname{min}}(\Sigma^{-1}_{t+1})-\lambda_{\operatorname{min}}(\Sigma^{-1}_{t})\gtrsim\lambda\eta\beta C_{d,\bc}\lambda_{\operatorname{min}}(\Sigma^{-1}_{t})^{\frac{2-\bc}{2}},
\end{equation}
where we again used that $\|\Sigma_t\|^{-1}$ equals the minimal eigenvalues of $\Sigma_t^{-1}$. Finally, by Lemma \ref{growth rate lamma}, we have $\lambda_{\operatorname{min}}(\Sigma^{-1}_{t})=\Omega(t^{\frac{2}{\bc}})$. 
To see this, one can naively let $\lambda_{\operatorname{min}}(\Sigma^{-1}_{t})\sim t^\alpha$ and then the left hand is in the order of $t^{\alpha-1}$ and the right $t^{\frac{(2-\bc)\alpha}{2}}$. Then we have $\alpha\ge\frac{2}{\bc}$ by solving that $\alpha-1\ge \frac{(2-\bc)\alpha}{2}.$

Finally, for the rigour of the proof, here we discuss how to handle $\1_{\{U\in\widetilde{\mathcal{K}}\}}.$ 
Since $U\sim\N(0,\idm_d)$, in order to obtain $C_{d,\bc}$ in Eq.~\eqref{ignore},
one should make sure that $\widetilde{\mathcal{K}}$ contains enough mass near the surface of $\mathbb{B}_d^d$, where $U$ is concentrated. We need the following lemma, whose proof can be found in Appendix \ref{proof for proportion}.
\begin{Lem}\label{proportion}
    If $\mathbb{B}_r^d\subset\mathcal{K}\subset\mathbb{B}_R^d$, then for all $x\in\mathcal{K}$
    , $\mathcal{K}$ contains a spherical cone, $\mathbb{C}_x^d$, of $\mathbb{B}_{r}^d(x)$, with the proportion of volume in $\mathbb{B}_{r}^d(x)$ not less than $\frac{1}{\sqrt{2\pi d}}\left(\frac{r}{\sqrt{2}R}\right)^{d-1}.$ 
\end{Lem}
If $\Sigma^{-1}_t\succeq\sigma^{-2}\idm_d$, then $\widetilde{\mathcal{K}}\supset\frac{\mathbb{C}_{\mu_t}^d-\mu_t}{(1-\lambda)\sigma}$. Since $\frac{U}{\|U\|_2}\indep\|U\|_2$, by Lemma \ref{proportion}, we have
\begin{equation}\label{widetilde k}
    \E_{t-1}[\|U\|_2^{\bc-2}\cdot \1_{\{U\in\widetilde{\mathcal{K}}\}}]\ge\frac{1}{\sqrt{2\pi d}}\left(\frac{r}{\sqrt{2}R}\right)^{d-1}\cdot\E_{t-1}[\|U\|_2^{\bc-2}\cdot \1_{\{\|U\|_2\leq\frac{r}{(1-\lambda)\sigma}\}}].
\end{equation}
And the last term is $\Omega(d^{\bc-2})$ when $\frac{r}{(1-\lambda)\sigma}\sim d$ because $\|U\|_2$ is concentrated with $d$. The formal version of the argument above is summarized in Appendix \ref{proof for growth rate >1}.\qed
\begin{Rema}
    Lemma \ref{nabla2} used in Eq.~\eqref{xtk} is also valid if $\|\cdot\|_p$ replaces $\|\cdot\|_2$ when $1<p<2$; we discuss the exponential dependence on the dimension of Lemma \ref{growth rate >1} in Appendix \ref{discussion of exponential}.
\end{Rema}

\paragraph{Technical challenge when $\bc=1$}\label{challenge q=1}
Things become strange when $\bc=1$. For instance, $(\beta,1)$-convex functions may not have a unique minimizer (consider the ReLU function, say $f(x)=x\cdot\1_{\{x\ge 0\}}$, which is clearly $(1,1)$-convex with $x_\star=0$). 
This will be fixed by Lemma \ref{1 qg} under the assumptions in Section \ref{main result q=1}. In addition, we will see that the arguments for $\bc>1$ above don't work. 

The first trouble occurs in tuning. Let's first take $d=1$ as an example, where $\nabla^2 \|x\|_2$ is just the Dirac delta function $\delta(x)$ and is relatively simple. Then similar to Eq.~\eqref{xtk}, informally, we have
\begin{equation}\label{d=1 new}
     \E_{t-1}[H_t]\approx
    \lambda\beta\E_{t-1}[\delta(\widetilde{X}_t-x_\star)]\approx \frac{2\beta}{\sqrt{2\pi}}e^{-\|\mu_t-x_\star\|_{\Sigma_t^{-1}}^2}\|\Sigma_t\|^{-1/2},
\end{equation}
where we used that $\widetilde{X}_t\sim\N(\mu_t,(1-\lambda)^2\Sigma_t)$ and the last term is just the density of $\N(\mu_t,(1-\lambda)^2\Sigma_t)$ at $x_\star$ by the property of the Dirac delta function (we also omit $\lambda$ for simplicity). This is very similar to the result in Eq.~\eqref{sigmat+1-sigmat} except $e^{-\|\mu_t-x_\star\|_{\Sigma_t^{-1}}^2}$. However, in the original \textbf{ONM}, $\|\mu_t-x_\star\|_{\Sigma_t^{-1}}^2\leq\frac{1}{2\lambda^2L^2}$ is $\widetilde{\mathcal{O}}(1)$, which is hard to improve to $\mathcal{O}(1)$ and will explode when put into the exponential function as in Eq.~\eqref{d=1 new}. The direct reason for this is that
the original analyses for Approximation Error in \cite{fokkema2024onlinenewtonmethodbandit} (say, Lemma \ref{reg qreg}) need that $\lambda\leq d^{-1}L^{-2}$ and we will solve this by a different decomposition (Eq.~\eqref{qg de 2 new}) from Eq.~\eqref{decomposition equation}, which makes use of that $\operatorname{sReg}_{\tau}(x_\star):=\sum_{t=1}^\tau\left(s_t(\mu_t)-s_t(x_\star)\right)\ge 0$ by Lemma \ref{smu>sx} under the assumptions in Section \ref{main result q=1}.
 
The second trouble occurs when $d\ge2$, where $\nabla^2 \|x\|_2=\frac{\idm_d-xx^\top/\|x\|_2^2}{\|x\|_2}$ is no longer positive definite, which makes its analysis much harder. 
One can't apply Lemma \ref{proportion} to handle the constrained case, which explains why we need the unconstrained assumption in Section \ref{main result q=1}. The details can be found in Appendix \ref{q=1 proof}, where we will show that for any $\kappa\in(0,1]$, $\Sigma_{t}^{-1}$ can grow at the rate of $t^{2-\kappa}$ if well-tuned. This is nearly optimal because $\Sigma_t^{-1}$ grows quadratically at most by Lemma \ref{sk upper bound for st''}.

%% file: tex/main/open.tex
In this paper we have studied a stochastic convex bandit problem with linear vanishing noise, and devised   a regularized online Newton method (\textbf{RONM}) for solving the problem. Our theoretical analysis has shown that \textbf{RONM} can reach a polylogarithmic regret in the time horizon when the loss function grows quadratically. We also analyze the convergence rate by capturing the growth rate of $\Sigma_t^{-1}$. There are several issues that remain open.   
  
First, we are not sure how necessary is the condition of quadratic growth (it is of course unnecessary,  because the regret is always $0$ if the loss function is a constant). As stated in \cite{lumbreras2024linear}, usual methods for deriving lower bounds on the minimax regret in this noise model fail because of the exploded KL divergence.
   
 Second, the experiments have shown that if $\|\cdot\|_2$ in the definition of $(\beta,\bc)$-convexity is replaced by 
 general $\ell_p$ norm $\|\cdot\|_p$, the growth rate of $\Sigma_t^{-1}$ seems unchanged. Actually, by Lemma \ref{nabla2}, our results still hold for all $1<p<2$. The experiments have also shown that the growth rate of $t^{2/\bc}$ may be true when $\bc>2$. We miss a general analysis for all $p\ge 1$ and $\bc\ge 1$.
 
Third, we have shown that it's possible to remove the exponential dependence on dimension when $\bc=2$, which makes little sense, however, because its contribution is far less than the regularized term. It remains unknown if this is possible for other $\bc<2$.

Finally, for $(\beta,1)$-convex functions, it would be desirable to develop an algorithm that reaches the convergence rate of $\frac{1}{t}$ and polylogarithmic regret simultaneously without  extra assumptions (see Section~\ref{main result q=1}). 
 

%% file: tex/append/Algorithm_and_proof.tex
\subsection{Proof for Theorem \ref{fqg}}\label{proof fgq}
We follow the proof in \cite{fokkema2024onlinenewtonmethodbandit} (Section \ref{pre onm}). Define the following quantities:
\[
\begin{aligned}&S_t=\sum_{u=1}^tH_u,&&\bar{S}_t=\sum_{u=1}^t\mathbb{E}_{u-1}[H_u]=\sum_{u=1}^t \nabla^2 s_u(\mu_u).\end{aligned}
\]
Let
\[
\bar{\Sigma}_t^{-1}=\Sigma_1^{-1}+\eta\sum_{u=1}^{t-1}\nabla^2q^\gamma(\mu_t)=\Sigma_1^{-1}+\eta\bar{S}_{t-1}/2+(t-1)\eta\gamma\idm_d,
\]
and $F_t=\frac12\left\|\mu_t-x_\star\right\|_{\Sigma_t^{-1}}^2$. We now make use of a stopping time to prove that $F_t$ won't be too large with high probability.
\begin{Def}\label{sk stopping}
    $\text{Let }\tau\text{ be the first round when one of the following does }\text{not hold:}$
    \[
    \begin{aligned}
    &\text{(a) } F_{\tau+1}\leq\frac{1}{2\lambda^{2}L^{2}};&&\\
    &\text{(b) } \Sigma_{\tau+1}\text{ is positive definite};&&\\
    &\text{(c) } \frac{1}{2}\bar{\Sigma}_{\tau+1}^{-1}\preceq\Sigma_{\tau+1}^{-1}\preceq\frac{3}{2}\bar{\Sigma}_{\tau+1}^{-1}.&&
    \end{aligned}
    \]
    $\text{In case(a)-(c) hold for all rounds }t\leq n\text{, then }\tau\text{ is defined to be }n.$
\end{Def}
\noindent
For all $t\leq\tau,$ we have
        \[
        \Sigma_t\succeq\frac{1}{2}\bar{\Sigma}_t\succeq\frac{1}{2\sigma^2}\idm_d+\frac{(t-1)\eta\gamma}2\idm_d,
        \]
        then noting that $\sigma^{-2}\ge\eta\gamma$ and by Lemma \ref{sk upper bound for st''}, we have
\begin{equation*}
    \frac{\eta\gamma t}{2}\idm_d\preceq\frac{1}{2}\bar{\Sigma}_{t}^{-1}\preceq\Sigma_{t}^{-1}\preceq\frac{3}{2}\bar{\Sigma}_{t}^{-1}\preceq\frac{3t^2h}{2}\idm_d.
\end{equation*}

\subsubsection*{Step 1: Concentration}
Define events $\operatorname{E}_1$ and $\operatorname{E}_2$ by 
\[
\begin{aligned}
    \operatorname{E}_1=\left\{\max_{1\leq t\leq\tau}\frac{|\xi_t|}{\|X_t-x_\star\|_2}\leq\frac{RL^{1/2}}{r}\right\},\quad
    \operatorname{E}_2=\left\{\max_{1\leq t\leq\tau}\|X_t-\mu_t\|_{\Sigma_t^{-1}}\leq d^{1/2}L^{1/2}\right\},
\end{aligned}
\]
where $\xi_t=\pi^+(X_t)\varepsilon_t$ and we denote that $\frac{|\xi_t|}{\|X_t-x_\star\|_2}=1$ when $X_t=x_\star$ since now $\xi_t=0, a.s.$ Note that $L=\Omega(\log(\max(n,1/\delta)))$, thus we have
\begin{Lem}\label{sk e1}
$\P(\operatorname{E}_1\cap\operatorname{E}_2)\ge1-2\delta/5.$
\end{Lem}
\noindent
Its proof can be found in Appendix \ref{proof for sk e1}. Since $\lambda\leq d^{-1/2}L^{-3/2}$, we have $\sqrt{dL}\leq\frac{1}{\lambda L}$. Then on $\operatorname{E}_1\cap\operatorname{E}_2$, for all $t\leq\tau$,
\begin{equation}\label{lambda1}
    \left\|X_t-x_\star\right\|_{\Sigma_t^{-1}}\leq \|X_t-\mu_t\|_{\Sigma_t^{-1}}+\left\|\mu_t-x_\star\right\|_{\Sigma_t^{-1}}\leq 2\lambda^{-1}L^{-1}.
\end{equation}
Recall that $\frac{\eta\gamma t}{2}\idm_d\preceq\Sigma_{t}^{-1}$. Then we have
\begin{equation}\label{x xstar l2}
    \left\|X_t-x_\star\right\|_{2}\leq \frac{2}{\sqrt{\eta\gamma t}}\left\|X_t-x_\star\right\|_{\Sigma_t^{-1}}\leq 
\frac{4}{\lambda L\sqrt{\eta\gamma t}}.
\end{equation}
Hence for $t\ge 2$,
\begin{equation}\label{Zt}
    \begin{aligned}
    |Z_t|\leq& |e(X_t)-e(x_\star)|+|e(X_{t-1})-e(x_\star)|+|\xi_t|+|\xi_{t-1}|\\
    \leq &\left(\lip(e)+\frac{RL^{1/2}}{r}\right)(\left\|X_t-x_\star\right\|_{2}+\left\|X_{t-1}-x_\star\right\|_{2})
    \leq \frac{H}{\lambda\sqrt{L\eta\gamma t}},
\end{aligned}
\end{equation}
where the final inequality used the definition of $H$ and Lemma \ref{extension}. For $t=1$, we also have
\begin{equation}\label{t=1}
    \begin{aligned}
    |Z_1|=|Y_1|&\leq |\pi^+(X_1)|\cdot|f(\frac{X_1}{\pi^+(X_1)})|+|\xi_1|\\
    &\leq1+\frac{1}{r}\|X_1-\mu_1\|_{2}+\frac{RL^{1/2}}{r}\|X_1-x_\star\|_2\leq \frac{H}{\lambda\sqrt{L\eta\gamma}},
\end{aligned}
\end{equation}
where the second inequality follows from that $|f(\frac{X_1}{\pi^+(X_1)})|\leq 1$ and $|\pi^+(X_1)|\leq\frac{1}{r}\|X_1-\mu_1\|_{2}+|\pi^+(\mu_1)|=\frac{1}{r}\|X_1-\mu_1\|_{2}+1$, since $\lip(\pi^+)\leq1/r$ and $\mu_1=0$. The third inequality used that $\frac{H}{\lambda\sqrt{L\eta\gamma}}\ge 3.$
 \noindent
 Therefore
 \begin{equation}\label{sum zt2}
     \sum_{t=1}^\tau Z_t^2\leq \frac{H^2}{L\lambda^2\eta\gamma}\sum_{t=1}^\tau 1/t\leq\frac{H^2}{\lambda^2\eta\gamma}.
 \end{equation}
 Similarly, we have the following lemma. 
 \begin{Lem}\label{sk Z V}
     Let $Z_{\max }=\max _{1 \leq t \leq \tau}\left(\left|Z_t\right|+\mathbb{E}_{t-1}\left[\left|Z_t\right|\right]\right)$ and $V_{\tau}=\sum_{t=1}^\tau \E_{t-1}[Z_{t}^2]$. If $\frac{H}{\lambda\sqrt{L\eta\gamma}}\ge 3$ and $\lambda\leq\frac{1}{2\sqrt{d}L}$, then on $\operatorname{E}_1\cap\operatorname{E}_2$,
     \[
     (a)\,Z_{\max}\leq \frac{H}{3\lambda\sqrt{L\eta\gamma}}
     ;
     \quad\quad\quad (b)\,V_{\tau}\leq 
    \frac{H^2}{9\lambda^2\eta\gamma}.
     \]
 \end{Lem}
 \noindent
The proof of Lemma \ref{sk Z V} is deferred in Appendix \ref{proof for sk Z V}. These bounds provide a nice concentration for $\widehat{\operatorname{q  Reg}_\tau}$. Define $\operatorname{E}_3$ by
\[
\operatorname{E}_3=\left\{\operatorname{q  Reg}_{\tau}\left(x_{\star}\right) \leq \widehat{\operatorname{q  Reg}}_{\tau}\left(x_{\star}\right)+1+
H\sqrt{\frac{L}{\lambda^4\eta\gamma}}\right\}.
\]
\begin{Lem}\label{sk e3}
If $\frac{H}{\lambda\sqrt{L\eta\gamma}}\ge 3$ and $\lambda\leq\frac{1}{2\sqrt{d}L}$, then
$\P(\operatorname{E}_1\cap\operatorname{E}_2\cap\operatorname{E}_3)\ge1-3\delta/5.$
\begin{proof}
    By the definition of $\widehat{q}  _t$, we have $\widehat{q}  _t(\mu_t)=q  _t(\mu_t)=0$, hence
    \[\operatorname{q  Reg}_{\tau}\left(x_{\star}\right) -\widehat{\operatorname{q  Reg}}_{\tau}\left(x_{\star}\right)=\sum_{t=1}^\tau\left(\hat{q}  _t\left(x_{\star}\right)-q  _t\left(x_{\star}\right)\right).
    \]
    Since $\max _{1 \leq t \leq \tau} \lambda\left\|x_\star-\mu_t\right\|_{\Sigma_t^{-1}} \leq L^{-1/2}$, by Lemma \ref{q con}
    , with probability at least $1-\delta/5$, 
    \[
\sum_{t=1}^\tau\left(\hat{q}  _t\left(x_{\star}\right)-q  _t\left(x_{\star}\right)\right) \leq 1+\frac{1}{\lambda}\left[\sqrt{V_{\tau}L}+Z_{\max } L\right] .
    \]
    Then when $\operatorname{E}_1$ and $\operatorname{E}_2$ both happen, by Lemma \ref{sk Z V},
    $$
    \sum_{t=1}^\tau\left(\hat{q}  _t\left(x_{\star}\right)-q  _t\left(x_{\star}\right)\right) \leq 1+H\sqrt{\frac{L}{\lambda^4\eta\gamma}}.
    $$
\end{proof}
\end{Lem}
\noindent
We also need that $\operatorname{Reg}_n(x)=\sum_{t=1}^n\left(f\left(\frac{X_t}{\pi^+(X_t)}\right)-f\left(x\right)\right)$ is well-concentrated around
\[
\widetilde{\Reg}_n(x):=\sum_{t=1}^n\left(\E_{t-1}\left[f\left(\frac{X_t}{\pi^+(X_t)}\right)\right]-f\left(x\right)\right),
\]
hence define $\operatorname{E}_4$ to be the event that
\[
\operatorname{E}_4=\left\{
\operatorname{Reg}_\tau(x_\star)\leq \widetilde{\Reg}_\tau(x_\star)+H\sqrt{\frac{L}{\lambda^2\eta\gamma}}
\right\},
\]
then similarly we have
\begin{Lem}\label{e4}
If $\frac{H}{\lambda\sqrt{L\eta\gamma}}\ge 3$ and $\lambda\leq\frac{1}{2\sqrt{d}L}$, then $\P(\operatorname{E}_1\cap\operatorname{E}_2\cap\operatorname{E}_4)\ge1-3\delta/5.$
\end{Lem}
\noindent
The proof of Lemma \ref{e4} can be found in Appendix \ref{proof for e4}. Finally, let $\operatorname{E}_5$ be the event
\[
\operatorname{E}_5=\left\{-\frac{Hd^2L^2}{\sqrt{\eta\gamma}}\cdot\frac32\bar{\Sigma}_\tau^{-1} \preceq S_\tau-\bar{S}_\tau \preceq \frac{Hd^2L^2}{\sqrt{\eta\gamma}} \cdot\frac32\bar{\Sigma}_\tau^{-1}\right\} .
\]
\begin{Lem}\label{sk e5}
If $\frac{H}{\lambda\sqrt{L\eta\gamma}}\ge 3$ and $\lambda\leq\frac{1}{2\sqrt{d}L}$, then    $\P(\operatorname{E}_1\cap\operatorname{E}_2\cap\operatorname{E}_5)\ge1-3\delta/5.$
    \begin{proof}
        By Lemma \ref{s'' con},
        with $\Sigma^{-1}=\frac{3}{2} \bar{\Sigma}_\tau^{-1}$, with probability at least $1-\delta/5$,

$$
\begin{aligned}
& \bar{S}_\tau-S_\tau \preceq \lambda L^2\left[1+\sqrt{d V_\tau}+d^2 Z_{\max }\right] \frac{3}{2} \bar{\Sigma}_\tau^{-1} \\
& S_\tau-\bar{S}_\tau \preceq \lambda L^2\left[1+\sqrt{d V_\tau}+d^2 Z_{\max }\right] \frac{3}{2} \bar{\Sigma}_\tau^{-1},
\end{aligned}
$$
then it suffices to apply Lemma \ref{sk Z V} when $\operatorname{E}_1$ and $\operatorname{E}_2$ both happen.
    \end{proof}
\end{Lem}
Let $\operatorname{E}=\operatorname{E}_1 \cap \operatorname{E}_2 \cap \operatorname{E}_3 \cap \operatorname{E}_4 \cap \operatorname{E}_5$ be the intersection of all these high probability events. Then, $\mathbb{P}(E) \geq 1-\delta$. For the remainder of the proof, we bound the regret on $\operatorname{E}$.
\subsubsection*{Step 2: Regret decomposition}
Now we present the explicit expressions of the $\operatorname{Error}$ term in Eq.~\eqref{equation decomposition 1} and Eq.~\eqref{equation decomposition 3}. First, by the definition of $\operatorname{E}_3$, the Estimation Error can be bounded by
\begin{equation}\label{estimation error}
    \operatorname{qReg}_\tau(x_\star)-\widehat{\operatorname{qReg}}_\tau(x_\star))\leq 1+
H\sqrt{\frac{L}{\lambda^4\eta\gamma}}.
\end{equation}
Then for the Approximation Error in Eq.~\eqref{equation decomposition 1}, it is $\operatorname{Reg}_\tau(x_\star)-\operatorname{qReg}_\tau(x_\star)$. By the extension,
\begin{equation}\label{r e}
    \begin{aligned}
    \operatorname{Reg}_{\tau}(x_\star)\leq \widetilde{\operatorname{Reg}}_{\tau}(x_\star)+H\sqrt{\frac{L}{\lambda^2\eta\gamma}}
    \leq \operatorname{e Reg}_{\tau}(x_\star)+H\sqrt{\frac{L}{\lambda^2\eta\gamma}}
    ,
\end{aligned}
\end{equation}
where $\operatorname{e Reg}_{\tau}(x_\star):=\sum_{t=1}^\tau \left(\E_{t-1}[e(X_t)]-e(x_\star)\right)$, the first inequality used the definition of $\operatorname{E}_4$ and the second inequality follows from Lemma \ref{e ep} (d). For small enough $\delta$, by Lemma \ref{reg qreg}, we have
\begin{equation}\label{qg eq}
    \operatorname{eReg}_{\tau}(x_\star)\leq
    \operatorname{q Reg}_{\tau}\left(x_{\star}\right)+\sum_{t=1}^\tau\frac{2}{\lambda}\operatorname{tr}(\nabla^2 s_t(\mu_t)\Sigma_t)+1.
\end{equation}
By combining Eq.~\eqref{r e} and Eq.~\eqref{qg eq}, the Approximation Error in Eq.~\eqref{equation decomposition 1} can be bounded by
\begin{equation}\label{approximation error 1}
    \operatorname{Reg}_\tau(x_\star)-\operatorname{qReg}_\tau(x_\star)\leq H\sqrt{\frac{L}{\lambda^2\eta\gamma}}+\sum_{t=1}^\tau\frac{2}{\lambda}\operatorname{tr}(\nabla^2 s_t(\mu_t)\Sigma_t)+1.
\end{equation}
Then by Eq.~\eqref{estimation error} and Eq.~\eqref{approximation error 1}, Eq.~\eqref{equation decomposition 1} becomes
\begin{equation}\label{qg de 1}
    \begin{aligned}
    \frac{1}{2}\left\|\mu_{\tau+1}-x_\star\right\|_{\Sigma_{\tau+1}^{-1}}^2 \leq&\frac{R^2}{2\sigma^2}+\frac{\eta^2}{2} \sum_{t=1}^\tau\left\|g_t\right\|_{\Sigma_{t+1}}^2+2\eta+\frac{H\sqrt{\eta L}}{\lambda^2\sqrt{\gamma }}+\frac{\eta\gamma}{2}\sum_{t=1}^\tau\|\mu_t-x_\star\|_2^2 \\
     & +\sum_{t=1}^\tau\frac{2\eta}{\lambda}\operatorname{tr}(\nabla^2 s_t(\mu_t)\Sigma_t)
    +\frac{H\sqrt{\eta L}}{\lambda\sqrt{\gamma }}-\eta\operatorname{Reg}_{\tau}(x_\star).
\end{aligned}
\end{equation}
For the Approximation Error in Eq.~\eqref{equation decomposition 3}, again, by Lemma \ref{reg qreg}, we have
$$\widetilde{\operatorname{e Reg}}_\tau(x_\star)-\operatorname{qReg}_\tau(x_\star)\leq\sum_{t=1}^\tau\frac{2}{\lambda}\operatorname{tr}(\nabla^2 s_t(\mu_t)\Sigma_t)+1.$$ Then similarly, Eq.~\eqref{equation decomposition 3} now becomes
\begin{equation}\label{qg de 2}
    \begin{aligned}
    \frac{1}{2}\left\|\mu_{\tau+1}-x_\star\right\|_{\Sigma_{\tau+1}^{-1}}^2 \leq& \frac{R^2}{2\sigma^2}+\frac{\eta^2}{2} \sum_{t=1}^\tau\left\|g_t\right\|_{\Sigma_{t+1}}^2+2\eta+
     \frac{H\sqrt{\eta L}}{\lambda^2\sqrt{\gamma }}+\sum_{t=1}^\tau\frac{2\eta}{\lambda}\operatorname{tr}(\nabla^2 s_t(\mu_t)\Sigma_t).
\end{aligned}
\end{equation}

\subsubsection*{Step 3: Basic bounds}
We first bound the gradient norm term $\left\|g_t\right\|_{\Sigma_{t+1}}^2$. For all $t\leq\tau$, by Definition \ref{sk stopping} (c), one can see that
$\left\|g_t\right\|_{\Sigma_{t+1}}^2 \leq 2\left\|g_t\right\|_{\bar{\Sigma}_{t+1}}^2 \leq 2\left\|g_t\right\|_{\bar{\Sigma}_t}^2 \leq 3\left\|g_t\right\|_{\Sigma_t}^2$ (it's true when $t=\tau$ by Eq. \eqref{condition 2} in Step 4).  
Then by Lemma \ref{Rt 3} and noting that $\lambda\leq1/2$, 
we have
\begin{equation}\label{gt}
    \sum_{t=1}^\tau\left\|g_t\right\|_{\Sigma_{t+1}}^2\leq3\sum_{t=1}^\tau\left\|g_t\right\|_{\Sigma_{t}}^2\leq  108\sum_{t=1}^\tau Z_t^2\|X_t-\mu_t\|_{\Sigma_{t}^{-1}}^2\leq  dL\sum_{t=1}^\tau Z_t^2\leq \frac{H^2dL}{\lambda^2\eta\gamma},
\end{equation}
where the final inequality used Eq.~\eqref{sum zt2}.

Then we apply Lemma \ref{trace} to bound the trace term $\operatorname{tr}(\nabla^2s_t(\mu_t)\Sigma_t)$. We should first check the condition that
\[
\eta\|\Sigma_t^{1/2}\nabla^2 s_t(\mu_t)\Sigma_t^{1/2}\|\leq 1,
\]
which is true because by Lemma \ref{upper bound for s''} (b)
\begin{equation}\label{condition 1}
    \eta\|\Sigma_t^{1/2}\nabla^2 s_t(\mu_t)\Sigma_t^{1/2}\|\leq\frac{\eta\lambda \operatorname{lip}(e)}{1-\lambda} \sqrt{d\|\Sigma_t\|}\leq H\eta\lambda\sigma\sqrt{d}\leq 1,
\end{equation}
where the second inequality used that  $\Sigma_t\preceq2\sigma^2\idm_d$.

Therefore,
\begin{equation}\label{trace dl}
    \sum_{t=1}^{\tau}\frac{\eta}{\lambda}\operatorname{tr}(\nabla^2 s_t(\mu_t)\Sigma_t)\leq 
\frac{8}{\lambda}\log\operatorname{det}\left(\sigma^2\bar{\Sigma}_{\tau+1}^{-1}\right)\leq \frac{8d}{\lambda}\log\left(\|\sigma^2\bar{\Sigma}_{\tau+1}^{-1}\|\right)
\leq \frac{dL}{\lambda},
\end{equation}
where the final inequality used Lemma \ref{sk upper bound for st''} and that $L\ge C
\max(
\log h, \log n)
,$ where $C$ is large enough.

\subsubsection*{Step 4: Proof for $\tau=n$, the regret bound and the convergence rate}
First, by the definition of $\operatorname{E}_5$,  since $\eta \frac{Hd^2L^2}{\sqrt{\eta\gamma}}\leq 2/3$, 
\begin{equation}\label{condition 2}
    \Sigma_{\tau+1}^{-1}=\Sigma_1^{-1}+\eta \left(\frac{1}{2}S_\tau+\gamma\tau\idm_d\right) \preceq \Sigma_1^{-1}+\eta \left(\frac{1}{2}\bar{S}_\tau+\gamma\tau\idm_d\right)+\frac12\eta \frac{Hd^2L^2}{\sqrt{\eta\gamma}} \cdot\frac32\bar{\Sigma}_\tau^{-1} \preceq \frac{3}{2} \bar{\Sigma}_{\tau+1}^{-1}.
\end{equation}
Similarly, $\Sigma_{\tau+1}^{-1}\succeq\frac{1}{2} \bar{\Sigma}_{\tau+1}^{-1}$, then Definition \ref{sk stopping} (b) and (c) still hold. Then we should make sure that Definition \ref{sk stopping} (a) is also valid. Combining Eq.~\eqref{qg de 2}, Eq.~\eqref{gt} and Eq.~\eqref{trace dl} leads to
\begin{equation}\label{solution 2}
    \begin{aligned}
    \frac{1}{2}\left\|\mu_{\tau+1}-x_\star\right\|_{\Sigma_{\tau+1}^{-1}}^2 \leq& \frac{R^2}{2\sigma^2}+ \frac{\eta H^2dL}{2\lambda^2\gamma}+2\eta
    +
     \frac{H\sqrt{\eta L}}{\lambda^2\sqrt{\gamma }}+\frac{dL}{\lambda},
\end{aligned}
\end{equation}
then by the definitions of constants, the right hand is less than $\frac{1}{2\lambda^2L^2}$ and then clearly $\tau=n$ on $\operatorname{E}$. Similarly, using Eq.~\eqref{qg de 1} instead of Eq.~\eqref{qg de 2}, we have
\begin{equation*}\label{gq de 11}
\begin{aligned}
    0\leq\frac{1}{2}\left\|\mu_{\tau+1}-x_\star\right\|_{\Sigma_{\tau+1}^{-1}}^2 
    \leq &\frac{R^2}{2\sigma^2}+ \frac{\eta H^2dL}{2\lambda^2\gamma}+2\eta
    +
     \frac{H\sqrt{\eta L}}{\lambda^2\sqrt{\gamma }}+\frac{dL}{\lambda}
    +\frac{H\sqrt{\eta L}}{\lambda\sqrt{\gamma }}+\frac{\eta\gamma}{2}\sum_{t=1}^\tau\|\mu_t-x_\star\|_2^2-\eta\operatorname{Reg}_{\tau}(x_\star)
    ,
\end{aligned}
\end{equation*}
then by the definitions of constants and recalling that $\frac{\eta\gamma}{2}\|\mu_t-x_\star\|_2^2\leq\frac{1}{2\lambda^2L^2t}$ for all $t\leq\tau$, we have
 $\operatorname{Reg}_{n}(x_\star)=\mathcal{O}(\frac{1}{\eta\lambda^2L})=\mathcal{O}(H^4d^6L^{10}/\rho).$ All of the constraints of constants are summarized in Appendix \ref{cons app}. Finally, note that for all $t\leq n$, $\pi^+(X_t)\ge 1$, then by Lemma \ref{xi lip},  we have
 \begin{equation}\label{true convergence rate}
          \left\|\frac{X_t}{\pi^+(X_t)}-x_{\star}\right\|_2\leq(1+1/r)\|X_t-x_\star\|_2
    .
    \end{equation}
 Then the convergence rate follows from Eq.~\eqref{x xstar l2}.
\subsection{Proof for Theorem \ref{rate >1}}\label{proof rate>1}
If $f(x)$ is $(\beta,\bc)$-convex, $1<\bc\leq 2$, then by Lemma \ref{q qg}, $f(x)$ is $2^{\bc-1}\beta
$-QG on $\mathcal{K}$ and we can apply \textbf{RONM} to $f(x)$. By Theorem \ref{fqg}, with probability at least $1-\delta$,
    \[
    \Reg_n=\mathcal{O}(H^4d^6L^{10}/\beta).
    \]
Then by Lemma \ref{growth rate >1}, in Algorithm \ref{algorithm 1}, 
    for all $t\leq\tau$, $\Sigma_t^{-1}\succeq\frac{\Theta}{16}t^{\frac2\bc}\idm_d\vee\frac{\eta\gamma t}{2}\idm_d,$ where $$\Theta=\left(\frac{\bc-1}{30}\right)^{2/\bc}\beta^{\frac2\bc}d^{-\frac1\bc}\left(\frac{r}{\sqrt{2}}\right)^{\frac{2(d-1)}{\bc}}\eta^{\frac2\bc}\lambda^{\frac6\bc-2}L^{\frac4\bc-2}.$$ 
    Recall that for all $t\leq \tau$,
    \[
    \|X_t-x_\star\|_2^2\leq \lambda_{\operatorname{min}}^{-1}(\Sigma_t^{-1})\|X_t-x_{\star}\|_{\Sigma_t^{-1}}^{2}\leq \lambda_{\operatorname{min}}^{-1}(\Sigma_t^{-1})\cdot\frac{1}{\lambda^2L^2},
    \]
    then the convergence rate follows from Eq.~\eqref{true convergence rate}.

%% file: tex/append/q=1.tex
Recall that when $\bc=1$ we need extra assumptions in Section \ref{main result q=1}. By Assumption 2, there's no need for an extension and hence we apply \textbf{ONM} for unconstrained convex bandits (Algorithm \ref{algorithm 2}). Clearly, now Eq.~\eqref{EtHt} becomes
\begin{equation}\label{EtHt new}
    \E_{t-1}[H_t]=\frac{\lambda}{(1-\lambda)^2}\E_{t-1}\left[f(\widetilde{X}_t)\left\{\Sigma_t^{-1}(\widetilde{X}_t-\mu_t)(\widetilde{X}_t-\mu_t)^{\top}\Sigma^{-1}_t/(1-\lambda)^2-\Sigma_t^{-1}\right\}\right],
\end{equation} 
where $\widetilde{X}_t\sim\N(\mu_t,(1-\lambda)^2\Sigma_t)$.
\subsection{\texorpdfstring{The growth rate of $\Sigma_t^{-1}$}{Growth rate of Sigmat-1}}
In this section, we analyze the growth rate of $\Sigma_t^{-1}$ under assumptions in Section \ref{main result q=1}.
When $d=1$, $f(x)-\beta|x-x_\star|$ is convex on $\R.$ Then by Eq.~\eqref{EtHt new}, we have
\[
\begin{aligned}
    \E_{t-1}[H_t]
    \ge&\frac{\beta\lambda}{(1-\lambda)^2}\E_{t-1}\left[|\widetilde{X}_t-x_\star|\left\{\Sigma_t^{-1}(\widetilde{X}_t-\mu_t)(\widetilde{X}_t-\mu_t)^{\top}\Sigma^{-1}_t/(1-\lambda)^2-\Sigma_t^{-1}\right\}\right]\\
    =&\frac{\beta}{1-\lambda} \E_{t-1}\left[\frac{\d}{\d x}|x-x_\star|\big{|}_{x=\widetilde{X}_t}\cdot\Sigma_t^{-1}(\widetilde{X}_t-\mu_t)\right]\\
    =&\frac{2\beta}{1-\lambda} \E_{t-1}\left[\1_{\{\widetilde{X}_t\ge x_{\star}\}}\cdot\Sigma_t^{-1}(\widetilde{X}_t-\mu_t)\right],
\end{aligned}
\]
where the first equality used Lemma \ref{gx^2} and the second equality used that $\frac{\d}{\d x}|x-x_\star|=2\1_{\{x\ge x_\star\}}-1$ and $\E_{t-1}[\widetilde{X}_t-\mu_t]=0.$ Since $\bar{X}_t:=(1-\lambda)^{-1}\Sigma^{-1/2}_t(\widetilde{X}_t-\mu_t)\sim\N(0,1)$, by Lemma \ref{dirac}, we have\footnote{It's worth noting that when $d=1$, this is also true even in the constrained case. Because by Lemma \ref{ex beta 1}, the convex extension $e(x)$ is also $(\beta,1)$-convex, which helps us get rid of the discussion of $\mathcal{K}$}
\[
\E_{t-1}[H_t]\ge2\beta\|\Sigma_t\|^{-1/2}\E_{t-1}\left[\1_{\{\bar{X}_t\ge (1-\lambda)^{-1}\Sigma^{-1/2}_t(x_\star-\mu_t)\}}\cdot\bar{X}_t\right] =2\beta\|\Sigma_t\|^{-1/2}\cdot\frac{1}{\sqrt{2\pi}}e^{-\frac{\|\mu_t-x_\star\|_{\Sigma_t^{-1}}^2}{2(1-\lambda)^2}}.
\]
Hence, if $\lambda\leq1-\frac{1}{\sqrt{2}}$, we have
\begin{equation}\label{d=1}
    \E_{t-1}[H_t]\ge\frac{2\beta}{\sqrt{2\pi}}e^{-\|\mu_t-x_\star\|_{\Sigma_t^{-1}}^2}\|\Sigma_t\|^{-1/2}.
\end{equation}
\noindent
When $d\ge2$, $\nabla^2 \|x\|_2=\frac{\idm_d-\frac{xx^\top}{\|x\|_2^2}}{\|x\|_2}$, then by some different and relatively difficult analyses, we have:
\begin{Lem}\label{growth rate =1 lemma}
        If $d\ge 2$, and $X\sim\N(\mu, \Sigma)$, then 
    \[
    \E\left[\nabla^2\|x-x_{\star}\|_2\big{|}_{x=X}
    \right]\succeq6^{-\frac{d}2}e^{-\|\mu-x_\star\|^2_{\Sigma^{-1}}}\|\Sigma\|^{-\frac12}\idm_d/2.
    \]
\end{Lem}
\noindent
The proof is deferred in Appendix \ref{proof for growth rate =1 lemma}. Similar to Eq.~\eqref{EtHt}, by Lemma \ref{growth rate =1 lemma} and Eq.~\eqref{EtHt new}, we have
\begin{equation}\label{all d}
    \E_{t-1}[H_t]\succeq\beta 6^{-\frac{d}2}e^{-\|\mu_t-x_\star\|_{\Sigma_t^{-1}}^2}\|\Sigma_t\|^{-\frac12}\idm_d/2,
\end{equation}
which is also true when $d=1$ by Eq.~\eqref{d=1}, if $\lambda\leq1-\frac{1}{\sqrt{2}}$. Though with Lemma \ref{growth rate lamma}, it seems that $\Sigma_t^{-1}$ grows quadratically, we can't make sure that $\|\mu_t-x_{\star}\|^2_{\Sigma_t^{-1}}=\mathcal{O}(1)$ when this really happens. Fortunately, we still have a slightly weaker result:
\begin{Lem}\label{growth rate=1}
    If assumptions in Section \ref{main result q=1} are satisfied, $\lambda\leq 1-\frac{1}{\sqrt{2}}$, $0<\kappa\leq1$ and $\sigma^{-2}\ge \max\{\Theta,1\}$, where
    $\Theta=\beta^{2-\kappa}\eta^{2-\kappa}6^{-\frac{d(2-\kappa)}{2}}e^{-\frac{2-\kappa}{\lambda^2L^2}}/32$, then in Algorithm \ref{algorithm 2}, for all $t\leq\tau$, we have
    \[
    \Sigma_{t}^{-1}\succeq\frac{\Theta}{16} t^{2-\kappa}\idm_d.
    \]
\end{Lem}
\noindent
The proof can be found in Appendix \ref{proof for growth rate=1}. In other words, if well-tuned, the order of the growth rate of $\Sigma_t^{-1}$ for $(\beta,1)$-convex 
loss functions can be arbitrarily close to $t^2$.

\subsection{Algorithm}
We apply \textbf{ONM} for unconstrained convex bandits. We set the constants as follows:
\begin{equation}\label{alg 2 con}
    \begin{aligned}
    \sigma=1,\quad\lambda=\frac{1}{2L},\quad\eta^{\kappa}=\frac{\beta^{2-\kappa}
    }{10^7 6^{\frac{d(2-\kappa)}{2}}H^2d^5L^{6}},
\end{aligned}
\end{equation}
where $L=C[1+\log \max (n, d, H, 1/\beta, 1/\delta)], \,\delta=\operatorname{Poly}(1/n,1/d,1/H)\in(0,1)
    ,H=C'\max(G,1)$ and $C$ and $C'$ are sufficiently large universal constants.
\begin{algorithm}[H]
    \caption{\textbf{ONM} for unconstrained convex bandits
    }\label{algorithm 2}
    Require: $\eta,\lambda,\sigma>0$

    Set $\mu_1=0$, $\Sigma_1=\sigma^2\idm_d$ and $Y_0=0$

    \For{$t=1,2,\cdots,n$}{
        \vspace{1mm}
        ~sample $X_t$ from $\N(\mu_t,\Sigma_t)$ with density $p_t$
        \vspace{1mm}
        
        ~obeserve $Y_t=f(X_t)+\varepsilon_t$
        \vspace{1mm}

        ~let~ $R_t=\frac{p_t\left(\frac{X_t-\lambda\mu_t}{1-\lambda}\right)}{(1-\lambda)^dp_t(X_t)}$ and $Z_t=Y_t-Y_{t-1}$
        \vspace{1mm}

        ~compute $g_t=\frac{R_tZ_t\Sigma_t^{-1}(X_t-\mu_t)}{(1-\lambda)^2}$
        \vspace{1mm}

        ~compute $H_t=\frac{\lambda R_tZ_t}{(1-\lambda)^2}\left[\frac{\Sigma_t^{-1}(X_t-\mu_t)(X_t-\mu_t)^\top\Sigma_t^{-1}}{(1-\lambda)^2}-\Sigma_t^{-1}\right]$
        \vspace{1mm}

        ~$\Sigma_{t+1}^{-1}\gets\Sigma_t^{-1}+
       \frac{\eta}{2}H_t$
       \vspace{1mm}

        ~$\mu_{t+1}\gets\arg\min_{\mu\in \mathcal{K}}\|\mu-[\mu_t-\eta\Sigma_{t+1}g_t]\|_{\Sigma_{t+1}^{-1}}$
       
    }
\end{algorithm}

\subsection{Proof for Theorem \ref{q=1 result}}
We will retain most of the notations from the previous proof in Appendix \ref{proof main}, and we will explicitly point out any differences. Let
\[
\bar{\Sigma}_t^{-1}=\Sigma_1^{-1}+\eta\sum_{u=1}^{t-1}\nabla^2q(\mu_t)=\Sigma_1^{-1}+\eta\bar{S}_{t-1}/2.
\]
By Lemma \ref{sk upper bound for st''} and Lemma \ref{growth rate=1}, for all $t\leq\tau$, we have
\begin{equation*}
        \frac{\Theta}{16} t^{2-\kappa}\idm_d\preceq\Sigma_{t}^{-1}\preceq\frac{3t^2h}{2}\idm_d.
\end{equation*}

\subsubsection*{Step 1: Concentration}
Define events $\operatorname{E}_1$ and $\operatorname{E}_2$ by 
\[
\begin{aligned}
    \operatorname{E}_1=\left\{\max_{1\leq t\leq\tau}\frac{|\varepsilon_t|}{\|X_t-x_\star\|_2}\leq L^{1/2}\right\},\quad
    \operatorname{E}_2=\left\{\max_{1\leq t\leq\tau}\|X_t-\mu_t\|_{\Sigma_t^{-1}}\leq d^{1/2}L^{1/2}\right\}.
\end{aligned}
\]
  Similar to Lemma \ref{sk e1}, we also have $\P(\operatorname{E}_1\cap\operatorname{E}_2)\ge1-\delta/2.$ Let $J=\max(\sqrt{dL},\frac{1}{\lambda L})$, then on $\operatorname{E}_1\cap\operatorname{E}_2$, for all $t\leq\tau$
\begin{equation*}
    \left\|X_t-x_\star\right\|_{\Sigma_t^{-1}}\leq \|X_t-\mu_t\|_{\Sigma_t^{-1}}+\left\|\mu_t-x_\star\right\|_{\Sigma_t^{-1}}\leq 2J.
\end{equation*}
Recall that $\frac{\Theta}{16} t^{2-\kappa}\idm_d\preceq\Sigma_{t}^{-1}$, then we have
\begin{equation}\label{x xstar l2 new}
    \left\|X_t-x_\star\right\|_{2}\leq \frac{4}{\sqrt{\Theta}t^{1-\frac{\kappa}{2}}}\left\|X_t-x_\star\right\|_{\Sigma_t^{-1}}\leq 
\frac{8J}{\sqrt{\Theta}t^{1-\frac{\kappa}{2}}}.
\end{equation}
Then for $t\ge 2$,
\begin{equation*}
    \begin{aligned}
    |Z_t|\leq& |f(X_t)-f(x_\star)|+|f(X_{t-1})-f(x_\star)|+|\varepsilon_t|+|\varepsilon_{t-1}|\\
    \leq &(G+L^{1/2})(\left\|X_t-x_\star\right\|_{2}+\left\|X_{t-1}-x_\star\right\|_{2})
    \leq \frac{HJ\sqrt{L}}{\sqrt{\Theta}t^{1-\frac{\kappa}{2}}},
\end{aligned}
\end{equation*}
where the final inequality used the definition of $H$. Similar to Eq.~\eqref{t=1}, this is also true when $t=1$ if $\frac{HJ\sqrt{L}}{\sqrt{\Theta}}\ge 3$. And one can also show that on $\operatorname{E}_1\cap\operatorname{E}_2,$
\begin{equation*}
    (a)\,\sum_{t=1}^\tau Z_t^2\leq\frac{H^2J^2L^2}{\Theta};\quad(b)\,\,Z_{\max}\leq \frac{HJ\sqrt{L}}{3\sqrt{\Theta}};\quad(c)\,\,V_{\tau}\leq \frac{H^2J^2L^2}{9\Theta},
\end{equation*}
where we also used that $\sum_{n=1}^{+\infty}n^{\kappa-2}<L$. Define $\operatorname{E}_3$ by
\[
\operatorname{E}_3=\left\{\operatorname{q Reg}_{\tau}\left(x_{\star}\right) \leq \widehat{\operatorname{q Reg}}_{\tau}\left(x_{\star}\right)+1+
\frac{HJL^{\frac32}}{\lambda\sqrt{\Theta}}\right\}.
\]
Define $\operatorname{E}_4$ by
\[
\operatorname{E}_4=\left\{-\frac{\lambda d^2HJL^{3}}{\sqrt{\Theta}}\cdot\frac32\bar{\Sigma}_\tau^{-1} \preceq S_\tau-\bar{S}_\tau \preceq \frac{\lambda d^2HJL^{3}}{\sqrt{\Theta}}\cdot\frac32\bar{\Sigma}_\tau^{-1}\right\} .
\]
Let $\operatorname{E}=\operatorname{E}_1 \cap \operatorname{E}_2 \cap \operatorname{E}_3 \cap \operatorname{E}_4
$ be the intersection of all these events. Then similar to Lemma \ref{sk e3}
and Lemma \ref{sk e5}, when $\frac{HJ\sqrt{L}}{\sqrt{\Theta}}\ge3$, we still have $\mathbb{P}(E)\geq 1-\delta$. 
For the remainder of the proof we bound the convergence rate on $\operatorname{E}$.
\subsubsection*{Step 2: Regret decomposition}
Similar to Eq.~\eqref{qg de 2}, we have\footnote{Remember that we have explained that the original decomposition fails for large $\lambda$ in Section \ref{section 2/q}. For the same reason \textbf{RONM} can't be applied if a fast convergence rate is needed.}
\begin{equation}\label{qg de 2 new}
    \begin{aligned}
    \frac{1}{2}\left\|\mu_{\tau+1}-x_\star\right\|_{\Sigma_{\tau+1}^{-1}}^2 \overset{\text{Lemma \ref{regret for onm}}}{\leq}& \frac{R^2}{2\sigma^2}+\frac{\eta^2}{2} \sum_{t=1}^\tau\left\|g_t\right\|_{\Sigma_{t+1}}^2-\eta \widehat{\operatorname{qReg}}_\tau(x_\star)\\
    \overset{\operatorname{E}_3}{\leq} &\frac{R^2}{2\sigma^2}+\frac{\eta^2}{2} \sum_{t=1}^\tau\left\|g_t\right\|_{\Sigma_{t+1}}^2+\eta+
     \frac{\eta HJL^{\frac32}}{\lambda\sqrt{\Theta}}
    -\eta \operatorname{qReg}_\tau(x_\star)\\
    \overset{\text{Lemma \ref{sq}}}{\leq} &\frac{R^2}{2\sigma^2}+\frac{\eta^2}{2} \sum_{t=1}^\tau\left\|g_t\right\|_{\Sigma_{t+1}}^2+2\eta
    +
     \frac{\eta HJL^{\frac32}}{\lambda\sqrt{\Theta}}
    - \eta \operatorname{sReg}_{\tau}(x_\star)
    \\
    \overset{\text{Lemma \ref{smu>sx}}}{\leq} &\frac{R^2}{2\sigma^2}+\frac{\eta^2}{2} \sum_{t=1}^\tau\left\|g_t\right\|_{\Sigma_{t+1}}^2+2\eta
    +
     \frac{\eta HJL^{\frac32}}{\lambda\sqrt{\Theta}},
\end{aligned}
\end{equation}
where $\operatorname{sReg}_{\tau}(x_\star):=\sum_{t=1}^\tau\left(s_t(\mu_t)-s_t(x_\star)\right)\ge 0.$ 

\subsubsection*{Step 3: Basic bounds}
Similar to Eq.~\eqref{gt}, when $\lambda\leq\frac12,$
\begin{equation}\label{gt new}
    \sum_{t=1}^\tau\left\|g_t\right\|_{\Sigma_{t+1}}^2\leq  108dL\sum_{t=1}^\tau Z_t^2\leq \frac{dH^2J^2L^3}{\Theta}.
\end{equation}
\subsubsection*{Step 4: Proof for $\tau=n$ and the convergence rate}
First, by the definition of $\operatorname{E}_4$,  since $\frac{\eta\lambda d^2HJL^{3}}{\sqrt{\Theta}}\leq 2/3$, 
\begin{equation}\label{new def bc}
    \Sigma_{\tau+1}^{-1}=\Sigma_1^{-1}+\frac{\eta}{2}S_\tau \preceq \Sigma_1^{-1}+\frac{\eta}{2}\bar{S}_\tau+\frac{\eta\lambda d^2HJL^{3}}{\sqrt{\Theta}}\cdot\frac32\bar{\Sigma}_\tau^{-1} \preceq \frac{3}{2} \bar{\Sigma}_{\tau+1}^{-1}.
\end{equation}
Similarly, $\Sigma_{\tau+1}^{-1}\succeq\frac{1}{2} \bar{\Sigma}_{\tau+1}^{-1}$, hence, Definition \ref{sk stopping} (b) and (c) still hold. Then, combining Eq.~\eqref{qg de 2 new} and Eq.~\eqref{gt new}, we have
\begin{equation}\label{alg 2}
    \frac{1}{2}\left\|\mu_{\tau+1}-x_\star\right\|_{\Sigma_{\tau+1}^{-1}}^2\leq \frac{R^2}{2\sigma^2}+\frac{d\eta^2H^2J^2L^3}{2\Theta}+2\eta
    +
     \frac{\eta HJL^{\frac32}}{\lambda\sqrt{\Theta}},
\end{equation}
then by the definitions of constants, we can show that the right hand is less than $\frac{1}{2\lambda^2L^2}$\footnote{This is impossible when $\kappa=0$ since $\eta$ will be canceled in the right hand which forces $\lambda$ to be very small.} and then clearly $\tau=n$ on $\operatorname{E}$. Hence, by Eq.~\eqref{x xstar l2 new}, we have
\[
\|X_t-x_\star\|_2=\widetilde{\mathcal{O}}(6^{\frac{d(2-\kappa)}{2\kappa}}t^{-1+\frac{\kappa}{2}}).
\]
 All of the constraints of constants are summarized in Appendix \ref{cons app}.

%% file: tex/main/noise_scaled.tex
Recall that in the setting of stochastic convex bandits with noise scaled to $\sigma(x)$ (Section \ref{scaled noise}), there exists $\sigma(x):\mathcal{K}\to\R^+$ such that at round $t$ the noise $\varepsilon_t=\sigma(X_t)\cdot\bar{\varepsilon}_t$, where $\{\bar{\varepsilon}_t\}_{t=1}^n$ are independent and identically distributed $1$-subgaussian non-degenerate random variables. Hence, for any fixed action $X$ in $\mathcal{K}$, if the player repeatedly chooses $X$ twice and gets feedback $Y^{(1)}$ and $Y^{(2)}$, then let 
    \[
    W=|Y^{(1)}-Y^{(2)}|=\sigma(X)|\bar{\varepsilon}^{(1)}-\bar{\varepsilon}^{(2)}|,
    \]
    which has conditional expectation $\E[|\bar{\varepsilon}^{(1)}-\bar{\varepsilon}^{(2)}|]\sigma(X)$. By Lemma \ref{x-y subgaussian}, there exists $C>0$ such that $W-\E[W]$ is conditionally $C\sigma(X)$-subgaussian.

    Imagine that there is a bandit player and an intermediary and every time the player tells the intermediary that the player's choice is $X$ in $\mathcal{K}$, the intermediary secretly picks $X$ twice and then computes $W$, i.e., the absolute value of the difference of two times of feedback observed by the intermediary. Finally, the intermediary tells the player that the player suffers the loss of $W$. Now, from the player's perspective, the loss function is $\E[|\bar{\varepsilon}^{(1)}-\bar{\varepsilon}^{(2)}|]\sigma(x)$ and the noise is $C\sigma(X)$-subgaussian. 

    Therefore, when $\sigma(x)$ has the $\rho$-QG property, the player can just implement \textbf{RONM}, which promises that with probability at least $1-\delta$,
    $
        \sum_{t=1}^{n}\left(\sigma(X_t)-\sigma(X_\star)\right)=\mathcal{O}(\operatorname{polylog}(n)),
    $ 
    and $\|X_t-x_\star\|_2=\widetilde{\mathcal{O}}(t^{-1/2})$. Though the real choices are $\bar{X}_t=X_{\lceil t/2\rceil}$ for $t=1,\cdots, 2n$, it's clear that we still have $\|\bar{X}_t-x_\star\|_2=\widetilde{\mathcal{O}}(t^{-1/2})$. And when $f(x)-f(x_\star)\leq C^{-1}\sigma(x)$ for all $x\in\mathcal{K}$, we have \[
    \Reg_{2n}=\sum_{t=1}^{2n}\left(f\left(\bar{X}_t\right)-f\left(x_\star\right)\right)\leq 2C^{-1}\sum_{t=1}^{n}\left(\sigma(X_t)-\sigma(X_\star)\right)=\mathcal{O}(\operatorname{polylog}(n)).
    \]    
    The case for $(\beta,\bc)$-convexity can be shown by the same argument, which completes the proofs for Corollary \ref{equi f sigma} and Corollary \ref{multiplicative noise}.

%% file: tex/append/useful_facts.tex
\begin{Lem}\label{upper rho}
    If $f(x)$ is $\rho$-QG on $\mathcal{K}$ and $\sup_{x\in\mathcal{K}}\|x\|_2=R$, then $\rho\leq 8/R^2.$
    \begin{proof}
        It suffices to note that for all $x\in\mathcal{K}$, we have 
        \[
        1\ge |f(x)-f(x_\star)|\ge \frac{\rho}{2}\|x-x_\star\|_2^2,
        \]
        and $\sup_{x\in\mathcal{K}}\|x-x_\star\|_2\ge R/2$.
    \end{proof}
\end{Lem}
\begin{Lem}\label{q qg}
    If $f(x)$ is $(\beta,\bc)$-convex on $\mathcal{K}$ and $\bc>1$, then $f(x)-\beta\|x-x_\star\|_2^\bc\ge 0$. Moreover, if $\bc\leq2$ and $\mathcal{K}\subset\mathbb{B}_R^d$, then $f(x)$ is also $2\beta(2R)^{\bc-2}$-QG on $\mathcal{K}$.
    \begin{proof}
        Let $g(x)=f(x)-f(x_\star)-\beta\|x-x_\star\|_2^\bc$, then $g(x)$ is convex on $\mathcal{K}$ and it suffices to show that $g(x)\ge 0.$ Otherwise, there exists $y\in\mathcal{K}$ such that $g(y)<0.$ Let $x_t=x_\star+t(y-x_\star), \forall t\in[0,1]$, then by convexity of $g$, we have
        \[
        g(x_t)\leq g(y)t, \quad\forall t\in[0,1].
        \]
        Hence
        \[
        f(x_t)\leq f(x_\star)+g(y)t+\beta t^\bc\|y-x_\star\|_2^\bc, \quad\forall t\in[0,1].
        \]
        Taking $0<t<\left(\frac{-g(y)}{\beta \|y-x_\star\|_2^\bc}\right)^{\frac{1}{\bc-1}}$ leads to $f(x_t)<f(x_\star)$, which contradicts with that $x_\star$ is the minimizer and hence $f(x)\ge f(x_\star)+\beta\|x-x_\star\|_2^\bc$ for all $x\in\mathcal{K}$. Since for all $x\in\mathcal{K}$, we have $\|x-x_\star\|_2\leq2R$, then clearly,
        \[
        f(x)\ge \beta\|x-x_\star\|_2^\bc\ge\beta(2R)^{\bc-2}\|x-x_\star\|_2^2.
        \]
    \end{proof}
\end{Lem}
\begin{Lem}\label{1 qg}
    If $f(x)$ satisfies Assumptions 1 and 3 in Section \ref{main result q=1}
    , then for all $x\in\R^d$, $f(x)\ge f(x_\star)+\beta\|x-x_\star\|_2$ and $\beta\leq2$.
    \begin{proof}
        Let $g(x)=f(x)-f(x_\star)-\beta\|x-x_\star\|_2$, then $g(x)$ is convex on $\mathcal{K}$ and it suffices to show that $g(x)\ge 0.$ Otherwise, there exists $y\in\R^d$ such that $g(x_\star+y)<0,$ then $g(x_\star-y)=g(x_\star+y)<0$, which implies that $g(x_\star)\leq\frac{g(x_\star-y)+g(x_\star+y)}{2}<0$ and contradicts with that $g(x_\star)=0$. By the proof of Lemma \ref{upper rho}, we have $\beta\leq \frac2R$. 
    \end{proof}
\end{Lem}
\begin{Lem}\label{xi lip}
    For every $x\in\R^d$, we have
    $\|x-\pi^+(x)x_\star\|_2\leq (1+R/r)\|x-x_\star\|_2.$
    \begin{proof}
        By Lemma \ref{lip pi},
        $\lip(\pi)\leq 1/r$, then we also have $\lip(\pi^+)\leq 1/r$, since it suffices to check that for all $x\notin\mathcal{K}$ and $y\in\mathcal{K}$, we have
        \[
        |\pi^+(x)-\pi^+(y)|=|\pi^+(x)-1|=|\pi(x)-\pi(y')|\leq \lip(\pi)\|x-y'\|_2\leq \lip(\pi)\|x-y\|_2,
        \]
        where $y'$ is the intersection of the line segment connected by $x$ and $y$ and $\partial\mathcal{K}$. 
        
        Let $h(x)=x-\pi^+(x)x_\star$, then since $x_\star\in\mathcal{K}$, it's clear that $h(x_\star)=0.$ Therefore
        \[
        \|x-\pi^+(x)x_\star\|_2=\|h(x)-h(x_\star)\|_2\leq \|x-x_\star\|_2+|\pi^+(x)-\pi^+(x_\star)|\cdot\|x_\star\|_2,
        \]
        then the result follows from $\|x_\star\|_2\leq R$ and $|\pi^+(x)-\pi^+(x_\star)|\leq \|x-x_\star\|_2/r.$
    \end{proof}
\end{Lem}

\begin{Lem}\label{ex beta 1}
    If $d=1$, $f(x)$ is $(\beta,1)$-convex on $\mathcal{K}$ then its convex extension, $e(x)$ defined in Eq.~\eqref{extension defi} is $(\beta,1)$-convex on $\R$.
    \begin{proof}
        For a convex function $h(x)$, we use $h^{'}_+(x)$ and $h^{'}_-(x)$ to denote its right and left derivatives at $x$. Since $f(x)-\beta|x-x_\star|$ is convex on $\mathcal{K}$, $f^{'}_{+}(x)\ge \beta$ for $x\ge x_\star$ and $f^{'}_{-}(x)\leq -\beta$ for $x\leq x_\star$. Let $g(x)=e(x)-\beta|x-x_\star|$, then one can see that $g^{'}_{+}(x)$ is increasing in $[x_\star,+\infty)$ and $g^{'}_{-}(x)$ is also increasing in $(-\infty,x_\star].$ Note that $g^{'}_{+}(x_\star)=f^{'}_{+}(x_\star)-\beta\ge 0\ge f^{'}_{-}(x_\star)+\beta=g^{'}_{-}(x_\star)$, then it's clear that the derivative of $g(x)$ is increasing on $\R$, thus $g(x)$ is also convex on $\R.$
    \end{proof}
\end{Lem}

%% file: tex/append/proofforlemmas.tex
\subsection{Proof for Lemma \ref{proportion}}\label{proof for proportion}
\begin{figure}
\centering
\includegraphics[width=0.6\textwidth]{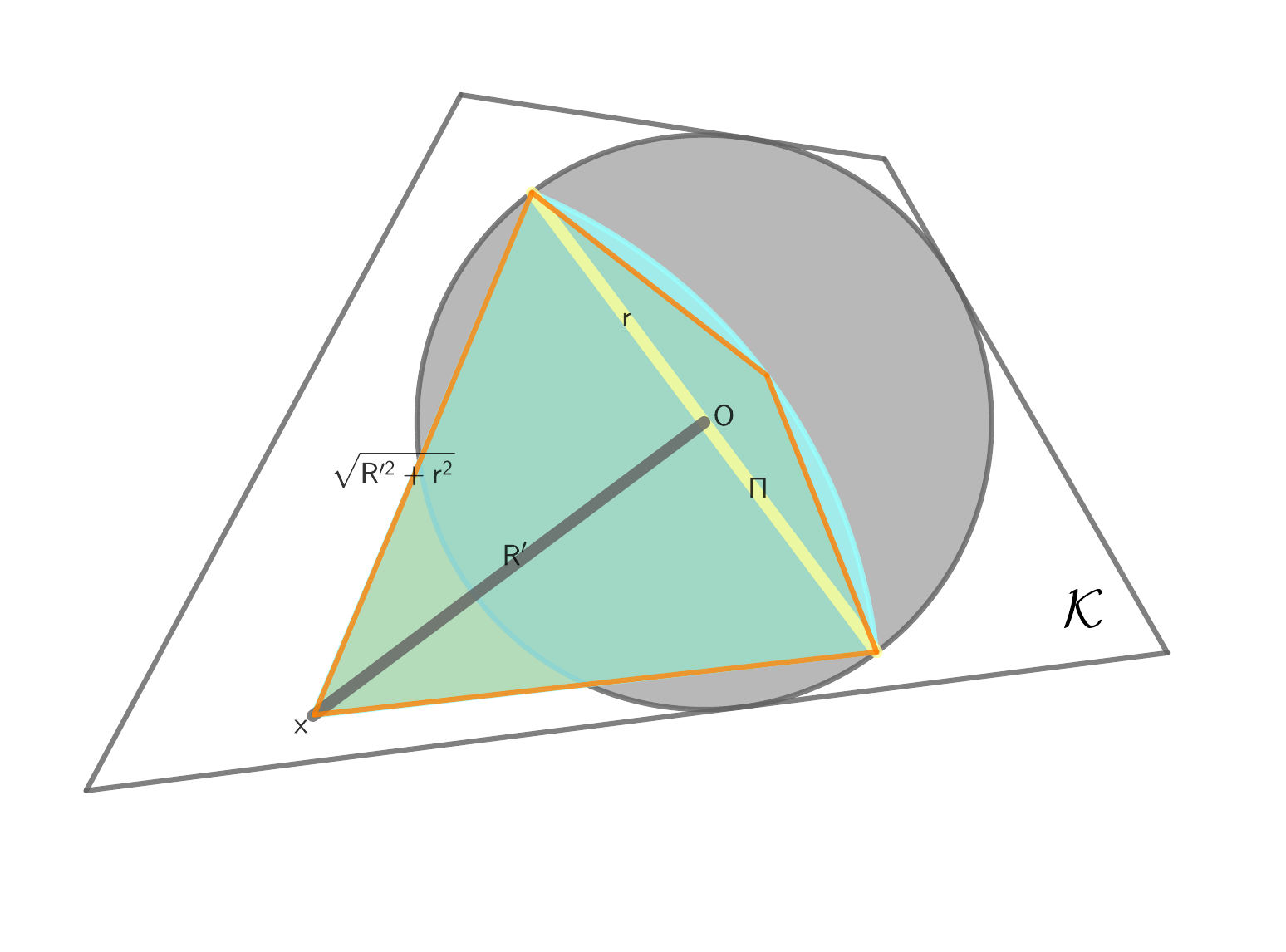}
\caption{This picture shows the case when $d=2$ in Lemma \ref{proportion}. The yellow segment is $\Pi$ and the blue sector is the spherical cone we are concerned with. The orange polygon is the combination of the two non-intersecting cones with the same bases $\Pi$.
}
\label{figure 1}
\end{figure}

\begin{proof}
    It's easy to see the result is true when $d=1$. For $d\ge 2$ and any non-zero $x\in\mathcal{K}$, let $\|x\|_2=R'\leq R$ and name the cross-section of $\mathbb{B}_r^d$, which is perpendicular to the line going through $x$ and $0$, $\Pi$ (see Figure \ref{figure 1}). Then by the convexity of $\mathcal{K}$, $\mathcal{K}$ contains the spherical cone induced by $\Pi$ in $\mathbb{B}_{\sqrt{R'^2+r^2}}^d(x)$. Clearly, this spherical cone consists of two non-intersecting cones with the same bases $\Pi$ and the sum of their heights is $\sqrt{R'^2+r^2}$. Therefore, its volume should be larger than
    \[
    V_1:=\frac{1}{d}\sqrt{R'^2+r^2}\operatorname{Vol}_{d-1}(\Pi)=\frac{\pi^{\frac{d-1}{2}}r^{d-1}\sqrt{R'^2+r^2}}{d\Gamma\left(\frac{d+1}{2}\right)},
    \]
    where we used Lemma \ref{vol of ball} because $\Pi$ is a $d-1$-dimensional ball with radius $r$. On the other hand, the volume of $\mathbb{B}_{\sqrt{R'^2+r^2}}^d(x)$ is
    \[
    V_2:=\frac{\pi^{d / 2}(R'^2+r^2)^{d/2}}{\Gamma\left(\frac{d}{2}+1\right)}.
    \]
    Hence, the proportion of the volume of this spherical cone in $\mathbb{B}_{\sqrt{R'^2+r^2}}^d(x)$ is larger than
    \[
    \frac{V_1}{V_2}=\pi^{-1/2}d^{-1}\cdot\frac{\Gamma\left(\frac{d}{2}+1\right)}{\Gamma\left(\frac{d+1}{2}\right)}\cdot\frac{r^{d-1}}{(R'^2+r^2)^{(d-1)/2}}.
    \]
    By Lemma \ref{gautschi}, we have 
    \[
    \frac{\Gamma\left(\frac{d}{2}+1\right)}{\Gamma\left(\frac{d+1}{2}\right)}\ge 
    \sqrt{\frac{d}{2}}.
    \]
    Since $R',r\leq R$,
    \[
    \frac{r^{d-1}}{(R'^2+r^2)^{(d-1)/2}}\ge \frac{r^{d-1}}{(2R^2)^{(d-1)/2}}=\left(\frac{r}{\sqrt{2}R}\right)^{d-1},
    \]
    then finally, we have
    \[
    \frac{V_1}{V_2}\ge \frac{1}{\sqrt{2\pi d}}\left(\frac{r}{\sqrt{2}R}\right)^{d-1}.
    \]
    By similarity, one can see that this is also true for all $\mathbb{B}_{r'}^d(x)$ if $r'\leq r\leq\sqrt{R'^2+r^2}$.
\end{proof}
\subsection{Proof for Lemma \ref{growth rate >1}}\label{proof for growth rate >1}
Recall that in our algorithm, for all $t\leq\tau$, $\Sigma_{t}^{-1}\succeq\frac12\bar{\Sigma}^{-1}_t\succeq\frac{1}{2\sigma^2}\idm_d$ and $\|\mu_t-x_\star\|_{\Sigma_t^{-1}}^2\leq\frac{1}{\lambda^2L^2}.$
Then combining Eq.~\eqref{widetilde k} and Eq.~\eqref{xtkc}, we have
\[
\E_{t-1}[H_t]\succeq\lambda\beta \bc(\bc-1)\frac{1}{\sqrt{2\pi d}}\left(\frac{r}{\sqrt{2}R}\right)^{d-1}\cdot\|\Sigma_t\|^{\frac{\bc-2}{2}}\idm_d\E_{t-1}\left[\left(\|U\|_2+\frac{1}{\lambda L}\right)^{\bc-2}\cdot\1_{\{\|U\|_2\leq\frac{r}{\sqrt{2}(1-\lambda)\sigma}\}}\right].
\]
Since $\frac{r}{\sqrt2 \sigma}\ge 5d$, $\1_{\{\|U\|_2\leq\frac{r}{\sqrt{2}(1-\lambda)\sigma}\}}\ge \1_{\{\|U\|_2\leq 5d\}}.$ Noting that $\bc\leq 2$ and by Jensen's inequality, we have
\[
\E\left[\left(\|U\|_2+\frac{1}{\lambda L}\right)^{\bc-2}\cdot\1_{\{\|U\|_2\leq 5d\}}\right]\ge \P(\|U\|_2\leq 5d)\cdot\left(\frac{\E[\|U\|_2\cdot\1_{\{\|U\|_2\leq 5d\}}]}{\P(\|U\|_2\leq 5d)}+\frac{1}{\lambda L}\right)^{\bc-2}\ge \frac12\left(10d+\frac1{\lambda L}\right)^{\bc-2},
\]
where we used that by Lemma \ref{chi}, $\P(\|U\|_2\leq 5d)\ge 1-e^{-d}>1/2.$ Since $\lambda\leq \frac{1}{10dL}$ and $\bc\ge1$, we have
\[
\frac12\left(10d+\frac1{\lambda L}\right)^{\bc-2}\ge\frac12\left(\frac2{\lambda L}\right)^{\bc-2}=\frac{2^{\bc-3}}{\lambda^{\bc-2} L^{\bc-2}}\ge \frac{1}{4\lambda^{\bc-2} L^{\bc-2}}.
\]
Therefore,
\[
\E_{t-1}[H_t]\succeq\frac{\beta \bc(\bc-1)}{4\sqrt{2\pi d}}\left(\frac{r}{\sqrt{2}R}\right)^{d-1}\lambda^{3-\bc}L^{2-\bc}\cdot\|\Sigma_t\|^{\frac{\bc-2}{2}}\idm_d.
\]
Then
\[
\bar{\Sigma}_{t+1}-\bar{\Sigma}_t\succeq\frac{\eta}{2}\E_{t-1}[H_t]\succeq\frac{\beta \bc(\bc-1)}{8\sqrt{2\pi d}}\left(\frac{r}{\sqrt{2}R}\right)^{d-1}\eta\lambda^{3-\bc}L^{2-\bc}\cdot\|\Sigma_t\|^{\frac{\bc-2}{2}}\idm_d.
\]
Recalling that for all $t\leq\tau$, 
$
\|\Sigma_t\|\leq 2\|\bar{\Sigma}_t\|
$, we have
\[
\|\Sigma_t\|^{\frac{\bc-2}{2}}\ge 2^{\frac{\bc-2}{2}}\|\bar{\Sigma}_t\|^{\frac{\bc-2}{2}}\ge \|\bar{\Sigma}_t\|^{\frac{\bc-2}{2}}/\sqrt{2},
\]
hence, now Eq.~\eqref{sigma min appro} becomes
\[
\lambda_{\operatorname{min}}(\bar{\Sigma}^{-1}_{t+1})-\lambda_{\operatorname{min}}(\bar{\Sigma}^{-1}_{t})\ge
\frac{\beta(\bc-1)}{30\sqrt{d}}
\left(\frac{r}{\sqrt{2}R}\right)^{d-1}\eta\lambda^{3-\bc}L^{2-\bc}\cdot\lambda_{\operatorname{min}}(\bar{\Sigma}^{-1}_{t})^{\frac{2-\bc}{2}}:=\Theta^{\frac{\bc}2}\lambda_{\operatorname{min}}(\bar{\Sigma}^{-1}_{t})^{\frac{2-\bc}{2}},
\]
where we used that $\frac{\bc}{16\sqrt{\pi}}>\frac{1}{30}$ and denoted that $\Theta=\left(\frac{\bc-1}{30}\right)^{2/\bc}\beta^{\frac2\bc}d^{-\frac1\bc}\left(\frac{r}{\sqrt{2}R}\right)^{\frac{2(d-1)}{q}}\eta^{\frac2\bc}\lambda^{\frac6\bc-2}L^{\frac4\bc-2}$. Since $\lambda_{\operatorname{min}}(\bar{\Sigma}^{-1}_{1})=\sigma^{-2}$ and $\sigma^{-2}\ge \Theta$,  by Lemma \ref{growth rate lamma},
\[
\lambda_{\operatorname{min}}(\bar{\Sigma}^{-1}_{t})\ge
\Theta t^{\frac2\bc}/8.
\]
Then for all $t\leq\tau$, 
\[
\Sigma_t^{-1}\succeq\frac12\bar{\Sigma}_{t}^{-1}\succeq\frac{\Theta}{16}t^{\frac2\bc}.
\]

\subsection{Proof for Lemma \ref{growth rate =1 lemma}}\label{proof for growth rate =1 lemma}
Clearly,
    \[
    \begin{aligned}
        \E\left[\frac{\idm_d-\frac{(X-x_\star)(X-x_\star)^\top}{\|X-x_\star\|_2^2}}{\|X-x_\star\|_2}\right]=&(2\pi)^{-\frac{d}{2}}\operatorname{det}(\Sigma)^{-\frac12}\int_{\R^d}\frac{\idm_d-\frac{xx^\top}{\|x\|_2^2}}{\|x\|_2}e^{-\frac12\left(x-(\mu-x_\star)\right)^\top \Sigma^{-1}\left(x-(\mu-x_\star)\right)}\d x\\
        \succeq&(2\pi)^{-\frac{d}{2}}e^{-\|\mu-x_\star\|^2_{\Sigma^{-1}}}\operatorname{det}(\Sigma)^{-\frac12}\int_{\R^d}\frac{\idm_d-\frac{xx^\top}{\|x\|_2^2}}{\|x\|_2}e^{-x^\top \Sigma^{-1}x}\d x\\
        =& 2^{-\frac{d}{2}}e^{-\|\mu-x_\star\|^2_{\Sigma^{-1}}}\E\left[\frac{\idm_d-\frac{YY^{\top}}{\|Y\|_2^2}}{\|Y\|_2}\right]:=2^{-\frac{d}{2}}e^{-\|\mu-x_\star\|^2_{\Sigma^{-1}}}I,
    \end{aligned}
    \]
    where we used that $\|x-(\mu-x_\star)\|^2_{\Sigma^{-1}}\leq 2\|x\|^2_{\Sigma^{-1}}+2\|\mu-x_\star\|^2_{\Sigma^{-1}}$ and $Y\sim\N(0,\Sigma/2)$. By rotation invariance, W.L.O.G., one can assume that $\Sigma$ is diagonal, say $\diag{\sigma_1^2,\cdots,\sigma_d^2}$, where $\|\Sigma\|^{\frac12}=\sigma_1\ge \cdots\ge\sigma_d>0.$ For all $1\leq i\ne j\leq d$, by symmetry
    \[
    I_{ij}=\E\left[\frac{-Y_iY_j}{\|Y\|_2}\right]=0,
    \]
    which implies that $I$ is also diagonal. Then, we show that $I_{11}$ is the smallest diagonal element and hence $\lambda_{\operatorname{min}}(I)=I_{11}$. 
    Since $Y_1, \ldots, Y_d$ are mutually independent, 
    \[
    \begin{aligned}
         \E\left[
    \frac{Y_2^2-Y_1^2}{\|Y\|_2^3}
    \big{|} Y_3,\cdots,Y_d
    \right]&=\frac{1}{\pi\sigma_1\sigma_2}\int_{\R^2}\frac{y_2^2-y_1^2}{\|y\|_2^3}e^{-\frac{y_1^2}{\sigma_1^2}-\frac{y_2^2}{\sigma_2^2}}\d y_1\d y_2\\&=\frac{1}{2\pi\sigma_1\sigma_2}\int_{\R^2}\frac{(y_2^2-y_1^2)\left(e^{-\frac{y_1^2}{\sigma_1^2}-\frac{y_2^2}{\sigma_2^2}}-e^{-\frac{y_2^2}{\sigma_1^2}-\frac{y_1^2}{\sigma_2^2}}\right)}{\|y\|_2^3}\d y_1\d y_2,
    \end{aligned}
    \]
    where the second equality used symmetry. Since $\sigma_1\ge\sigma_2$, we have
    \[
    (y_2^2-y_1^2)\left(\left[-\frac{y_1^2}{\sigma_1^2}-\frac{y_2^2}{\sigma_2^2}\right]-\left[-\frac{y_2^2}{\sigma_1^2}-\frac{y_1^2}{\sigma_2^2}\right]\right)=(y_2^2-y_1^2)^2(\frac{1}{\sigma_1^2}-\frac{1}{\sigma_2^2})\leq 0,
    \]
    and by that $e^x$ is increasing, we also have
    \[
    (y_2^2-y_1^2)\left(e^{-\frac{y_1^2}{\sigma_1^2}-\frac{y_2^2}{\sigma_2^2}}-e^{-\frac{y_2^2}{\sigma_1^2}-\frac{y_1^2}{\sigma_2^2}}\right)\leq0,
    \]
    which implies that
    \[
    I_{11}-I_{22}=\E\left[
    \E\left[
    \frac{Y_2^2-Y_1^2}{\|Y\|_2^3}
    \big{|} Y_3,\cdots,Y_d
    \right]
    \right]\leq 0.
    \]
    Similarly, for all $i> 1$, $I_{11}\leq I_{ii}$. Therefore, it suffices to lower bound $I_{11}$. Let $S=Y_2^2+\cdots+Y_d^2$. Then
    \[
    \E[I_{11}\mid Y_2,\cdots,Y_d]=\E\left[\frac{1-\frac{Y_1^2}{\|Y\|_2^2}}{\|Y\|_2}\mid Y_2,\cdots,Y_d\right]=\E\left[\frac{S}{(Y_1^2+S)^{\frac32}}\mid S\right].
    \]
    Note that $Y_1\sim\N(0,\sigma_1^2/2)$, then this is larger than
    \[
    \int_{-\sqrt{S}}^{\sqrt{S}}\frac{1}{\sigma_1\sqrt{\pi}}\frac{S}{(y_1^2+S)^{\frac32}}e^{-\frac{y_1^2}{\sigma_1^2}}\d y_1\ge \int_{-\sqrt{S}}^{\sqrt{S}}\frac{1}{\sigma_1\sqrt{\pi}}\frac{S}{(S+S)^{\frac32}}e^{-\frac{S}{\sigma_1^2}}\d y_1\ge \frac{1}{\sigma_1\sqrt{2\pi}}e^{-\frac{S}{\sigma_1^2}}.
    \]
    Hence,
    \[
    \lambda_{\operatorname{min}}(I)=I_{11}\ge \frac{1}{\sigma_1\sqrt{2\pi}}\E[e^{-\frac{S}{\sigma_1^2}}]= \frac{1}{\sigma_1\sqrt{2\pi}}\prod_{i=2}^d\E[e^{-\frac{Y_i^2}{\sigma_1^2}}]=\frac{1}{\sigma_1\sqrt{2\pi}}\prod_{i=2}^d(1+\frac{2\sigma_i^2}{\sigma_1^2})^{-\frac12}\ge \frac{3^{-\frac{d-1}{2}}}{\sigma_1\sqrt{2\pi}},
    \]
    where we used that $Y_2,\cdots,Y_d$ are mutually independent and $\sigma_i\leq\sigma_1.$ In summary, we have
    \[
    \E\left[\frac{\idm_d-\frac{(X-x_\star)(X-x_\star)^\top}{\|X-x_\star\|_2^2}}{\|X-x_\star\|_2}\right]\succeq2^{-\frac{d}{2}}e^{-\|\mu-x_\star\|^2_{\Sigma^{-1}}}\frac{3^{-\frac{d-1}{2}}}{\sqrt{2\pi}}\|\Sigma\|^{-\frac12}\idm_d\succeq 6^{-\frac{d}2}e^{-\|\mu-x_\star\|^2_{\Sigma^{-1}}}\|\Sigma\|^{-\frac12}\idm_d/2,
    \]
    where we used that $\sqrt{\frac{3}{2\pi}}>\frac12.$

\subsection{Proof for Lemma \ref{growth rate=1}}\label{proof for growth rate=1}
\begin{proof}
    The proof is similar to that of Lemma \ref{growth rate >1}. Recall that for all $t\leq\tau$, $\Sigma_{t}^{-1}\succeq\frac12\bar{\Sigma}^{-1}_t\succeq\frac{1}{2\sigma^2}\idm_d$ and $\|\mu_t-x_\star\|_{\Sigma_t^{-1}}^2\leq\frac{1}{\lambda^2L^2}.$ Note that $\bar{\Sigma}_t\preceq\sigma^2\idm_d$, for all $t\leq \tau$, $\|\bar{\Sigma}_t\|\leq\sigma^2\leq 1$ and that $-\frac12\leq\frac{\kappa-1}{2-\kappa}\leq 0$. Then we have
    \[
    \|\bar{\Sigma}_t\|^{-\frac12}\ge\|\bar{\Sigma}_t\|^{\frac{\kappa-1}{2-\kappa}}.
    \]
    Hence by Eq.~\eqref{all d} and  $\|\Sigma_t\|\leq2\|\bar{\Sigma}_t\|,$ we have
    \[
    \E_{t-1}[H_t]\succeq\beta 6^{-\frac{d}2}e^{-\|\mu_t-x_\star\|_{\Sigma_t^{-1}}^2}\|\Sigma_t\|^{-\frac12}\idm_d/2\succeq\beta 6^{-\frac{d}2}e^{-\|\mu_t-x_\star\|_{\Sigma_t^{-1}}^2}\|\bar{\Sigma}_t\|^{\frac{\kappa-1}{2-\kappa}}\idm_d/2\sqrt{2},
    \]
    and hence, 
    \[
    \bar{\Sigma}_{t+1}-\bar{\Sigma}_t\succeq\frac{\eta}{2}\E_{t-1}[H_t]\succeq\beta\eta 6^{-\frac{d}2}e^{-\frac{1}{\lambda
    ^2L^2}}\|\bar{\Sigma}_t\|^{\frac{\kappa-1}{2-\kappa}}\idm_d/4\sqrt{2},
    \]
    which implies that 
    \[
    \lambda_{\operatorname{min}}(\bar{\Sigma}^{-1}_{t+1})-\lambda_{\operatorname{min}}(\bar{\Sigma}^{-1}_{t})\ge\beta\eta 6^{-\frac{d}2}e^{-\frac{1}{\lambda
    ^2L^2}}\lambda_{\operatorname{min}}(\bar{\Sigma}^{-1}_{t})^{\frac{1-\kappa}{2-\kappa}}/4\sqrt{2}\ge\Theta^{\frac{1}{2-\kappa}}\lambda_{\operatorname{min}}(\bar{\Sigma}^{-1}_{t})^{\frac{1-\kappa}{2-\kappa}},
    \]
    where we denoted  $\Theta=\beta^{2-\kappa}\eta^{2-\kappa}6^{-\frac{d(2-\kappa)}{2}}e^{-\frac{2-\kappa}{\lambda^2L^2}}/32$ and used  $(4\sqrt{2})^{2-\kappa}\leq32$.
    
    Finally, since $\Theta\leq \sigma^{-2}=\lambda_{\operatorname{min}}(\bar{\Sigma}_1^{-1})$ and by Lemma \ref{growth rate lamma}, we have
    \[
    \lambda_{\operatorname{min}}(\bar{\Sigma}_t^{-1})\ge\Theta t^{2-\kappa}/8,
    \]
    and hence
    \[
    \Sigma_t\succeq\frac{1}{2}\bar{\Sigma}_t\succeq\frac{\Theta}{16} t^{2-\kappa}.
    \]
\end{proof}

\subsection{Proof for Lemma \ref{sk e1}}\label{proof for sk e1}
\begin{proof}
        By Eq.~\eqref{extension noise}, $\frac{|\xi_t|}{\|X_t-x_\star\|_2}$ is $2R/r$-subgaussian conditioning on $\mathscr{F}_{t-1}$ and $X_t$. Note that
        \[
        \max_{1\leq t\leq\tau}\frac{|\xi_t|}{\|X_t-x_\star\|_2}\leq \max_{1\leq t\leq n}\frac{|\xi_t|}{\|X_t-x_\star\|_2}.
        \]
        By Lemma \ref{k1k2}, there exists $C>0$ such that 
        \[
        \P(\operatorname{E}_1^c)\leq\sum_{t=1}^n 2\exp\left(-C\left[\frac{(RL^{1/2})}{r}\right]^2/(R/r)^2\right),
        \]
        which is less than $\delta/5$ when $L$ is large enough.

        Since $X_t|\mathscr{F}_{t-1}\sim\N(\mu_t,\Sigma_t)$,  $\|X_t-\mu_t\|_{\Sigma_t^{-1}}^2|\mathscr{F}_{t-1}\sim\chi^2_d.$ By Lemma \ref{norm subgaussian}, there exists $C'>0$ such that for $t>0$,
        \[
        \P(\|X_t-\mu_t\|_{\Sigma_t^{-1}}-\sqrt{d}\ge t)\leq 2\exp(-t^2/C').
        \]
        Hence
        \[
        \P\left(
        \max_{1\leq t\leq n}\|X_t-\mu_t\|_{\Sigma_t^{-1}}\ge \sqrt{d}+\sqrt{C'\log(\frac{10n}{\delta})}
        \right)\leq \delta/5,
        \]
        and similarly we obtain $\P(\operatorname{E}_2)\ge 1-\delta/5.$
    \end{proof}
\subsection{Proof for Lemma \ref{sk Z V}}\label{proof for sk Z V}
\begin{proof}
    First, on $\operatorname{E}_1\cap\operatorname{E}_2$, similar to Eq.~\eqref{x xstar l2}, for all $t\leq \tau$, we have $\|\mu_{t}-x_\star\|_2\leq \frac{2}{L\lambda\sqrt{\eta\gamma t}}$. Since $\|X_t-x_{\star}\|_2\leq \|X_t-\mu_t\|_2+\|\mu_t-x_{\star}\|_2$ and $\|\mu_t-x_{\star}\|_2$ is $\mathscr{F}_{t-1}$-measurable, we have
    \[
    \E_{t-1}[\|X_{t}-x_\star\|_2^2]\leq 2\E_{t-1}[\|X_t-\mu_t\|_2^2]+2\|\mu_{t}-x_\star\|^2_2\leq 2\E_{t-1}[\|X_t-\mu_t\|_2^2]+\frac{8}{L^2\lambda^2\eta\gamma t}.
    \]
    Conditioning on $\mathscr{F}_{t-1}$, $X_t\sim\N(\mu_t,\Sigma_t)$, then we obtain
    \[ 
    \E_{t-1}[\|X_t-\mu_t\|_2^2]=\tr{\Sigma_t}\leq \frac{2d}{\eta\gamma t}\leq\frac{1}{L^2\lambda^2\eta\gamma t},
    \]
    where the second inequality used that $\Sigma^{-1}_t\succeq\frac{\eta\gamma t}{2}\idm_d$ and the final inequality holds because $\lambda\leq \frac{1}{2\sqrt{d}L}.$ Then $\E_{t-1}[\|X_{t}-x_\star\|^2_2]\leq\frac{10}{L^2\lambda^2\eta\gamma t},
    $ 
    and for all $t\ge 2$, there exist $C,C'>0$ such that
    \begin{equation}\label{tge2}
            \begin{aligned}
        \E_{t-1}[|Z_t|^2]=&\E_{t-1}\left[3\lip(e)^2\cdot\|X_t-x_\star\|^2_2+3|\xi_t|^2\right]+3\left(\lip(e)\cdot\|X_{t-1}-x_\star\|_2+|\xi_{t-1}|\right)^2\\
        \leq &\left(3\lip(e)^2+\frac{12CR^2}{r^2}\right)\E_{t-1}[\|X_{t}-x_\star\|_2^2]+\left(\lip(e)+\frac{RL^{1/2}}{r}\right)^2\cdot\frac{12}{L^2\lambda^2\eta\gamma t}\\
        \leq& \left(\lip(e)+\frac{R}{r}\right)^2\cdot\frac{C'}{L\lambda^2\eta\gamma t},
    \end{aligned}
    \end{equation}
    where the first inequality used the definition of $\operatorname{E}_1$ and that by Lemma \ref{k1k2},
    $
    \E_{t-1}[|\xi_t|^2]\leq\E_{t-1}[\frac{2CR^2}{r^2}\|X_t-x_{\star}\|_2^2]
    $, because $\xi_t$ is conditionally $\frac{2R}{r}\|X_t-x_{\star}\|_2$-subgaussian. For $t=1$, recall Eq.~\eqref{t=1}. 
    Then it's easy to check that Eq.~\eqref{tge2} is also true.
    Also, by the Cauchy-Schwarz inequality, $\E_{t-1}[|Z_t|]\leq\left(\lip(e)+\frac{R}{r}\right)\cdot\frac{\sqrt{C'}}{\sqrt{L\lambda^2\eta\gamma t}}.$
    Therefore, combining Eq.~\eqref{Zt}, for all $t\leq \tau$,
    \[
    \left|Z_t\right|+\mathbb{E}_{t-1}\left[\left|Z_t\right|\right]\leq\frac{H}{3\sqrt{L\lambda^2\eta\gamma t}},
    \]
    we have that $Z_{\max}\leq \frac{H}{3\sqrt{L\lambda^2\eta\gamma}}$ and 
    \[
    V_\tau\leq\sum_{t=1}^\tau\left(\lip(e)+\frac{R}{r}\right)^2\cdot\frac{C'}{L\lambda^2\eta\gamma t}\leq \frac{H^2}{9\lambda^2\eta\gamma},
    \]
    by the definition of $H$.
\end{proof}

\subsection{Proof for Lemma \ref{e4}}\label{proof for e4}
\begin{proof}
    Our plan is to apply Lemma \ref{concentration}. Let $\Delta_t=f\left(\frac{X_t}{\pi^+(X_t)}\right)-f(x_\star)$. Then by definitions,
    \[
    \begin{aligned}
        \operatorname{Reg}_\tau(x_\star)-\widetilde{\Reg}_\tau(x_\star)=&\sum_{t=1}^{\tau}\left(
    f\left(\frac{X_t}{\pi^+(X_t)}\right)-\E_{t-1}\left[f\left(\frac{X_t}{\pi^+(X_t)}\right)\right]
    \right)\\
    =&\sum_{t=1}^{\tau}\left(
    \Delta_t-\E_{t-1}\left[\Delta_t\right]
    \right).
    \end{aligned}
    \]
    By Lemma \ref{extension}, it's clear that
    \[
    0\leq \Delta_t\leq e(X_t)-f(x_\star)=e(X_t)-e(x_\star)\leq\lip(e)\|X_t-x_\star\|_2.
    \]
    Hence by the proof of Lemma \ref{sk Z V}, on $\operatorname{E}_1\cap\operatorname{E}_2$, we have
    \[
    \widetilde{Z}_{\max}:=\max _{1 \leq t \leq \tau}\left(\left|\Delta_t\right|+\mathbb{E}_{t-1}\left[\left|\Delta_t\right|\right]\right)\leq \frac{H}{2\sqrt{L\lambda^2\eta\gamma}}, \widetilde{V}_\tau:=\sum_{t=1}^\tau \E_{t-1}[\Delta^2_t]\leq \frac{H^2}{4\lambda^2\eta\gamma}.
    \]
    \noindent
    By lemma \ref{concentration} and noting that $\frac{H}{\lambda\sqrt{L\eta\gamma}}\ge 2$, with probability at least $1-\delta/5$,
    \[
    \sum_{t=1}^{\tau}\left(
    \Delta_t-\E_{t-1}\left[\Delta_t\right]
    \right)\leq \sqrt{\widetilde{V}_{\tau}L}+\widetilde{Z}_{\max } L.
    \]
    \noindent
    Then when $\operatorname{E}_1$ and $\operatorname{E}_2$ both happen, 
    \[
    \operatorname{Reg}_\tau(x_\star)-\widetilde{\Reg}_\tau(x_\star) \leq H\sqrt{\frac{L}{\lambda^2\eta\gamma}}.
    \]
\end{proof}

%% file: tex/append/exchange.tex
In this section, we generalize Stein's Lemma \citep{stein1981estimation} to Lipschitz convex functions. Adapting the proof of Lemma 1.1 from \cite{demaret2019stein} and recalling that, by Monotone Differentiation Theorem, all monotone functions are almost everywhere differentiable, we have:
\begin{Lem}[modification of Stein's Lemma]\label{xgx}
    If $g(x):\R\to\R$ is monotonically increasing in $\R$, and $X\sim\N(0,1)$, then
    \[
    \E[g^\prime(X)]\leq \E[Xg(X)].
    \]
    If $g(x)$ is absolutely continuous (not necessarily monotone), the equality holds.
    \begin{proof}
    By Monotone differentiation theorem, 
    \begin{equation}\label{gy-g0}
        \int_0^y g^{\prime}(x) \mathrm{d} x\leq g(y)-g(0), \forall y\ge 0, \quad\int_y^0 g^{\prime}(x) \mathrm{d} x\leq g(0)-g(y), \forall y\leq 0.
    \end{equation}
        The density function $\phi(x)$ of the standard Gaussian law, as it will be noted from now on, is such that $\phi^{\prime}(x)=-x \phi(x)$. Note also how, using $\int_{\mathbb{R}} y \phi(y) \mathrm{d} y=0$,

$$
\int_{-\infty}^x-y \phi(y) \mathrm{d} y=\int_x^{+\infty} y \phi(y) \mathrm{d} y, \quad \forall x \in \mathbb{R}.
$$
\noindent
We then have

$$
\begin{aligned}
\mathbb{E}\left[g^{\prime}(X)\right] & =\int_{\mathbb{R}} g^{\prime}(x) \phi(x) \mathrm{d} x \\
& =\int_{\mathbb{R}} g^{\prime}(x)\left(\int_{-\infty}^x(-y \phi(y)) \mathrm{d} y\right) \mathrm{d} x \\
& =\int_0^{+\infty} g^{\prime}(x)\left(\int_x^{+\infty} y \phi(y) \mathrm{d} y\right) \mathrm{d} x-\int_{-\infty}^0 g^{\prime}(x)\left(\int_{-\infty}^x y \phi(y) \mathrm{d} y\right) \mathrm{d} x \\
& =\int_0^{+\infty} y \phi(y)\left(\int_0^y g^{\prime}(x) \mathrm{d} x\right) \d y-\int_{-\infty}^0 y \phi(y)\left(\int_y^0 g^{\prime}(x) \mathrm{d} x\right) \mathrm{d} y \\
& \leq\int_{\mathbb{R}} y \phi(y)(g(y)-g(0)) \mathrm{d} y \\
& =\mathbb{E}[X g(X)]-g(0) \mathbb{E}[X]=\mathbb{E}[X g(X)],
\end{aligned}
$$
where the inequality used Eq.~\eqref{gy-g0}.
    \end{proof}
\end{Lem}
\begin{Lem}\label{dirac}
    If $X\sim\N(0,1)$ and $g(x)=\1_{\{x\ge a\}}$, where $a\in\R$, then
    \[
    \E[Xg(X)]=\frac{1}{\sqrt{2\pi}}e^{-\frac{a^2}{2}}.
    \]
    \begin{proof}
       Clearly, 
        \[
        \E[Xg(X)]=\int_{a}^{+\infty}\frac{x}{\sqrt{2\pi}}e^{-\frac{x^2}{2}}\d x=-\frac{1}{\sqrt{2\pi}}e^{-\frac{x^2}{2}}\big{|}_{a}^{+\infty}=\frac{1}{\sqrt{2\pi}}e^{-\frac{a^2}{2}}.
        \]
    \end{proof}
\end{Lem}
\begin{Thm}[Alexandrov]\label{ale}
    If $f(x)$ is a convex function over $U$, which is an open subset of $\R^d$, then $f(x)$ has a second derivative almost everywhere.
\end{Thm}
\begin{Lem}\label{gx^2}
    If $g(x):\R\to\R$ is convex and absolutely continuous on $\R$, $X\sim\N(0,1)$ and $\E[|g(X)|(X^2+1)]<+\infty$, then
    \[
    \E[g^{\prime\prime}(X)]\leq\E[Xg'(X)]=\E[g(X)(X^2-1)].
    \]
    \begin{proof}
        Since $g\prime$ is monotonically increasing, then by Lemma \ref{xgx}, we have
        \[
    \E[g^{\prime\prime}(X)]\leq \E[Xg^\prime(X)].
    \]
    Note that $xg(x)$ is also absolutely continuous and $(xg(x))'=xg'(x)+g(x)$, then similarly, by Lemma \ref{xgx}, we have
    \[
    \E[Xg^\prime(X)]=\E[(Xg(X))^\prime]-\E[g(X)]=\E[g(X)(X^2-1)].
    \]
    \end{proof}
\end{Lem}
\begin{Lem}[Generalized Stein’s Lemma]\label{convex lip}
    If $g(x):\R^d\to\R$ is convex and Lipschitz continuous on $\R^d$ and $X\sim\N(\mu,\Sigma)$, then
    \[
    \E[\nabla^2 g(X)]\preceq\E\left[g(X)\left\{\Sigma^{-1}(X-\mu)(X-\mu)^{\top}\Sigma^{-1}-\Sigma^{-1}\right\}\right].
    \]
    \begin{proof}
        Let $Z=\Sigma^{-1/2}(X-\mu)\sim\N(0,\idm)$ and $h(z)=g\left(\Sigma^{1/2}z+\mu\right)$. Then it's clear that
        \[
        \nabla^{2} h(z)=\Sigma^{1/2}\nabla^2 g\left(\Sigma^{1/2}z+\mu\right)\Sigma^{1/2},
        \]
        hence
        \[
        \E[\nabla^2 g(X)]=\E\left[\nabla^2 g\left(\Sigma^{1/2}Z+\mu\right)]=\Sigma^{-1/2}\E[\nabla^2 h(X)\right]\Sigma^{-1/2}.
        \]
        And the right hand is equal to
        \[
        \Sigma^{-1/2}\E\left[g\left(\Sigma^{1/2}Z+\mu\right)\cdot\left\{ZZ^{\top}-\idm\right\}\right]\Sigma^{-1/2}=\Sigma^{-1/2}\E\left[h(Z)\cdot\left\{ZZ^{\top}-\idm\right\}\right]\Sigma^{-1/2}.
        \]
        
        Hence it suffices to show that
        \[
        \E[\nabla^2 h(Z)]\preceq\E[h(Z)(ZZ^{\top}-\idm)], 
        \]
        namely for all $u\in\R^d$ such that $\|u\|_2=1$, we have
        \[
        \E[u^\top\nabla^2 h(Z)u]\leq\E[h(Z)(u^\top ZZ^{\top}u-1)].
        \]
        Take an orthogonal matrix $P$ such that its first column is just $u$ and let $r(w)=h(Pw)$. Then $\nabla^2 r(w)=P^{\top}\nabla^2 h(Pw) P$. Let $W=P^{\top}Z\sim\N(0,\idm)$, then the left hand becomes 
        \[
        \E[\nabla^2 r(W)]_{11}=\E\left[\frac{\pd^2}{\pd w_1\pd w_1} r(W)\right],
        \]
        and the right hand becomes
        \[
        \E[r(W)(W_1^2-1)].
        \]
        Finally, note that $r(w)$ is also Lipschitz continuous and hence is absolutely continuous about $w_1$ for any fixed $w_2,\cdots, w_d$. Then by Lemma \ref{gx^2}, we have
        \[
        \E\left[\frac{\pd^2}{\pd w_1\pd w_1} r(W)\big{|} W_2,\cdots,W_d\right]\leq \E[r(W)(W_1^2-1)\mid W_2,\cdots,W_d],
        \]
        then the result follows by taking expectation for both sides.
    \end{proof}
\end{Lem}

%% file: tex/append/hessian.tex
\begin{Lem}\label{nabla2}
    For all $2\ge p>1$, $\bc>1$ and $x\ne 0$, we have
    \[
    \nabla^2\|x\|_p^\bc\succeq \bc(\bc\wedge p -1)d^{-\frac{(2-\bc)(2-p)}{2p}} \|x\|_{2}^{\bc-2}\idm_d.
    \]
    For all $2\ge p>1$ and $x\ne 0$, we have
    \[
    \nabla^2\|x\|_p\succeq(p-1) d^{-\frac{3(2-p)}{2p}}\|x\|_2^{-1}\cdot\left(\idm_d-\frac{xx^{\top}}{\|x\|_2^2}\right).
    \]
    \begin{proof}
        By direct computation, we have
        \[
        \nabla \|x\|_p=\|x\|_p^{1-p}\left(
        |x_1|^{p-1}\sgn{x_1},\cdots,|x_d|^{p-1}\sgn{x_d}
        \right)^\top:=
        x_{(p)},
        \]
        then
        \[
        \begin{aligned}
            \nabla^2 \|x\|_p=&(1-p)\|x\|_p^{-1}x_{(p)}x_{(p)}^\top+(p-1)\|x\|_p^{1-p}\diag{|x_1|^{p-2},\cdots,|x_d|^{p-2}}\\
            =&(p-1)\|x\|_p^{-1}\cdot\left(
            \Lambda_{(p)}-x_{(p)}x_{(p)}^\top
            \right),
        \end{aligned}
        \]
        where we denoted that $\Lambda_{(p)}=\diag{\left(\frac{|x_1|}{\|x\|_p}\right)^{p-2},\cdots,\left(\frac{|x_d|}{\|x\|_p}\right)^{p-2}}.$ Then by Lemma \ref{xp} and Lemma \ref{eui norm}, we have
        \[
        \nabla^2 \|x\|_p\succeq(p-1)d^{\frac12-\frac1p}\|x\|_2^{-1}\cdot\left(
            \Lambda_{(p)}-x_{(p)}x_{(p)}^\top
            \right)\succeq (p-1)d^{-\frac{3(2-p)}{2p}}\|x\|_2^{-1}\cdot\left(\idm_d-\frac{xx^{\top}}{\|x\|_2^2}\right).
        \]
        For all $\bc\ge1$, we have
        \[
        \begin{aligned}
            \nabla^2\|x\|_p^\bc=&\bc\left(
        (\bc-1)\|x\|_p^{\bc-2}\nabla\|x\|_p\nabla^\top\|x\|_p+\|x\|_p^{\bc-1}\nabla^2\|x\|_p
        \right)\\
        =&\bc\|x\|_p^{\bc-2}\cdot\left(
        (p-1)\Lambda_{(p)}+(\bc-p)x_{(p)}x_{(p)}^\top
        \right).
        \end{aligned}
        \]
        Therefore, note that by Lemma \ref{xp}, $\Lambda_{(p)}-x_{(p)}x_{(p)}^\top\succeq 0$, then we have
        \[
        \nabla^2\|x\|_p^\bc\succeq \bc(\bc\wedge p -1)\|x\|_{p}^{\bc-2}\Lambda_{(p)}\succeq \bc(\bc\wedge p -1)d^{-\frac{(2-\bc)(2-p)}{2p}} \|x\|_{2}^{\bc-2}\idm_d,
        \]
        where the last inequality used that $\Lambda_{(p)}\succeq\idm_d$ 
when $p\leq2$ and Lemma \ref{eui norm}.
    \end{proof}
\end{Lem}
\begin{Lem}\label{xp}
   Under the same definitions of $\Lambda_{(p)}$ and $x_{(p)}$ in Lemma \ref{nabla2}, then for all $x\ne 0$ and $1<p\leq 2$, we have
    \[
    \Lambda_{(p)}-x_{(p)}x_{(p)}^\top\succeq d^{-\frac{2-p}{p}}\left(\idm_d-\frac{xx^{\top}}{\|x\|_2^2}\right),
    \]
    which is clearly positive semidefinite.
    \begin{proof}
        For all $u\in\R^d$, we have
        \[
        \begin{aligned}
          u^\top \left(
            \Lambda_{(p)}-x_{(p)}x_{(p)}^\top
            \right)u
        =\sum_{i=1}^d u_i^2\cdot\left(\frac{|x_i|}{\|x\|_p}\right)^{p-2}-\left(\sum_{i=1}^d u_i\cdot \left(\frac{|x_i|}{\|x\|_p}\right)^{p-1}\sgn{x_i}\right)^2.
        \end{aligned}
        \]
        Choosing $a_i=u_i\cdot\left(\frac{|x_i|}{\|x\|_p}\right)^{\frac{p-2}{2}}$ and $b_i=\left(\frac{|x_i|}{\|x\|_p}\right)^{\frac{p}{2}}\sgn{x_i}$ and noting that $\sum_{i=1}^d b_i^2=1$, by Lemma \ref{cauchy binet}, the difference is just 
        \[
        \sum_{i,j=1}^n\left(a_ib_j-a_jb_i\right)^2=\frac{1}{\|x\|_p^2}\sum_{i,j=1}^n\left(\frac{|x_i|}{\|x\|_p}\right)^{p-2}\left(\frac{|x_j|}{\|x\|_p}\right)^{p-2}(u_ix_j-u_jx_i)^2,
        \]
        where we used that $|x_i|\sgn{x_i}=x_i$. By Lemma \ref{eui norm} and noting that $p> 1$, this is larger than
        \[
        \frac{d^{-\frac{
        2-p
        }{p}}}{\|x\|_2^2}\sum_{i,j=1}^n\left(\frac{|x_i|}{\|x\|_p}\right)^{p-2}\left(\frac{|x_j|}{\|x\|_p}\right)^{p-2}(u_ix_j-u_jx_i)^2.
        \]
        Noting that $\frac{|x_j|}{\|x\|_p},\frac{|x_j|}{\|x\|_p}\leq 1$ and $p-2\leq 0$, this is larger than
        \[
        \frac{d^{-\frac{2-p}{p}}}{\|x\|_2^2}\sum_{i,j=1}^n(u_ix_j-u_jx_i)^2=d^{-\frac{2-p}{p}}\cdot u^{\top}\left(\idm_d-\frac{xx^{\top}}{\|x\|_2^2}\right)u,
        \]
        which just implies
        \[
        \Lambda_{(p)}-x_{(p)}x_{(p)}^\top\succeq d^{-\frac{2-p}{p}}\left(\idm_d-\frac{xx^{\top}}{\|x\|_2^2}\right).
        \]
    \end{proof}
\end{Lem}

%% file: tex/append/Surrogatefunction.tex
In this section, we present some useful results for the surrogate loss functions. Most of them can be found in \cite{lattimore40bandit}.
\subsection{Preliminary}
In the following, we suppose that $s(x)$ is the surrogate function with $f(x)$, $\N(\mu,\Sigma)$ and $\lambda\in(0,1)$,
\begin{Lem}[Lemma 12.3(b) in \cite{lattimore40bandit}]\label{sleqe}
    For all $x\in\R^d, s(x)\leq f(x).$
\end{Lem}
\begin{Lem}[Proposition 12.5 in \cite{lattimore40bandit}]\label{upper bound for s''}
    For all $z\in\R^d$ and $t\leq n$, 
    \[
    \begin{aligned}
        &(a) \left\|\nabla^2 s(z)\right\| \leq \frac{\lambda \operatorname{lip}(f)}{1-\lambda} \sqrt{d\left\|\Sigma^{-1}\right\|};
&&(b)\left\|\Sigma^{1 / 2} \nabla^2 s(z) \Sigma^{1 / 2}\right\| \leq \frac{\lambda \operatorname{lip}(f)}{1-\lambda} \sqrt{d\|\Sigma\|}.
    \end{aligned}
    \]
\end{Lem}
\noindent
The lemma below is only used in the proof for Theorem \ref{q=1 result}.
\begin{Lem}\label{smu>sx}
    If $f(x)$ 
    satisfies Assumption 1 and 3 in Section \ref{main result q=1}, then we have
    \[
    s(\mu)\ge s(x_\star).
    \]
    \begin{proof}
        By the definition of $s(x)$, we have
        \[
        s(\mu)-s(x_\star)=\frac{1}{\lambda}\E\left[
        f((1-\lambda) X+\lambda\mu)-f((1-\lambda) X+\lambda x_\star)
        \right],
        \]
        where $X\sim\N(\mu,\Sigma)$. Let $Y=(1-\lambda)(X-\mu)$ and $g(x)=f(x+x_\star)$, then it's clear that
        \[
        s(\mu)-s(x_\star)=\frac{1}{\lambda}\E\left[
        g(Y+\mu-x_\star)-g(Y+(1-\lambda)(\mu-x_\star))
        \right].
        \]
        For all $v\in\R^d$, let $u(t)=\E[g(Y+t\cdot v)]$, then it's easy to see that $u(t)$ is convex on $\R$. Since $g(x)$ is Lipschitz continuous on $\R^d$, then we have
        \[
        \frac{\d u}{\d t}\mid_{t=0}=\E[\nabla g(Y+tv)\cdot v]\mid_{t=0}=\E[\nabla g(Y)]\cdot v=0,
        \]
        where the final inequality follows from Assumption 3 and the symmetry of $Y$. Hence for all $\lambda\in[0,1]$, we have $\E[g(Y+\mu-x_\star)]-\E[g(Y+(1-\lambda)(\mu-x_\star))]\ge 0$, which implies that $s(\mu)\ge s(x_\star).$
    \end{proof}
\end{Lem}
\noindent
Here and later, $s_t$ and $q_t$ are those defined in our algorithms. 
\begin{Lem}\label{unbias}
    The following hold:
    \[
    \begin{aligned}
(a) \mathbb{E}_{t-1}\left[\nabla\hat{s}_t(z)\right]=\nabla s_t(z);&&&&
(b) \mathbb{E}_{t-1}\left[\nabla^2\hat{s}_t(z)\right]=\nabla^2 s_t(z).
    \end{aligned}
    \]
    \begin{proof}
        For (a), by Exercise 12.10 in \cite{lattimore40bandit}, the result is true when $Z_t$ is replaced by $Y_t$, hence it suffices to notice that
        \[
        \E_{t-1}\left\{Y_{t-1}\frac{R_t(z)}{1-\lambda} \Sigma_t^{-1}\left[\frac{X_t-\lambda z}{1-\lambda}-\mu_t\right] \right\}=Y_{t-1}\E_{t-1}\left\{
        \frac{R_t(z)}{1-\lambda} \Sigma_t^{-1}\left[\frac{X_t-\lambda z}{1-\lambda}-\mu_t\right]
        \right\}=0,
        \]
        where the second equality follows from the definition of $R_t$ in Eq. \eqref{Ht} and it's similar for $\nabla^2\hat{s}_t(z).$
    \end{proof}
\end{Lem}
\begin{Lem}[Lemma 12.15 in \cite{lattimore40bandit}]\label{Rt 3}
    $R_t(\mu_t)\leq 3$ for all $t\leq n.$
\end{Lem}

\subsection{Concentration}

 \noindent
 There are many concentration properties for $s_t$ and $q_t$ in \cite{lattimore40bandit} and they need the condition that for all $t\leq\tau$,
\begin{equation}\label{condition concentration}
    \max \left(d, \operatorname{lip}(e), \sup _{x \in \mathcal{K}}|e(x)|,\|\Sigma_t\|,\left\|\Sigma_t^{-1}\right\|, 1 / \lambda\right) \leq \frac{1}{\delta}.
\end{equation}
 We first show that this is satisfied in our algorithms. To begin with, we show that in \textbf{RONM} (\textbf{ONM}), $\Sigma_t^{-1}$ grows quadratically 
 at most:
 \begin{Lem}\label{sk upper bound for st''}
    If $\eta<4$ and $\eta\gamma\leq \sigma^{-2}$, then for all $t\leq \tau$ and $z\in\R^d$, $
\left\|\nabla^2 s_t(z)\right\|\leq th,
$ 
$\left\|\Sigma_t^{-1}\right\|\leq \frac{3t^2h}{2}$, where $h=\max\left(\frac{1}{\sigma^2},\frac{4\lambda^2 \operatorname{lip}(e)^2d}{(1-\lambda)^2}\right).$ Especially, $\left\|\bar{\Sigma}_{\tau+1}^{-1}\right\|\leq (n+1)^2h$.
\begin{proof}
    By the definition of $\bar{\Sigma}_t^{-1}$ and that $\eta\gamma\leq \sigma^{-2}$, for $t\leq \tau$,
\begin{equation}\label{32ts}
    \left\|\Sigma_t^{-1}\right\|\leq\frac{3}{2}\left\|\bar{\Sigma}_t^{-1}\right\|\leq \frac{3}{2}\left(\frac{t}{\sigma^2}+\sum_{k=1}^{t-1} \eta\left\|\nabla^2 s_k(\mu_k)\right\|/2\right).
\end{equation}
   Then by Lemma \ref{upper bound for s''}, for all $t\leq\tau,$
we have
\[
\left\|\nabla^2 s_t(\mu_t)\right\| \leq \frac{\lambda \operatorname{lip}(e)}{1-\lambda} \sqrt{d\left\|\Sigma_t^{-1}\right\|}
\leq \frac{2\lambda \operatorname{lip}(e)d^{1/2}}{1-\lambda} \sqrt{\frac{t}{\sigma^2}+\sum_{k=1}^{t-1} \eta\left\|\nabla^2 s_k(\mu_k)\right\|/2}.
\]
Hence by Lemma \ref{induction}, 
 since $\eta<4$, we have 
$
\left\|\nabla^2 s_t(\mu_t)\right\|\leq th$, which is true for all $z\in\R^d$ because we can replace the first term with $\left\|\nabla^2 s_t(z)\right\|$ and also proves that  $\left\|\Sigma_t^{-1}\right\|\leq \frac{3t^2h}{2}$. Finally, for $\bar{\Sigma}_{\tau+1}^{-1}$, it suffices to note that the final inequality in Eq.~\eqref{32ts} is also true for $t=\tau+1$. 
\end{proof}
\end{Lem} 
\noindent
Also, this gives an exact upper bound for $\Sigma_t^{-1}$, in \cite{fokkema2024onlinenewtonmethodbandit}. Note that by the constants we choose, $h\leq d^2H^2$ and $\sigma^2\leq1$, which implies that for all $t\leq\tau$, we have
\begin{equation}\label{delta poly}
    \max(\|\Sigma_t\|,\|\Sigma_t^{-1}\|)\leq\operatorname{poly}(n,d,H).
\end{equation}
Then recalling that we have chosen $\delta=\operatorname{poly}(1/n,1/d,1/H)$ small enough, it's clear that Eq. \eqref{condition concentration} is met. Hence by \cite{lattimore40bandit}, in our algorithms, we have the following results.

\begin{Lem}[Proposition 12.22 in \cite{lattimore40bandit}]\label{q con}
    For all $x\in\R^d$, if $\max _{1 \leq t \leq \tau} \lambda\left\|x-\mu_t\right\|_{\Sigma_t^{-1}} \leq L^{-1/2}$ almost
surely, then with probability at least $1-\delta$,
    \[
    \left|\sum_{t=1}^\tau\left(q_t(x)-\hat{q}_t(x)\right)\right| \leq 1+\frac{1}{\lambda}\left[\sqrt{V_\tau L}+Z_{\max } L\right].
    \]
\end{Lem}
\begin{Lem}[Proposition 12.25 in \cite{lattimore40bandit}]\label{s'' con}
    Let $\mathscr{S} 
    $ be the (random) set of positive definite matrices such that $\Sigma_t^{-1} \preceq \Sigma^{-1}$ for all $t \leq \tau$ and $S_t=\sum_{u=1}^t \nabla^2\hat{s}_u\left(\mu_u\right)$ and $\bar{S}_t=$ $\sum_{u=1}^t \nabla^2 s_u\left(\mu_u\right)$. Then with probability at least $1-\delta$, for all $\Sigma^{-1} \in \mathscr{S}$,

$$
-\lambda L^2\left[1+\sqrt{d V_\tau}+d^2 Z_{\max }\right] \Sigma^{-1} \preceq S_\tau-\bar{S}_\tau \preceq\lambda L^2\left[1+\sqrt{d V_\tau}+d^2 Z_{\max }\right] \Sigma^{-1}.
$$

\end{Lem}
\begin{Lem}[modification of Lemma 10.15 in \cite{lattimore40bandit}]\label{trace}
    If for all $t\leq\tau$, $\eta\left\|\Sigma_t^{1 / 2} \nabla^2 s_t\left(\mu_t\right) \Sigma_t^{1 / 2}\right\|\leq 1$, then
    \[
    \frac{1}{\lambda} \sum_{t=1}^\tau \operatorname{tr}\left(\nabla^2 s_t\left(\mu_t\right) \Sigma_t\right)\leq \frac{8}{\lambda\eta}\operatorname{det}\left(\sigma^2\bar{\Sigma}_{\tau+1}^{-1}\right).
    \]
    \begin{proof}
        By lemma \ref{trace a}, 
        \[
        \begin{aligned}
            \frac{1}{\lambda} \sum_{t=1}^\tau \operatorname{tr}\left(\nabla^2 s_t\left(\mu_t\right) \Sigma_t\right)=\frac{4}{\lambda\eta} \sum_{t=1}^\tau \operatorname{tr}\left(\frac{\eta}{4}\nabla^2 s_t\left(\mu_t\right) \Sigma_t\right)\leq \frac{8}{\lambda\eta}\sum_{t=1}^\tau \log\operatorname{det}\left(\idm_d+\frac{\eta}{4}\nabla^2 s_t\left(\mu_t\right) \Sigma_t\right).
        \end{aligned}\]
        Note that $\Sigma_t\preceq2\bar{\Sigma}_t$ and $\frac{\eta}{2}\nabla^2 s_t\left(\mu_t\right)\preceq\frac{\eta}{2}\nabla^2 s_t\left(\mu_t\right)+\eta\gamma\idm_d=\bar{\Sigma}_{t+1}^{-1}-\bar{\Sigma}_t^{-1}$, then by lemma \ref{det plus},
        \[
        \begin{aligned}
              \operatorname{det}\left(\idm_d+\frac{\eta}{4}\nabla^2 s_t\left(\mu_t\right) \Sigma_t\right)\leq\operatorname{det}\left(\bar{\Sigma}_{t+1}^{-1}\bar{\Sigma}_t\right),
        \end{aligned}
        \]
        and the result follows from telescoping.
    \end{proof}
\end{Lem}
\begin{Lem}[Proposition 12.7 in \cite{lattimore40bandit}]\label{sq}
    For all $t\leq\tau$, we have 
    \[
    s_t(\mu_t)-s_t(x_\star)\leq q_t(\mu_t)-q_t(x_\star)+\frac{\delta}{\lambda^2}.
    \]
\end{Lem}

\begin{Lem}[Proposition 10.6 in \cite{lattimore40bandit}]\label{reg qreg}
     Suppose that $x \in \mathbb{R}^d$ satisfies $\lambda\|x-\mu_t\|_{\Sigma^{-1}_t} \leq \frac{1}{L}$ and $\lambda\leq d^{-1}L^{-2}$. Then for all $t\leq \tau$, 
$$
    \begin{aligned}
        e(\mu_t)-e(x_\star) \leq& \E_{t-1}[e(X_t)]-e(x_\star)\\
        \leq& q_t(\mu_t)-q_t(x_\star)+\frac{2}{\lambda} \operatorname{tr}\left(\nabla^2s_t(\mu_t) \Sigma_t\right)+\delta\left[\frac{2 d}{\lambda}+\frac{1}{\lambda^2}\right].
    \end{aligned}
$$
\end{Lem}

%% file: tex/append/concentration.tex
\begin{Lem}[Lemma 1 in \cite{10.1214/aos/1015957395}]\label{chi}
    If $X\sim\chi^2(k)$, then for all $x>0$,
    \[
    \mathbb{P}(X \geq k+2 \sqrt{k x}+2 x) \leq e^{-x}.
    \]
\end{Lem}
\begin{Lem}[Theorem B.17 in \cite{lattimore40bandit}]\label{concentration}
    Let $X_1, \ldots, X_n$ be a sequence of random variables adapted to filtration $\left(\mathscr{F}_t\right)$ and $\tau$ be a stopping time with respect to $\left(\mathscr{F}_t\right)_{t=1}^n$ with $\tau \leq n$ almost surely. Let $\mathbb{E}_t[\cdot]=\mathbb{E}\left[\cdot \mid \mathscr{F}_t\right]$. Then, with probability at least $1-\delta$,

$$
\left|\sum_{t=1}^\tau\left(X_t-\mathbb{E}_{t-1}\left[X_t\right]\right)\right| \leq 3 \sqrt{V_\tau \log \left(\frac{2 \max \left(B, \sqrt{V_\tau}\right)}{\delta}\right)}+2 B \log \left(\frac{2 \max \left(B, \sqrt{V_\tau}\right)}{\delta}\right),
$$
where $V_\tau=\sum_{t=1}^\tau \mathbb{E}_{t-1}\left[\left(X_t-\mathbb{E}_{t-1}\left[X_t\right]\right)^2\right]$ is the sum of the predictable variations and 
$B=\max \left(1,\right.$
$\left. \max _{1 \leq t \leq \tau}\left|X_t-\mathbb{E}_{t-1}\left[X_t\right]\right|\right)$.
\end{Lem}

\begin{Lem}[Proposition 2.5.2 in \cite{vershynin2018high}]\label{k1k2}
    If $W$ is a zero-mean random variable, then the following properties are equivalent; the parameters $K_i>0$ appearing in these
properties differ from each other by at most an absolute constant factor.
    \begin{enumerate}
        \item[(i)]There exists $K_1>0$ such that $\P(|W|\ge t)\leq 2\exp(-t^2/K_1^2)$, for all $t\ge 0;$
        \item[(ii)]There exists $K_2>0$ such that  $\mathbb{E} [\exp \left(\lambda^2 W^2\right)] \leq \exp \left(K_2^2 \lambda^2\right)$ for all $\lambda$ such that $|\lambda| \leq \frac{1}{K_2};$
        \item[(iii)]There exists $K_3>0$ such that $\|W\|_{L_p}=(\E[|W|^p])^{1/p}\leq K_3\sqrt{p}$ for all $p\ge 1;$
        \item[(iv)]There exists $K_4>0$ such that 
$\mathbb{E}[ \exp (\lambda W)] \leq \exp \left(K_4^2 \lambda^2\right)$ for all $\lambda \in \mathbb{R}$.
    \end{enumerate}
\end{Lem}

\begin{Lem}[Theorem 3.1.1 in \cite{vershynin2018high}]\label{norm subgaussian}
    Let $W=\left(W_1, \ldots, W_d\right) \in \mathbb{R}^d$ be a random vector with independent, subgaussian coordinates $W_i$ that satisfy $\mathbb{E} W_i^2=1$. Then
$$
\left\|\|W\|_2-\sqrt{d}\right\|_{\psi_2} \leq C K^2,
$$
where $K=\max _i\left\|W_i\right\|_{\psi_2}$ and $C$ is an absolute constant.
\end{Lem}

\begin{Lem}[Exercise 2.7.11 in \cite{vershynin2018high}]\label{norm}
    $\|\cdot\|_{\psi_2}$ is a norm on the space $\{W:\|W\|_{\psi_2}<+\infty\}.$
\end{Lem}

\begin{Lem}[Lemma B.6 in \cite{lattimore40bandit}]\label{norm expect}
    For any random variable $W$, $\|W-\mathbb{E}[W]\|_{\psi_2} \leq\left(1+\frac{1}{\log (2)}\right)\|W\|_{\psi_2}.$
\end{Lem}

\begin{Lem}\label{x-y subgaussian}
    If $X$ and $Y$ are independent and identically distributed $1$-subgaussian variables, then there exists $C<+\infty$ such that $|X-Y|-\E[|X-Y|]$ is also $C$-subgaussian.
    \begin{proof}
        By Lemma \ref{k1k2}, there exists $C_1>0$ such that $\|X\|_{\psi_2}=C_1<+\infty$. Then by Lemma \ref{norm} and Lemma \ref{norm expect}, it's clear that 
        \[
        \||X-Y|-\E[|X-Y|]\|_{\psi_2}\leq \left(1+\frac{1}{\log (2)}\right)\|X-Y\|_{\psi_2}\leq 2\left(1+\frac{1}{\log (2)}\right)C_1<+\infty.
        \]
        Hence, again, by Lemma \ref{k1k2}, there exists $C_2>0$ such that $|X-Y|-\E[|X-Y|]$ is $C_2$-subgaussian.
    \end{proof}
\end{Lem}

%% file: tex/append/auxi.tex
\begin{Lem}[Lemma 3.3(g) in \cite{lattimore40bandit}]\label{lip pi}
    Let $\mathcal{K}$ be a convex body and $\pi$ the associated Minkowski functional. Then $\lip(\pi)\leq 1/r$ whenever $\mathcal{K}\supset\mathbb{B}_r^d.$
\end{Lem}

\begin{Lem}[Lemma A.5 in \cite{lattimore40bandit}]\label{trace a}
    Suppose that $A$ is positive semidefinite and $A\preceq\idm$. Then $\tr{A}\leq 2\log\operatorname{det}(\idm+A).$
\end{Lem}

\begin{Lem}[Corollary \uppercase\expandafter{\romannumeral3}.1.2 in \cite{bhatia97}]\label{lambda k}
    Suppose that $A,B$ are both positive semidefinite in $\R^{d\times d}$ and $A\succeq B$. Then for all $1\leq k\leq d$, the $k$-th smallest eigenvalue of $A$ is also larger than $B$'s.
\end{Lem}

\begin{Lem}\label{min deter}
    If $A,B$ are both positive semidefinite and $A\succeq B$, then  
    $\operatorname{det}(A)\ge\operatorname{det}(B)$.
    \begin{proof}
        Note that the determinant is just the product of all eigenvalues, which are all non-negative for positive semidefinite matrices, then it follows from Lemma \ref{lambda k}.
    \end{proof}
\end{Lem}
\begin{Lem}\label{det plus}
    Suppose that $A,B$ and $C$ are positive semidefinite. If $B\preceq C$, then
    $\operatorname{det}(\idm+BA)=\operatorname{det}(\idm+AB)\leq \operatorname{det}(\idm+AC)=\operatorname{det}(\idm+CA)$.
    \begin{proof}
        The equality follows from the folklore that $AB$ and $BA$ have the same non-zero eigenvalues. For the inequality, first suppose that $A$ is positive definite, then 
        \[
        \operatorname{det}(\idm+AB)=\operatorname{det}\left(A^{-1/2}(\idm+AB)A^{1/2}\right)=\operatorname{det}(\idm+A^{1/2}BA^{1/2}).
        \]
        Note that $A^{1/2}BA^{1/2}\preceq A^{1/2}CA^{1/2}$, then by Lemma \ref{min deter}, clearly,
        \[
        \operatorname{det}(\idm+A^{1/2}BA^{1/2})\leq \operatorname{det}(\idm+A^{1/2}CA^{1/2})=\operatorname{det}(\idm+AC).
        \]
        When $A$ is positive semidefinite, note that for all $t>0$,
        \[
        \operatorname{det}(\idm+(A+t\idm)B)\leq \operatorname{det}(\idm+(A+t\idm)C),
        \] 
        then it suffices to let $t\to0.$
    \end{proof}
\end{Lem}
\begin{Lem}[Proposition A.4(b) in \cite{lattimore40bandit}]\label{4 moment}
    If $W\sim\N(0,\Sigma)$, then $\E[\|W\|_2^4]=\tr{\Sigma}^2+2\tr{\Sigma^2}$
\end{Lem}

\begin{Lem}\label{induction}
    Given a sequence of positive numbers $x_n$ that satisfies that
    \[
    \,x_n\leq a\sqrt{bn+c\sum_{k=1}^{n-1} x_k}, \forall n\ge 1,
    \]
    where $a, b, c>0$ and $c\leq 2$, and letting $h=\max\{a^2,b\}$, we have $x_n\leq hn,\forall n\ge 1.$
    \begin{proof}
        We prove this result by induction. It's true for $n=1$ since $x_1\leq a\sqrt{b}\leq h$. Assume that it holds for all $n\leq m$, $m>1$. Then
        \[
        x_{m+1}\leq a\sqrt{b(m+1)+c\sum_{k=1}^{m} hk}=a\sqrt{b(m+1)+\frac{chm(m+1)}{2}}\leq h\sqrt{m+1+m(m+1)}= h(m+1),
        \]
        which completes the proof.
    \end{proof}
\end{Lem}
\begin{Lem}\label{growth rate lamma}
    Assume that a sequence of positive numbers $x_n$ satisfies that
    \[
    \,x_{n+1}-x_n\ge ax_n^{\frac{2-b}{2}}, \forall n\ge 1,
    \]
    where $a>0$ and $1\leq b\leq2$. Then if $x_1^{b/2}\ge a$, we have $x_n\ge (an)^{\frac{2}{b}}/8,\forall n\ge 1.$
    \begin{proof}
        Because $x+ax^{\frac{2-b}{2}}$ is increasing in $\R^+$, W.L.O.G., we can assume that the equality always holds for all $n\ge 1$. Let $h(x)=x^{b/2}$, then by Lagrange's mean value theorem and that $h'(x)=\frac{b}{2}x^{\frac{b-2}{2}}$ is decreasing and $x_n$ is increasing, we have
        \[
        h(x_{n+1})-h(x_n)\ge h'(x_{n+1})(x_{n+1}-x_n)=\frac{ab}{2}x_{n+1}^{\frac{b-2}{2}}x_n^{\frac{2-b}{2}}=\frac{ab}{2}(1+ax_n^{-\frac{b}{2}})^{\frac{b-2}{2}}.
        \]
        Since $x_n\ge x_1$ and $1\leq b\leq 2$, the right hand is larger than
        $
        \frac{ab}{2}(1+ax_1^{-\frac{b}{2}})^{\frac{b-2}{2}}.
        $ Therefore, $h(x_n)-h(x_1)\ge (n-1)\frac{ab}{2}(1+ax_1^{-\frac{b}{2}})^{\frac{b-2}{2}},$ which implies that
        \[
        x_n^{b/2}\ge \frac{abn}{2}(1+ax_1^{-\frac{b}{2}})^{\frac{b-2}{2}},
        \]
        where we used that $\frac{ab}{2}(1+ax_1^{-\frac{b}{2}})^{\frac{b-2}{2}}\leq \frac{ab}{2}\leq a\leq x_1^{b/2}$. 
        Then note that $\frac{ab}{2}(1+ax_1^{-\frac{b}{2}})^{\frac{b-2}{2}}\ge\frac{ab}{2}\cdot 2^{\frac{b-2}{2}}$ and $b^{\frac2b}\ge1$. We thus  have
        \[
        x_n\ge \left(\frac{ab}{2}(1+ax_1^{-\frac{b}{2}})^{\frac{b-2}{2}}\right)^{\frac{2}{b}}\cdot  n^{\frac{2}{b}}\ge2^{1-\frac4b}(abn)^{\frac{2}{b}}\ge (abn)^{\frac{2}{b}}/8\ge (an)^{\frac2b}/8.
        \]
    \end{proof}
\end{Lem}

\begin{Lem}[Cauchy–Binet formula]\label{cauchy binet} It holds that
    \[
    \left(\sum_{i=1}^na_i^2\right)\left(\sum_{i=1}^nb_i^2\right)-\left(\sum_{i=1}^na_ib_i\right)^2=\sum_{i,j=1}^n\left(a_ib_j-a_jb_i\right)^2.
    \]
\end{Lem}

\begin{Lem}[Equivalence of norms]\label{eui norm}
    For all $1\leq p\leq 2$ and $x\in\R^d$, we have
    \[
    \|x\|_2\leq \|x\|_p\leq d^{\frac{1}{p}-\frac{1}{2}}\|x\|_2. 
    \]
\end{Lem}

\begin{Lem}\label{vol of ball}
    $\operatorname{Vol}_d(\mathbb{B}_R^d)=\frac{\pi^{d / 2}R^d}{\Gamma\left(\frac{d}{2}+1\right)}$, where $\Gamma$ is the gamma function.
\end{Lem}
\begin{Lem}[Gautschi's inequality]\label{gautschi}
    $\text {Let } x \text { be a positive real number, and let } s \in(0,1) \text {. Then, }x^{1-s}<\frac{\Gamma(x+1)}{\Gamma(x+s)}<(x+1)^{1-s} .$
\end{Lem}

\begin{Lem}[Proposition 3.19. in \cite{lattimore40bandit}]\label{extension}
   If $e(x)$ is the convex extension defined in Section \ref{extension section}, then it satisfies the following:
    \[
    \begin{aligned}
        &(a)\, e(x)=f(x) \text{ for all } x \in K;  
        &&(b)\,e \text{ is convex};\\
&(c)\, \lip(e) \leq \frac{2 G R}{r}+G+\frac{1}{r}; 
&&(d)\, \text{For all } x \notin \mathcal{K}, e(x / \pi(x)) \leq e(x).
    \end{aligned}\]
\end{Lem}

\begin{Lem}[Theorem 10.2 in \cite{lattimore40bandit}]\label{regret for onm}
    Suppose that $\Sigma_t^{-1}$ is positive definite for all $1 \leq t \leq n+1$. Then in the online Newton method defined in Section \ref{pre onm}, for all $x \in \mathcal{K}$,

$$
\frac{1}{2}\left\|\mu_{n+1}-x\right\|_{\Sigma_{n+1}^{-1}}^2 \leq \frac{1}{2}\left\|\mu_1-x\right\|_{\Sigma_1^{-1}}^2+\frac{\eta^2}{2} \sum_{t=1}^n\left\|g_t\right\|_{\Sigma_{t+1}}^2-\eta \widehat{\operatorname{qReg}}_n(x) .
$$
\end{Lem}

%% file: tex/append/constraint.tex
In Algorithm \ref{algorithm 1}, for $\rho$-QG function, as Eq.~\eqref{alg1 1} we need $\gamma=\rho$ and
\[
\begin{aligned}
    &\bullet\quad\eta\leq4, \sigma^{-2}\ge\eta\gamma,  \text{ Lemma \ref{sk upper bound for st''}},&&\bullet\quad
    h=\max\left(\frac{1}{\sigma^2},\frac{4\lambda^2 \operatorname{lip}(e)^2d}{(1-\lambda)^2}\right)\leq d^2H^2, \sigma^2\leq1,\text{ Eq.~\eqref{delta poly}},\\
    &\bullet\quad\lambda\leq d^{-1/2}L^{-3/2},\text{Eq.~\eqref{lambda1}},&&\bullet\quad\lambda\leq d^{-1}L^{-2},\text{ Lemma \ref{reg qreg}},\\
    &\bullet\quad
    \lip(e)\leq 
    \frac{2 G R}{r}+G+\frac{1}{r}\leq H,\text{ Eq.~\eqref{Zt}},
    &&\bullet\quad\lambda\leq 1/2,\text{ Eq.~\eqref{gt}},\\
    &\bullet\quad H\eta\lambda\sigma\sqrt{d}\leq 1,\text{ Eq.~\eqref{condition 1}},
    &&\bullet\quad \eta \frac{Hd^2L^2}{\sqrt{\eta\gamma}}\leq 2/3,\text{ Eq.~\eqref{condition 2}},
        \\&\bullet\quad
    \frac{R^2}{2\sigma^2}, \frac{\eta H^2dL}{2\lambda^2\gamma},2\eta
    ,
     \frac{H\sqrt{\eta L}}{\lambda^2\sqrt{\gamma }},\frac{dL}{\lambda}\leq\frac{1}{10\lambda^2L^2},\text{ Eq.~\eqref{solution 2}},&&\bullet\quad \frac{H}{\lambda\sqrt{L\eta\gamma}}\ge 3, \lambda\leq\frac{1}{2\sqrt{d}L}
     \text{ Lemma \ref{sk Z V}}.
\end{aligned}
\]
For $(\beta,\bc)$-convex function, where $1<\bc\leq 2$, in addition, 
we need
\[
\begin{aligned}
    &\bullet\quad\gamma=
    2^{\bc-1}\beta
    ,\text{ Lemma \ref{q qg}},\quad&&\bullet\quad\sigma^{-2}\ge\left(\frac{\bc-1}{30}\right)^{2/\bc}\beta^{\frac2\bc}d^{-\frac1\bc}\left(\frac{r}{\sqrt{2}}\right)^{\frac{2(d-1)}{\bc}}\eta^{\frac2\bc}\lambda^{\frac6\bc-2}L^{\frac4\bc-2},\text{ Lemma \ref{growth rate >1}},\\
    &\bullet\quad\frac{r}{\sqrt{2}\sigma}\ge 5d, \lambda\leq\frac{1}{10dL},\text{ Lemma \ref{growth rate >1}}.
\end{aligned}\]
\noindent
Note that the reader should also use that $\gamma\leq \frac{8}{R^2}
=8
$ by Lemma \ref{upper rho} to check them. 

In Algorithm \ref{algorithm 2}, for $(\beta,1)$-convex function, as Eq.~\eqref{alg 2 con}
, we need
\[
\begin{aligned}
    &\bullet\quad\eta\leq4,\text{ Lemma \ref{sk upper bound for st''}},&&\bullet\quad\lambda\leq1-\frac{1}{\sqrt2},\text{ Lemma \ref{growth rate=1}},\\
    &\bullet\quad\sigma^{-2}\ge\Theta=\beta^{2-\kappa}\eta^{2-\kappa}6^{-\frac{d(2-\kappa)}{2}}e^{-\frac{2-\kappa}{\lambda^2L^2}}/32,\sigma^{-2}\ge 1,\text{Lemma \ref{growth rate=1}},
    &&\bullet\quad
    \frac{HJ\sqrt{L}}{\sqrt{\Theta}}\ge 3
    ,
    \operatorname{E}_3,\operatorname{E}_4
    ,\\
    &\bullet\quad \frac{R^2}{2\sigma^2},\frac{d\eta^2H^2J^2L^3}{2\Theta},2\eta
    ,
     \frac{\eta HJL^{\frac32}}{\lambda\sqrt{\Theta}}\leq\frac{1}{8\lambda^2L^2},\text{ Eq.~\eqref{alg 2}},
    &&\bullet\quad \frac{\eta\lambda d^2HJL^{3}}{\sqrt{\Theta}}\leq 2/3,\text{ Eq.~\eqref{new def bc}},\\
    &\bullet\quad
    h=\max\left(\frac{1}{\sigma^2},\frac{4\lambda^2 \operatorname{lip}(e)^2d}{(1-\lambda)^2}\right)\leq d^2H^2, \sigma^2\leq1,\text{ Eq.~\eqref{delta poly}}.&&
\end{aligned}
\]
Note that the reader should also use that $\beta\leq \frac{2}{R}=2$ by Lemma \ref{1 qg}
 to check them.

%% file: tex/append/exponential_dependence.tex
\subsection{\texorpdfstring{Exponential dependence on $d$ can be removed when $\bc=2$
}{Exponential dependence on d can be removed when q=2
}}\label{exponential dependence}
    Though we have established the growth rate for $\Sigma_t^{-1}$ when $f(x)$ is $(\beta,\bc)$-convex, where $1< \bc\leq 2$, its exponential dependence on $d$ is very undesirable, which is from Lemma \ref{proportion}. However, it's easy to see that Lemma \ref{proportion} is impossible to get improved in the sense of removing the exponential dependence on $d$ because one can just take a hypercube as an example. Any small ball centered at any vertex of the hypercube will at most has its $\frac{1}{2^d}$ inside the cube. 
    
    Therefore, the only way to remove the exponential dependence on $d$ is to improve Eq.~\eqref{xtk}. If that
    \[
    \nabla^2 e(x)\succeq C\|x-x_{\star}\|_2^{\bc-2}
    \]
    is true not only on $\mathcal{K}$ but a larger area, one can skip Lemma \ref{proportion} and achieve the goal. 
    
    Actually, this is possible when $\bc=2$, say $f(x)$ is $\alpha$-strongly convex, because we find a new extension of $f(x)$, which is also strongly convex in the neighbourhood of $\mathcal{K}$. The idea is very straightforward. Since $g(x)=f(x)-\frac{\alpha}{2}\|x\|_2^2$ is convex, we can just extend $g(x)$ with the method introduced in Section \ref{extension section} and then add the quadratic term back. The explicit expression of it can be found in Appendix \ref{exten app}. 

    With such a strongly convex extension, we have:
    \begin{Lem}\label{strongly convex}
    When $f(x)$ is $\alpha$-strongly convex, given any $\epsilon, \sigma>0$, if $\frac{\epsilon}{\sigma}\ge 10d$ and $\Sigma\preceq\sigma^2\idm_d$, then for all $\mu\in\mathbb{B}_{R}^d$, the surrogate function $s(x)$ with $\N(\mu,\Sigma)$,  $e(x)$, which is the strongly convex extension with respect to $f(x)$ and $\epsilon$ defined in Appendix \ref{exten app} and $\lambda\in(0,1)$  satisfies that 
    \[
    \nabla^2s(z)\succeq \frac{\lambda\alpha}{2}\idm_d,
    \]
    for all $z\in\mathbb{B}_R^d$.
\end{Lem}
Its proof can be found in Appendix \ref{proof for strongly convex} and is similar to Eq.~\eqref{EtHt}. This lemma shows that, informally, in every time step $t$, $\Sigma^{-1}_t$ will be added by $\lambda\eta\alpha\idm_d$ without any exponential dependence on $d$. Unfortunately, one should note that when $f(x)$ is $\alpha$-strongly convex, it's also $\alpha$-QG and hence one can also apply \textbf{RONM}. Recall that in \textbf{RONM}, in every time step $t$, $\Sigma_t^{-1}$ is added by $\eta\alpha\idm_d$ (we omit other constants temporarily), which is much larger than $\lambda\eta\alpha\idm_d$ because $\lambda$ will be very small. Hence though we remove the exponential dependence on dimension by the strongly convex extension, the contribution to $\E_{t-1}[H_t]$ is still less than the regularized term and then this won't improve the orders of the regret and the convergence rate.

It remains unknown whether it is possible to get a similar extension when $\bc<2$, which is difficult 
since the idea of the strongly convex extension fails. This is because in the strongly convex extension, we rely on the fact that if $f(x)$ is $(\beta,2)$-convex, i.e., $f(x)-\beta\|x-x_\star\|_2^2$ is convex, then $f(x)-\beta\|x\|_2^2$ is also convex. However, this is not true for general $\bc$. 

\subsection{Strongly convex extension}\label{exten app}
In this section, we introduce the \emph{strongly convex extension}, which is strongly convex in the neighbourhood of $\mathcal{K}$, to remove the exponential dependence on $d$. Later, we assume that $f(x)$ is $\alpha$-strongly convex on $\mathcal{K}$.

Let $g(x)=f(x)-\frac{\alpha}{2}(\|x\|_2^2-R^2),\,\forall x\in\mathcal{K}$. Then one can see that
\[
\lip_{\mathcal{K}}(g)\leq \lip_{\mathcal{K}}(f)+\alpha R\leq G+\alpha R, \text{and }0\leq g(x)\leq 1+\frac{\alpha R^2}{2}.
\]
Then for all $\epsilon>0$, we define that
\[
\widetilde{e}_\epsilon(x)=\pi^+(x) g\left(\frac{x}{\pi^+(x)}\right)+(G+\alpha R)R(\pi^+(x)-1)+\frac{\alpha}{2}(\|x\|_2^2-R^2), 
\]
$\text{for all $x$ such that }\|x\|_2\leq R+\epsilon.$ It can be shown that:
\begin{Lem}\label{wt e}
    $\widetilde{e}_\epsilon(x)$ satisfies the following:
    \[
    \begin{aligned}
        &(a)\widetilde{e}_\epsilon(x)=f(x) \text{ for all $x\in\mathcal{K}$}; \\&(b) \widetilde{e}_\epsilon \text{ is $\alpha$-strongly convex on $\mathbb{B}_{R+\epsilon}^d$};\\
        &(c)\lip_{\mathbb{B}_{R+\epsilon}^d}(\widetilde{e}_\epsilon)\leq \frac{2R(G+\alpha R)}{r}+(G+\alpha R)+\frac{1}{r}\left(1+\frac{\alpha R^2}{2}\right)+\alpha(R+\epsilon):=G_{\epsilon}/4;\\&(d) \text{For all $x\in\mathcal{K}^{c}\cap\mathbb{B}_{R+\epsilon}^d$, } \widetilde{e}_\epsilon(x / \pi(x)) \leq \widetilde{e}_\epsilon(x);
        \\&(e)0\leq \widetilde{e}_\epsilon(x)\leq (R+\epsilon)G_{\epsilon}/4:=M_\epsilon.
    \end{aligned}
    \]
\end{Lem}
The proof for Lemma \ref{wt e} is deferred in Appendix \ref{proof for wte}.
Finally we extend $\widetilde{e}_\epsilon$ to $\R^d$ by the original way in Section \ref{extension section}:
\[
\begin{aligned}
    e_\epsilon(x)=\pi_\epsilon^+(x) \widetilde{e}_\epsilon\left(\frac{x}{\pi_\epsilon^+(x)}\right)+
M_\epsilon\left(\pi_\epsilon^+(x)-1\right),
\end{aligned}
\]
where $\pi_\epsilon(x)=\frac{\|x\|_2
}{R+\epsilon}$ is the Minkowski functional of $\mathbb{B}_{R+\epsilon}^d.$ Then similarly, we have
\begin{Lem}\label{e ep}
    $e_\epsilon(x)$ satisfies the following:
    \[
    \begin{aligned}
        &(a)e_\epsilon(x)=f(x) \text{ for all $x\in\mathcal{K}$}; &&(b) e_\epsilon \text{ is $\alpha$-strongly convex on $\mathbb{B}_{R+\epsilon}^d$};\\
        &(c)\lip(e_\epsilon)\leq G_{\epsilon};&&(d) \text{For all $x\in\mathcal{K}$, } e_\epsilon(x / \pi(x)) \leq e_\epsilon(x).
    \end{aligned}
    \]
\end{Lem}
Its proof can be found in Appendix \ref{proof for e ep}. By Lemma \ref{pipi} and direct computation in Appendix \ref{Computation of eepsilon}, we have
\[
\begin{aligned}
    e_\epsilon(x)=\pi^+(x) f\left(\frac{x}{\pi^+(x)}\right)+(G+\frac32\alpha R)R(\pi^+(x)-\pi_\epsilon^+(x))+\frac{\alpha}{2}\|x\|_2^2(1/\pi_\epsilon^+(x)-1/\pi^+(x))+M_\epsilon\left(\pi_\epsilon^+(x)-1\right).
\end{aligned}
\]
When $X$ is chosen by the learner, the learner actually picks $\frac{X}{ \pi^+(X)}$ and the bandit gives $f(\frac{X}{ \pi^+(X)})+\varepsilon$. Simply substituting $f(\frac{X}{ \pi^+(X)})$ with $f(\frac{X}{ \pi^+(X)})+\varepsilon$, we can feed the player with the loss
\[
\begin{aligned}
    Y=
    e_{\epsilon}(X)+\pi^+(X)\varepsilon.
\end{aligned}
\]
Then from the learner's perspective, the loss function is $e_\epsilon(x)$ and noise is $\xi:=\pi^+(X)\varepsilon$. Following the same arguments, 
we have that Eq.~\eqref{extension noise} is also true for the strongly convex extension.

\subsubsection{Some useful facts}
\begin{Lem}\label{lip e ep}
    If $\epsilon<R$, then $\lip(e_\epsilon)\leq \frac{12RG}{r}+\frac{48R^2}{r^3}.$
    \begin{proof}
        By Lemma \ref{e ep},
    \[
    \begin{aligned}
        \lip(e_\epsilon)\leq G_\epsilon=&\frac{8R(G+\alpha R)}{r}+4(G+\alpha R)+\frac{4}{r}\left(1+\frac{\alpha R^2}{2}\right)+4\alpha(R+\epsilon)\\
        \leq &\frac{8R(G+2R/r^2)}{r}+4(G+2R/r^2)+\frac{4}{r}(1+\frac{R^2}{r^2})+\frac{16R}{r^2}\\
        \leq &\frac{12RG}{r}+\frac{48R^2}{r^3},
    \end{aligned}
    \]
        where the second inequality is by Lemma \ref{upper rho} and the final inequality used that $R\ge r.$
    \end{proof}
\end{Lem}
\begin{Lem}\label{pipi}
    For all $x\in\R^d$,  we have
    \[
    \pi_\epsilon^+(x)\pi^+\left(\frac{x}{\pi_\epsilon^+(x)}\right)=\pi^+(x).
    \]
    \begin{proof}
        If $x\notin\mathcal{K}$, it follows from positive homogeneity of $\pi(x).$ Otherwise, $\pi_\epsilon^+(x)=1$ because $\mathcal{K}\subset\mathbb{B}_{R+\epsilon}^d$. Accordingly, the equality clearly holds.
    \end{proof}
\end{Lem}
\subsection{Proofs for Lemmas}\label{proof for lemma in sc}
\subsubsection{\texorpdfstring{Computation of $e_\epsilon$}{Computation of eepsilon}}\label{Computation of eepsilon}
\[
\begin{aligned}
    &e_\epsilon(x)=\pi_\epsilon^+(x) \widetilde{e}_\epsilon\left(\frac{x}{\pi_\epsilon^+(x)}\right)+
M_\epsilon\left(\pi_\epsilon^+(x)-1\right)\\
&=\pi_\epsilon^+(x) \left[
\pi^+(x/\pi_\epsilon^+(x)) g\left(\frac{x/\pi_\epsilon^+(x)}{\pi^+(x/\pi_\epsilon^+(x))}\right)+(G+\alpha R)R(\pi^+(x/\pi_\epsilon^+(x))-1)+\frac{\alpha}{2}(\|x/\pi_\epsilon^+(x)\|_2^2-R^2)
\right]
\\  & \qquad +M_\epsilon\left(\pi_\epsilon^+(x)-1\right)\\
&=\pi^+(x) g\left(\frac{x}{\pi^+(x)}\right)+(G+\alpha R)R(\pi^+(x)-\pi_\epsilon^+(x))+\frac{\alpha}{2}(\|x\|_2^2/\pi_\epsilon^+(x)-R^2\pi_\epsilon^+(x))+M_\epsilon\left(\pi_\epsilon^+(x)-1\right)\\
&=\pi^+(x) f\left(\frac{x}{\pi^+(x)}\right)+\left(G+\frac32\alpha R\right)R(\pi^+(x)-\pi_\epsilon^+(x))+\frac{\alpha}{2}\|x\|_2^2(1/\pi_\epsilon^+(x)-1/\pi^+(x))+M_\epsilon\left(\pi_\epsilon^+(x)-1\right),
\end{aligned}
\]
where the second equality follows from Lemma \ref{pipi} and the final equality follows from that $g(x)=f(x)-\frac{\alpha}{2}(\|x\|_2^2-R^2)$.
\subsubsection{Proof for Lemma \ref{wt e}}\label{proof for wte}
\begin{proof}
    Since
    $
\lip_{\mathcal{K}}(g)\leq \lip_{\mathcal{K}}(f)+\alpha R\leq G+\alpha R \text{ and }0\leq g(x)\leq 1+\frac{\alpha R^2}{2},
$
let $\bar{g}(x)=\frac{g(x)}{1+\frac{\alpha R^2}{2}}$, then 
\[0\leq \bar{g}(x)\leq 1, \lip_{\mathcal{K}}(\bar{g})=\frac{\lip_{\mathcal{K}}(g)}{1+\frac{\alpha R^2}{2}}\leq \frac{G+\alpha R}{1+\frac{\alpha R^2}{2}}.\]
Thus we can extend $\bar{g}(x)$ by Lemma \ref{extension} as
\[
\widetilde{\bar{g}}(x)=\pi^+(x) \bar{g}\left(\frac{x}{\pi^+(x)}\right)+\frac{G+\alpha R}{1+\frac{\alpha R^2}{2}}\cdot R(\pi^+(x)-1).
\]
Let $$\widetilde{g}(x)=\left(1+\frac{\alpha R^2}{2}\right)\widetilde{\bar{g}}(x)=\pi^+(x) g\left(\frac{x}{\pi^+(x)}\right)+(G+\alpha R) R(\pi^+(x)-1),$$
and note that $\widetilde{e}_\epsilon(x)=\widetilde{g}(x)+\frac{\alpha}{2}(\|x\|_2^2-R^2)$. Then by Lemma \ref{extension} (a), for all $x\in\mathcal{K}$,
\[
\widetilde{e}_\epsilon(x)=\left(1+\frac{\alpha R^2}{2}\right)\widetilde{\bar{g}}(x)+\frac{\alpha}{2}(\|x\|_2^2-R^2)=\left(1+\frac{\alpha R^2}{2}\right)\bar{g}(x)+\frac{\alpha}{2}(\|x\|_2^2-R^2)=g(x)+\frac{\alpha}{2}(\|x\|_2^2-R^2)=f(x),
\]
which yields (a). By Lemma \ref{extension} (b), $\widetilde{g}(x)$ is convex and then $\widetilde{e}_\epsilon(x)$ is $\alpha$-strongly convex, which gives (b). For part (c), for all $x\in\mathbb{B}_{R+\epsilon}^d$,
\[
\begin{aligned}
    \lip_{\mathbb{B}_{R+\epsilon}^d}(\widetilde{e}_\epsilon)\leq &\lip_{\mathbb{B}_{R+\epsilon}^d}(\widetilde{g})+\alpha(R+\epsilon)\\
    =&\left(1+\frac{\alpha R^2}{2}\right)\lip_{\mathbb{B}_{R+\epsilon}^d}(\widetilde{\bar{g}})+\alpha(R+\epsilon)\\
    \leq &\left(1+\frac{\alpha R^2}{2}\right)\cdot\left(
    \frac{2 \frac{G+\alpha R}{1+\frac{\alpha R^2}{2}}R}{r}+\frac{G+\alpha R}{1+\frac{\alpha R^2}{2}}+\frac{1}{r}
    \right)+\alpha(R+\epsilon)\\
    \leq &\left(
    \frac{2 (G+\alpha R)R}{r}+G+\alpha R+\frac{1+\frac{\alpha R^2}{2}}{r}
    \right)+\alpha(R+\epsilon)=G_\epsilon/4,
\end{aligned}
\]
where the second inequality used Lemma \ref{extension} (c). For part (d), for all $x\in\mathcal{K}^{c}\cap\mathbb{B}_{R+\epsilon}^d$ we have
\[
\begin{aligned}
   \widetilde{e}_\epsilon\left(\frac{x}{\pi(x)}\right)=&\left(1+\frac{\alpha R^2}{2}\right)\widetilde{\bar{g}}\left(\frac{x}{\pi(x)}\right)+\frac{\alpha}{2}\left(\left\|\frac{x}{\pi(x)}\right\|_2^2-R^2\right)\\
   \leq &\left(1+\frac{\alpha R^2}{2}\right)\bar{g}\left(x\right)+\frac{\alpha}{2}\left(\left\|x\right\|_2^2-R^2\right)=\widetilde{e}_\epsilon(x),
\end{aligned}
\]
where the inequality used Lemma \ref{extension} (d) and $\pi(x)\ge 1$. Finally, for part (d), by part (c) it suffices to show that $\widetilde{e}_\epsilon(x)\ge 0$ for all $x\in\mathbb{B}_{R+\epsilon}^d$, which is true if $x\in\mathcal{K}$ by part (a). For all $x\in\mathcal{K}^{c}\cap\mathbb{B}_{R+\epsilon}^d$, by part (d), $\widetilde{e}_\epsilon\left(\frac{x}{\pi(x)}\right)\ge \widetilde{e}_\epsilon(x)\ge  0.$
\end{proof}

\subsubsection{Proof for Lemma \ref{e ep}}\label{proof for e ep}
\begin{proof}
 We apply the same reasoning in Lemma \ref{wt e}. Let $\overline{\widetilde{e}}_\epsilon(x)=\widetilde{e}_\epsilon(x)/M_\epsilon$, with its extension
 \[
 \widetilde{\overline{\widetilde{e}}}_\epsilon(x)=\pi^+_{\epsilon}(x)\overline{\widetilde{e}}_\epsilon\left(\frac{x}{\pi^+_\epsilon(x)}\right)+\frac{G_\epsilon/4}{M_\epsilon}(R+\epsilon)(\pi^+_{\epsilon}(x)-1).
 \]
Then $e_\epsilon(x)=M_\epsilon\widetilde{\overline{\widetilde{e}}}_\epsilon(x).$ Part (a) follows from Lemma \ref{extension} (a) and Lemma \ref{wt e} (a). Part (b) follows from Lemma \ref{extension} (a) and Lemma \ref{wt e} (b). For part (c), by Lemma \ref{extension} (c), we have
\[
\lip(e_\epsilon)=M_\epsilon\lip(\widetilde{\overline{\widetilde{e}}}_\epsilon)\leq M_\epsilon\left(
\frac{2 G_\epsilon (R+\epsilon)/4}{ M_\epsilon(R+\epsilon)}+G_\epsilon/4M_\epsilon+\frac{1}{R+\epsilon}
\right)=G_\epsilon.
\]
\end{proof}

\subsubsection{Proof for Lemma \ref{strongly convex}}\label{proof for strongly convex}
\begin{proof}
         Let $\bar{\mu}=(1-\lambda)\mu+\lambda z$, then clearly $\mathbb{B}_{\epsilon}^d(\bar{\mu})\subset\mathbb{B}_{R+\epsilon}.$
        Following the same arguments in Eq.~\eqref{EtHt}, we have
        \[
        \begin{aligned}
            \nabla^2 s(z)&\succeq\lambda\E\left[\nabla^2e\left(\widetilde{X}\right)\right],\quad &&\text{where }\widetilde{X}\sim\N(\bar{\mu},(1-\lambda)^2\Sigma),
            \\&\succeq \lambda\alpha\idm_d\cdot\P(\widetilde{X}\in\mathbb{B}_{\epsilon/2}^d(\bar{\mu})),&&\text{since $e(x)$ is $\alpha$-strongly convex on $\mathbb{B}_{\epsilon/2}^d(\bar{\mu})$}.
        \end{aligned}
        \]
        
        Now it suffices to lower bound $\P(\widetilde{X}\in\mathbb{B}_{\epsilon/2}^d(\bar{\mu}))$, which is just $\P(\bar{X}\in\Sigma^{-1/2}\mathbb{B}_{\frac{\epsilon}{2(1-\lambda)}}^d)$, where $\bar{X}\sim\N(0,\idm_d)$. Since $\Sigma^{-1/2}\succeq\idm_d/\sigma$, then $\Sigma^{-1/2}\mathbb{B}_{\frac{\epsilon}{2(1-\lambda)}}^d\supset\mathbb{B}_{\epsilon/\sigma}^d$
        , hence we have
        \[
\P(\bar{X}\in\Sigma^{-1/2}\mathbb{B}_{\frac{\epsilon}{2(1-\lambda)}}^d)\ge \P(\bar{X}\in\mathbb{B}_{\epsilon/2\sigma}^d)\ge\P(\bar{X}\in\mathbb{B}_{5d}^d),
        \]
which is larger than $1-e^{-d}\ge 1/2$ by Lemma \ref{chi}. 
    \end{proof}

%% file: main.bbl
\begin{thebibliography}{25}
\providecommand{\natexlab}[1]{#1}
\providecommand{\url}[1]{\texttt{#1}}
\expandafter\ifx\csname urlstyle\endcsname\relax
  \providecommand{\doi}[1]{doi: #1}\else
  \providecommand{\doi}{doi: \begingroup \urlstyle{rm}\Url}\fi

\bibitem[Agarwal et~al.(2013)Agarwal, Foster, Hsu, Kakade, and Rakhlin]{AFHK13}
Alekh Agarwal, Dean~P Foster, Daniel~J Hsu, Sham~M Kakade, and Alexander Rakhlin.
\newblock Stochastic convex optimization with bandit feedback.
\newblock \emph{SIAM Journal on Optimization}, 23\penalty0 (1):\penalty0 213--240, 2013.

\bibitem[Bhatia(1997)]{bhatia97}
Rajendra Bhatia.
\newblock \emph{Matrix Analysis}, volume 169.
\newblock Springer, 1997.

\bibitem[Bubeck and Eldan(2018)]{BE18}
S{\'e}bastien Bubeck and Ronen Eldan.
\newblock Exploratory distributions for convex functions.
\newblock \emph{Mathematical Statistics and Learning}, 1\penalty0 (1):\penalty0 73--100, 2018.

\bibitem[Bubeck et~al.(2015)Bubeck, Dekel, Koren, and Peres]{BDKP15}
S{\'e}bastien Bubeck, Ofer Dekel, Tomer Koren, and Yuval Peres.
\newblock Bandit convex optimization: $\sqrt{T}$ regret in one dimension.
\newblock In \emph{Conference on Learning Theory}, pages 266--278, 2015.

\bibitem[Bubeck et~al.(2021)Bubeck, Eldan, and Lee]{bubeck2021kernel}
S{\'e}bastien Bubeck, Ronen Eldan, and Yin~Tat Lee.
\newblock Kernel-based methods for bandit convex optimization.
\newblock \emph{Journal of the ACM (JACM)}, 68\penalty0 (4):\penalty0 1--35, 2021.

\bibitem[Demaret et~al.(2019)]{demaret2019stein}
Tom Demaret et~al.
\newblock About stein's estimators: the original result and extensions.
\newblock 2019.

\bibitem[Drusvyatskiy and Lewis(2018)]{drusvyatskiy2018error}
Dmitriy Drusvyatskiy and Adrian~S Lewis.
\newblock Error bounds, quadratic growth, and linear convergence of proximal methods.
\newblock \emph{Mathematics of Operations Research}, 43\penalty0 (3):\penalty0 919--948, 2018.

\bibitem[Flaxman et~al.(2005)Flaxman, Kalai, and McMahan]{FK05}
Abraham~D Flaxman, Adam~Tauman Kalai, and H~Brendan McMahan.
\newblock Online convex optimization in the bandit setting: Gradient descent without a gradient.
\newblock In \emph{SODA'05: Proceedings of the sixteenth annual ACM-SIAM symposium on Discrete algorithms}, pages 385--394, 2005.

\bibitem[Fokkema et~al.(2024)Fokkema, van~der Hoeven, Lattimore, and Mayo]{fokkema2024onlinenewtonmethodbandit}
Hidde Fokkema, Dirk van~der Hoeven, Tor Lattimore, and Jack~J. Mayo.
\newblock Online newton method for bandit convex optimisation, 2024.
\newblock URL \url{https://arxiv.org/abs/2406.06506}.

\bibitem[Hazan and Levy(2014)]{HL14}
Elad Hazan and Kfir Levy.
\newblock Bandit convex optimization: Towards tight bounds.
\newblock In \emph{Advances in Neural Information Processing Systems}, pages 784--792, 2014.

\bibitem[Hazan et~al.(2007)Hazan, Agarwal, and Kale]{hazan2007logarithmic}
Elad Hazan, Amit Agarwal, and Satyen Kale.
\newblock Logarithmic regret algorithms for online convex optimization.
\newblock \emph{Machine Learning}, 69\penalty0 (2):\penalty0 169--192, 2007.

\bibitem[Hodgkinson and Mahoney(2020)]{hodgkinson2020multiplicativenoiseheavytails}
Liam Hodgkinson and Michael~W. Mahoney.
\newblock Multiplicative noise and heavy tails in stochastic optimization, 2020.
\newblock URL \url{https://arxiv.org/abs/2006.06293}.

\bibitem[Karimi et~al.(2020)Karimi, Nutini, and Schmidt]{karimi2020linearconvergencegradientproximalgradient}
Hamed Karimi, Julie Nutini, and Mark Schmidt.
\newblock Linear convergence of gradient and proximal-gradient methods under the polyak-\l{}ojasiewicz condition, 2020.
\newblock URL \url{https://arxiv.org/abs/1608.04636}.

\bibitem[Kleinberg(2005)]{Kle04}
Robert Kleinberg.
\newblock Nearly tight bounds for the continuum-armed bandit problem.
\newblock In \emph{Advances in Neural Information Processing Systems}, pages 697--704. MIT Press, 2005.

\bibitem[Lattimore(2020)]{Lat20-cvx}
Tor Lattimore.
\newblock Improved regret for zeroth-order adversarial bandit convex optimisation.
\newblock \emph{Mathematical Statistics and Learning}, 2\penalty0 (3/4):\penalty0 311--334, 2020.

\bibitem[Lattimore(2024)]{lattimore40bandit}
Tor Lattimore.
\newblock Bandit convex optimisation, 2024.
\newblock URL \url{https://arxiv.org/abs/2402.06535v3}.

\bibitem[Lattimore and Gy\"orgy(2021)]{LG21a}
Tor Lattimore and Andr\'as Gy\"orgy.
\newblock Improved regret for zeroth-order stochastic convex bandits.
\newblock In \emph{Conference on Learning Theory}, pages 2938--2964, 2021.

\bibitem[Lattimore and György(2023)]{lattimore2023secondordermethodstochasticbandit}
Tor Lattimore and András György.
\newblock A second-order method for stochastic bandit convex optimisation, 2023.
\newblock URL \url{https://arxiv.org/abs/2302.05371}.

\bibitem[Laurent and Massart(2000)]{10.1214/aos/1015957395}
Beatrice Laurent and Pascal Massart.
\newblock {Adaptive estimation of a quadratic functional by model selection}.
\newblock \emph{The Annals of Statistics}, 28\penalty0 (5):\penalty0 1302 -- 1338, 2000.
\newblock \doi{10.1214/aos/1015957395}.
\newblock URL \url{https://doi.org/10.1214/aos/1015957395}.

\bibitem[Lumbreras and Tomamichel(2024)]{lumbreras2024linear}
Josep Lumbreras and Marco Tomamichel.
\newblock Linear bandits with polylogarithmic minimax regret.
\newblock \emph{arXiv preprint arXiv:2402.12042}, 2024.

\bibitem[Lumbreras et~al.(2022)Lumbreras, Haapasalo, and Tomamichel]{Lumbreras_2022}
Josep Lumbreras, Erkka Haapasalo, and Marco Tomamichel.
\newblock Multi-armed quantum bandits: Exploration versus exploitation when learning properties of quantum states.
\newblock \emph{Quantum}, 6:\penalty0 749, June 2022.
\newblock ISSN 2521-327X.
\newblock \doi{10.22331/q-2022-06-29-749}.
\newblock URL \url{http://dx.doi.org/10.22331/q-2022-06-29-749}.

\bibitem[Mhammedi(2024)]{mhammedi2024online}
Zakaria Mhammedi.
\newblock Online convex optimization with a separation oracle.
\newblock \emph{arXiv preprint arXiv:2410.02476}, 2024.

\bibitem[Stein(1981)]{stein1981estimation}
Charles~M Stein.
\newblock Estimation of the mean of a multivariate normal distribution.
\newblock \emph{The annals of Statistics}, pages 1135--1151, 1981.

\bibitem[Suggala et~al.(2024)Suggala, Sun, Netrapalli, and Hazan]{suggala2024second}
Arun Suggala, Y~Jennifer Sun, Praneeth Netrapalli, and Elad Hazan.
\newblock Second order methods for bandit optimization and control.
\newblock \emph{arXiv preprint arXiv:2402.08929}, 2024.

\bibitem[Vershynin(2018)]{vershynin2018high}
Roman Vershynin.
\newblock \emph{High-dimensional probability: An introduction with applications in data science}, volume~47.
\newblock Cambridge university press, 2018.

\end{thebibliography}
